\documentclass[11pt,a4paper]{ansarticle}      
\usepackage{setspace,graphicx}     
\usepackage{enumerate}  

\usepackage{amssymb}
\usepackage{mystyle}    
\usepackage{mythesisStyle}   
\usepackage{setspace}
\usepackage[numbers]{natbib}
\usepackage{subfigure}
\usepackage[toc,page]{appendix}
\DeclareGraphicsExtensions{.pdf,.png}

\title{Formulation, Implementation and Validation of the Horizontal Coupling Method for 1D/2D Shallow Water Flow Models}
\author{C. Nwiagwe \thanks{nwaigwe@warwick.ac.uk}}
\author{A. S. Dedner \thanks{a.s.dedner@warwick.ac.uk}}
\affil{Centre for Scientific Computing and Warwick Mathematics Institute, University of Warwick, United Kingdom}

\begin{document}
\maketitle

%
One dimensional (1D) simulations of the flow and flooding of open channels are known to be inaccurate
as the flow is multi-dimensional in nature, especially at the flooded regions. However, multi-dimensional
simulations, even in two dimensions (2D), are computationally expensive, hence the problem of efficiently
coupling 2D and 1D simulations for the flow and flooding of open channels has been the subject of much
research and is investigated in this paper.
We adopt a 1D model with coupling term for the channel flow and the 2D shallow water flow model for the floodplain.
The 1D model with coupling term is derived by integrating the 3D Free Surface Euler equations 
but without imposing any restriction on the channel width variations. 
Finite volume methods are formulated for both the 2D and 1D models including a discrete coupling term in closed form.
Coupling is achieved through the discrete coupling term in the 1D model and the lateral numerical fluxes in the 2D model.
Since the lateral discharge in the channel cannot be guaranteed to be zero
during flooding, we aim to recover the lateral variation by computing two lateral discharges over 
each cross section and propose to use an ad-hoc model based on the $y$-discharge equation in the 2D model for this 
purpose. We then propose the numerical scheme for this ad-hoc model following the hydrostatic reconstruction
philosophy. 
Then, we show that the resulting method, named Horizontal Coupling Method (HCM), is well-balanced;
we introduce the no-numerical flooding property and also show that the method satisfies the property.
Three numerical test cases are used to verify the performance of the method.
The results show that the method performs well in both accuracy and efficiency and also approximates
the channel lateral discharges with very good accuracy and little computational overhead.

\section{Introduction}
Flows in open channels, such as rivers, in which the vertical and lateral variations 
in velocities can be assumed negligible, can be accurately simulated using the 1D Saint Venant Equations. 
During flooding, the channel overflows and the flow becomes high dimensional, rendering
1D simulation inadequate. These claims have also been observed numerically, see \cite{chineduthesis}
for example. But then, even a 2D simulation of the entire flow is computationally
expensive. This leads to the difficulty of choosing between an expensive but more accurate
high dimensional simulations and an inexpensive but less accurate 1D simulation. To tackle this problem,
a 1D simulation can be used along the channel while a 2D simulation is used for the floodplains.
The problem of how to couple the two simulations then arises.
This has led to many research work and also the subject of this paper.

A lot of research has been carried out to propose methods to couple 1D channel model with 2D floodplain
flow model.
In \cite{bladeetal2012}, a method, which numerically couples the 1D and 2D models by including the
lateral numerical fluxes in the 1D numerical scheme, is proposed. 
They referred to the method as the Flux-Based Method (FBM).
The theory of characteristics was employed to couple  1D/2D models in \cite{yongcanetal2012}; 
matching conditions are defined at the 2D/1D interfaces, then a prediction and correction algorithm
was used to ensure that these conditions are satisfied. 
The 1D river model and the 2D non-inertia model were also coupled in \cite{solomonetal2012} 
to simulate the interaction of a sewer system with over-land flow. 
Here, the water level differences between the flows in the two domains are used to calculate the 
interacting discharges in the sub-domains.

In \cite{moralesetal2013}, see also \cite{moralesthesis}, two methods which are based on 
post-processing of separately computed solutions of the existing 1D and 2D models, are proposed.
In these methods, the separately computed solutions are used to calculate the total water volume
in a 1D cell and all its adjacent floodplain 2D cells. Then, a common water level is found for all the cells, finally
the wetted cross sectional area for the 1D cell and the water height for all the adjacent 2D cells are found. 
These methods have been applied to Tiber River, Rome in \cite{moralesapplied2016}.
The superposition approach, proposed in \cite{monnierMarin2009super}, classically derives the exchange
terms in the 1D model from the full 3D Inviscid Euler's equations; an optimal control 
process is used to couple the models.
The superposition approach of \cite{monnierMarin2009super}
was extended to finite volume methods in \cite{nietocoupled}, proposing a discrete exchange term that leads to globally
well-balanced scheme. This approach superposes a 2D grid over the 1D channel
grid and convergence is achieved using a Schwartz-like iterative algorithm. For practical cases,
the iterative algorithm can jeopardise the overall efficiency
of the method \cite{nicoleetal20142d}.

A great difficulty for coupling methods is how to calculate the lateral discharge 
along the river channel because the 1D model does not have an equation to compute it.
The channel lateral discharge is set to zero and used to calculate the 2D numerical 
fluxes at the 2D/1D interfaces in \cite{comparison2015}.
In \cite{nicoleetal20142d},  the exchange terms derived in \cite{monnierMarin2009super} were adopted and
a strategy to estimate the lateral discharge without superposition or overlapping, was proposed.
The approach is an iterative technique which uses the solution of successive Riemann problems 
to estimate the transverse velocity.  Another approach which decomposes the  channel 1D discharge
into lateral and frontal components, using the angle which the channel axis makes with the $x$-axis,
can be found in \cite{bladeetal2012,moralesetal2013}. However, if this angle
is zero, then this approach would be inadequate whenever the channel is full because it would always compute
a zero discharge which is unrealistic. Therefore, the problem of computing the channel lateral or
transverse discharge remains challenging.

In addition to the above difficulty, another fundamental 
issue is the 1D assumption on the channel flow, namely that both the
free surface elevation and lateral velocity are laterally constant.
By physical intuition, during overflow  like flooding or draining, water flows out of or into
the channel from both of its lateral boundaries. 
This means that the lateral velocities (or discharges) at both sides are in opposite directions and
very likely to differ in magnitude.
Therefore, the lateral discharge will rarely be constant across the channel cross sections,
even when the free-surface elevation is assumed constant over the cross section. 
This means that the 1D assumption is inadequate if overflowing.
However, most existing coupling methods retain this 1D assumption even during overflowing.
We, therefore, propose that different discharges for each lateral boundary of
a cross section, need to be computed; we propose to use the 2D $y$-component shallow water equation,
as an ad-hoc model, to compute these lateral discharges.

The method we propose here, which we the Horizontal Coupling Method (HCM), follows the lines of \cite{monnierMarin2009super}
to derive a similar but 
slightly different coupling terms however, we do not impose or use any restriction on the channel width
variations.
The essence of this paper is, therefore, to propose a strategy (i) to overcome the difficulty in calculating
the lateral discharges, (ii) eliminate  the limitations of the 1D assumption on the channel lateral discharge during flooding,
(iii) to derive a more general variant of the coupling terms of \cite{monnierMarin2009super} (iv) prove the properties
of the resulting method and (v) validate the method using some numerical test cases.

The paper is organised as follows. In section \ref{hsecmodels}, we present the background 
for the problem of the flow and flooding of open channels and derive the channel flow model with  coupling term. 
We also present the 2D shallow water flow models for the floodplain flows.
The numerical schemes for the uncoupled 1D  and 2D  models are presented in section  \ref{hsecschemes};
the numerical scheme for the lateral discharge model and the detailed formulation of the discrete coupling term
are given in sections \ref{horsecqyscheme} and \ref{horseccopscheme} respectively. 
The algorithm for the HCM is summarised in section \ref{summary-hcm}; a flow chart for the implementation
is also given.
The properties of the method are discussed in section \ref{horiprops} where we introduce 
the no-numerical flooding property and show that the method preserves the property and is
also well-balanced. 
In section \ref{hsecnumresults}, we present some numerical test cases to access the performance of the method and
discuss the results; we conclude the paper in section \ref{hsecconc}.

\section{Mathematical Models}\label{hsecmodels}
\subsection{Background}

\begin{figure}[ht!]
 \includegraphics[width=0.98\textwidth, height=70mm]{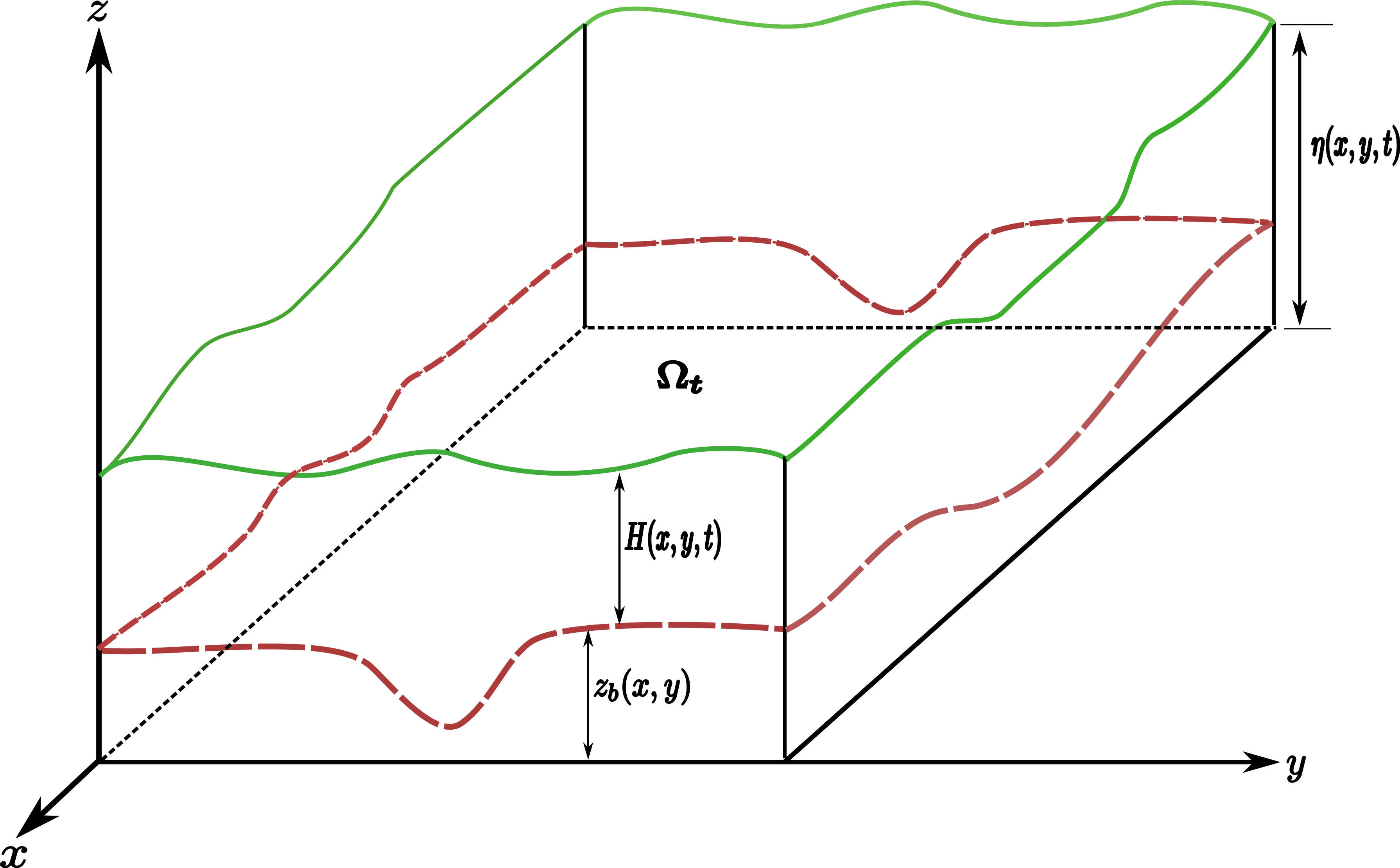}
 \caption{
 	Flow over a domain with bottom topography $z_b(x,y)$ (dashed line) comprising of a channel and floodplain.
    The channel length is along the $x$-axis and the width, along the $y$-axis; $H(x,y,t)$ is water depth,
    $\eta(x,y,t)$, the free surface elevation, $t$ is time variable and $(x,y,z) \in \real^3$.}
\label{fig3dflowdomain}
\end{figure}

\begin{figure}[ht!]
 \subfigure[Flow cross section at a fixed point $x$. The bottom topography comprises of the channel with bank elevations, $z_{bl}(x)$ and $z_{br}(x)$, and  
  		      the floodplain which occupies the remaining regions.]
          {\label{figFlowCrossSection}\includegraphics[width=0.45\textwidth, height=0.35\textwidth]{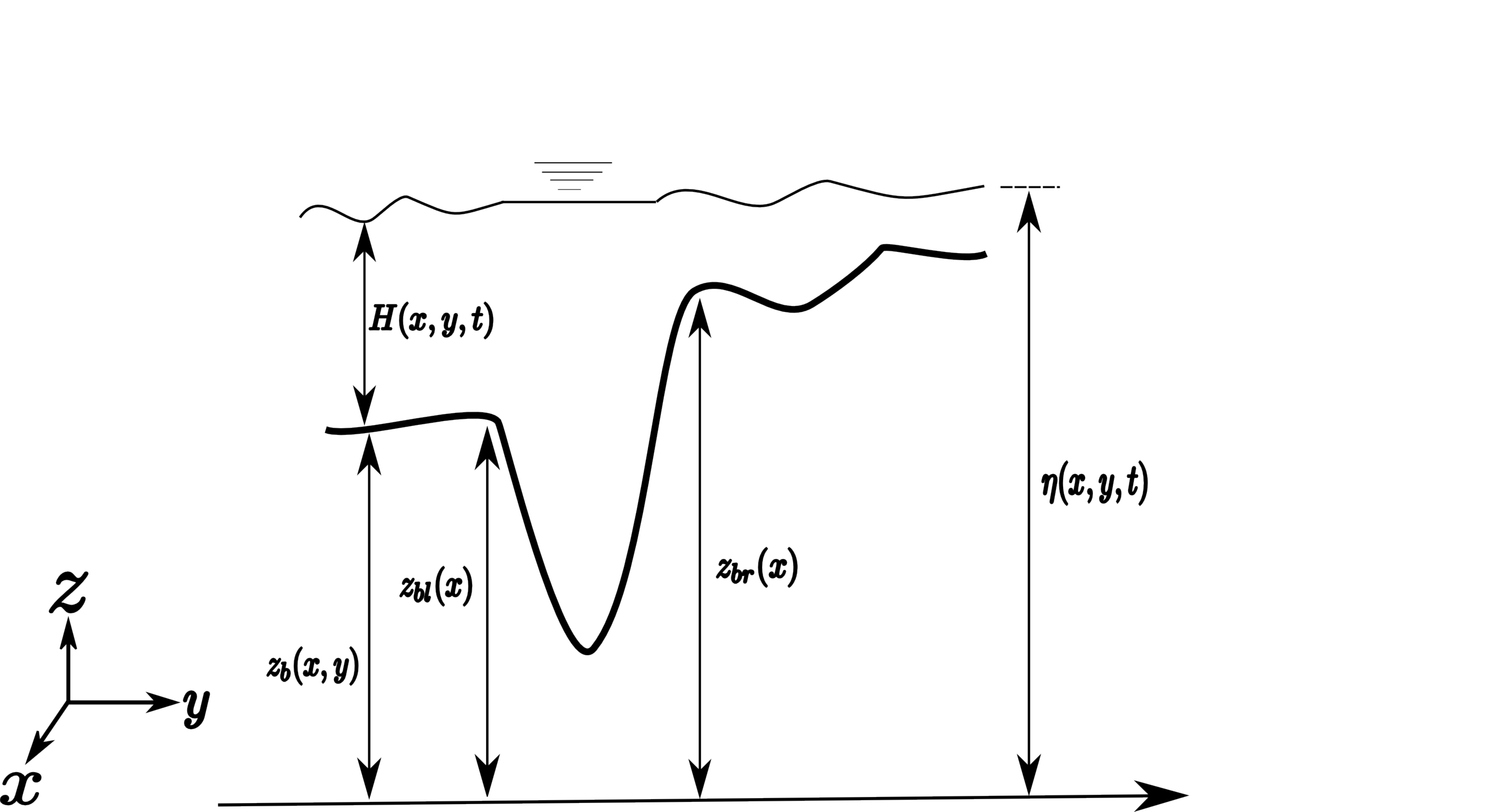}}
\hspace{0.03\textwidth}
 \subfigure[Channel cross section, showing the channel wall elevation $\zwall$, the top width $B(x,\zwall)$, the bottom elevation in 1D sense $Z_b(x)$,  
           laterally flat free-surface elevation $\etabar$, 
           and the $y-$coordinates  $\ylwall := y_l(x,\zwall)$ and $\yrwall :=y_r(x,\zwall)$ respectively of the left and right lateral walls at 
           the channel top.]
          {\label{figChannelCrossSection}    \includegraphics[width=0.5\textwidth, height=0.4\textwidth]{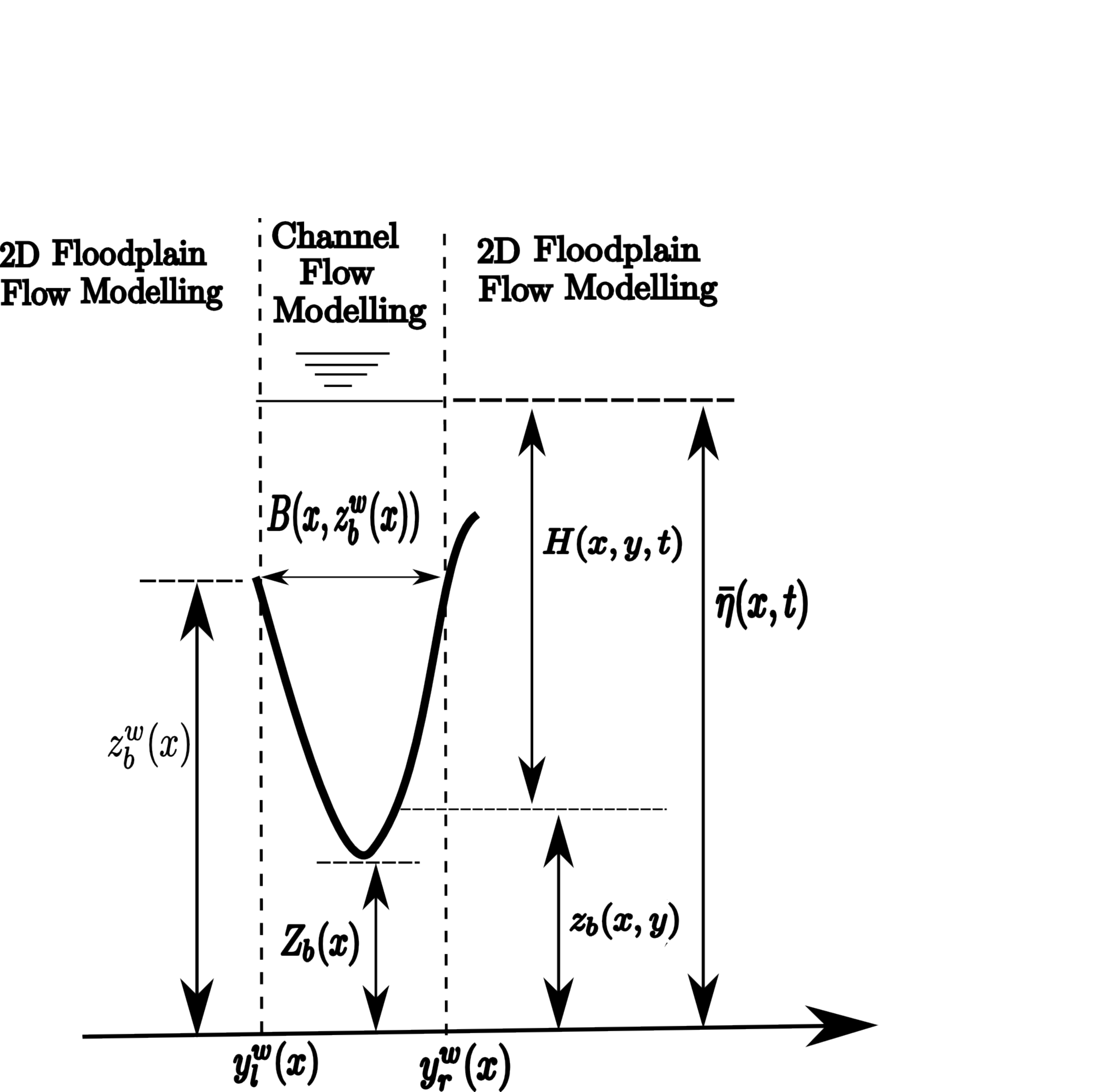}}
\caption{Flow cross sections}
\end{figure}

In this section, we present the model equations for the problem under consideration.
Let us begin by considering the flow of water over the fixed horizontal 2D domain, $\Omega_H \subset \real^2$ which consists
of a channel and floodplains (see figure \ref{fig3dflowdomain}), such that 
the flow at time, $t$ occupies the 3D domain, 
$\Omega_t$ defined by
\begin{equation} \label{moddomain3D}
	\Omega_t 
     	   =\{
			(\xy,z) \in \real^3 : \xy =(x,y) \in \Omega_H, \, 
			  z_b(\xy) \leq z \leq \eTa		
		   \},
\end{equation}
bounded below by a fixed bottom, $z_b(\xy)$ and above by the water free-surface elevation, $\eTa$ given by
\begin{equation}
	 \eTa = z_b(\xy) + H(\xy,t)
\end{equation}
and $H(\xy,t)$ is the depth of fluid at time, $t$. A cross section of the flow domain is shown
in figure \ref{figFlowCrossSection}.
The length of the channel lies along the $x$-axis (frontal direction) and the width, along the $y$-axis (lateral direction), while
$z_{b,l}(x)$ and $z_{br}(x)$ are the left and right bank elevation of the channel, see figure \ref{figFlowCrossSection}.
An important quantity is the maximum channel wall elevation or simply, channel wall elevation.
\begin{mydefinition}[Channel wall elevation, $z_b^w(x)$]\label{modzwalldef}
	The channel wall elevation at cross section x, is the minimum elevation of the channel banks
	above which flooding is said to have occurred. We denote it by $z_b^w(x)$, that is
	\begin{equation}
		z_b^w(x)  = \min( z_{bl}(x), z_{br}(x) ),
	\end{equation}
see figures \ref{figFlowCrossSection} and \ref{figChannelCrossSection}.
\end{mydefinition}
That is, $z_b^w(x)$ is the channel top. 
Figure \ref{figChannelCrossSection} shows the channel cross section, depicting its geometry, including
the 1D laterally constant channel bottom topography given by  
	$$Z_b(x) = \min_{y}z_b(x,y).$$
It also depicts the wall elevation $\zwall$, the top width $B(x, \zwall)$ and the $y$-coordinates, 
$\ylwall:= y_l(x,\zwall)$ and $\yrwall := y_r(x,\zwall)$  of the left and right
boundaries at the top, where  $B(x,z)$ gives the channel width at an elevation $z$ above a reference elevation $z$,
and $y_l(x,z)$ and $y_r(x,z)$ are the  $y-$coordinates of the left and right lateral boundaries, respectively, at elevation $z$.
So that
	\begin{equation}\label{modwidthequation}
	  B(x,z) = \yrz - \ylz \quad \forall z,
   \end{equation}
such that
\begin{equation}\label{hor-eqn-bank-1a}
 \yrz =\ylz,\quad	B(x,z) = 0  \quad  \mbox{ for all } z < Z_b(x);
\end{equation}
and the bottom elevation satisfies
\begin{align}\label{hor-eqn-bank-1b}
 	\zb|_{y=\ylz,\yrz} = z  \quad  \forall z \in [ Z_b(x), z_b^w(x)   ],
\end{align} 
see figure \ref{figChannelCrossSection}.
Furthermore, we extend the definition of the width functions above the channel top ($z>\zwall$) 
in the following:
\begin{align}\label{hor-eqn-bank-1c}
 	y_{l,r}(x, z) = y_{l,r}^w(x), \quad B(x,z) = B(x,\zwall) \quad \forall z \ge \zwall,
\end{align}
see figure \ref{figChannelCrossSection}.

The flow cross section, $-\infty < y < \infty$ in \ref{figFlowCrossSection} has been partitioned into
(i) the channel cross section, $\ylwall \le y \le \yrwall$
and (ii) the floodplains, $-\infty < y \le \ylwall $ and $\yrwall \le y < -\infty$, see figure \ref{figChannelCrossSection}. 
The flow in the floodplains is simulated with the standard 2D shallow water models (see section \ref{horsec2dswes}), therefore we
focus on deriving the model equations 
for the flow in the channel.

With the channel geometry completely defined, we now consider the initial flow condition. In general, the free surface elevation $\eTa$ is 2D, see 
figures \ref{fig3dflowdomain} and \ref{figFlowCrossSection}. However, we assume that it is always 1D (laterally constant and given by
 $\etabar$) within the channel, see figure \ref{figChannelCrossSection}). Hence, we say the channel is full whenever
\begin{equation}
  \etabar > \zwall.
\end{equation}
Note that in general, the channel flow lateral boundaries are at the coordinates, $y_{l,r}(x, \etabar)$, not $y_{l,r}^w(x)$,
see figures \ref{figChannelCrossSection} and \ref{fig-non-full-channel}. These two sets of coordinates are only equal
if the channel is full, see \eqref{hor-eqn-bank-1c} and also figure \ref{figChannelCrossSection}. They are not equal
if the channel is not full, see figure \ref{fig-non-full-channel}.
If the channel is not full ($\etabar \le \zwall$), then the water height and velocities are zero at
the top lateral boundaries, that is
\begin{align}\label{svmbanks}
    \bigg( u(\xyzt), v(\xyzt), w(\xyzt), H(\xyt) \bigg )\bigg|_{y=y_{l,r}^w(x)}= 0 \quad \mbox{ whenever } \etabar \le \zwall,
\end{align}
(see figure \ref{fig-non-full-channel}) where $u,v,w$ are velocity components along the $x,y,z$ directions respectively.

\subsection{Derivation of the channel flow model with coupling term}\label{sec-derive-models}
Under the assumption of compressible and inviscid fluid, the flow of water in the channel is governed by the following 3D Free-Surface Euler Equations
\cite{lannes2013}:
%
\begin{align}
     \begin{split} \label{fseecont}
		\pdiff{x}{u(\xyzt)} + \pdiff{y}{v(\xyzt)} + \pdiff{z}{w(\xyzt)} = 0.
	 \end{split}
	    \\
	\begin{split}\label{fseexmom}
		 \pdiff{t}{u(\xyzt)} 
		+ u(\xyzt)\pdiff{x}{u(\xyzt)} 
		+ v(\xyzt)\pdiff{y}{u(\xyzt)}
		+ w(\xyzt)\pdiff{z}{u(\xyzt)}
		= -\frac{1}{\rho} \pdiff{x}{P(\xyzt)}.
	\end{split}
	\\
	\begin{split}\label{fseeymom}
		 \pdiff{t}{v(\xyzt)} 
		+ u(\xyzt)\pdiff{x}{v(\xyzt)} 
		+ v(\xyzt)\pdiff{y}{v(\xyzt)}
		+ w(\xyzt)\pdiff{z}{v(\xyzt)}
		= -\frac{1}{\rho} \pdiff{y}{P(\xyzt)}.
	\end{split}
	\\
	\begin{split}\label{fseezmom}
		 \pdiff{t}{w(\xyzt)} 
		+ u(\xyzt)\pdiff{x}{w(\xyzt)} 
		+ v(\xyzt)\pdiff{y}{w(\xyzt)}
		+ w(\xyzt)\pdiff{z}{w(\xyzt)}
		\\= 
		-\frac{1}{\rho} \pdiff{x}{P(\xyzt)} - g.
	\end{split}
\\
\begin{split}
	\left( 
		   u(\xyzt) \pdiff{x}{z_b(\xy)} 	+ v(\xyzt) \pdiff{y}{z_b(\xy)}	- 	   
		 w(\xyzt) 
	\right)\bigg|_{z=z_b(\xy)}	
	= 0.  \label{modkinbed}
\end{split}
\\
\begin{split}
	\left( 
		\pdiff{t}{\eTa} + u(\xyzt) \pdiff{x}{\eTa} 	+ v(\xyzt) \pdiff{y}{\eTa}	- w(\xyzt)
	\right)\bigg|_{z=\eTa}	 
	= 0.   \label{modkinsurf}
\end{split}
\end{align}
%
\begin{equation}\label{moddyn}
	P(\xyzt) = P_{atm}  \quad  \mbox{ on } z = \eTa,
\end{equation}
where $\rho, (u,v,w)^T$ and $P$ are the fluid density, velocity vector and pressure at point $(\xyz)$ at time, $t$
and $P_{atm}$ is the atmospheric pressure, which is usually conveniently taken to be zero.

The flow quantities of interest in the 1D channel model  are the wetted cross sectional area, $A(x,t)$ and the section 
averaged discharge, $Q(x,t)$ given by the following averages:
\begin{align}
	Q(\xt) &= \dint{z}{y}{\yle}{\yre}{\zb}{\etaxt}{ u(\xyzt)}, \label{modQ}
	\\
	A(\xt) &= \dint{z}{y}{\yle}{\yre}{\zb}{\etaxt}{} = \sint{y}{\yle}{\yre}{H(\xyt)} \label{modA}.
\end{align}
So that the section-averaged velocity, $\uoned$ is given as
\begin{equation}\label{moduxt}
	\uoned (\xt)= \frac{Q(\xt)}{A(\xt)}
		         = \frac{1}{A(\xt)}\dint{z}{y}{\yle}{\yre}{\zb}{\etaxt}{ u(\xyzt)}.
\end{equation}
First, we note that $y$-independence of the free-surface, $\etaxt$ means that
the sum,
\begin{equation}\label{etahzb}
	H(\xyt) + \zb = \etaxt, \quad \forall  y_l(x,\etabar) \leq y \leq y_r(x,\etabar),
\end{equation}
is constant in $y$, even though each of $H(\xyt)$ and $\zb$ depends on $y$, see figures \ref{figChannelCrossSection}
and \ref{fig-non-full-channel}.

The shallow water assumption that water depth is small compared to horizontal length leads to neglect
vertical acceleration, hence the z-momentum equation, \eqref{fseezmom} and the dynamic boundary 
condition, \eqref{moddyn},  lead to the following hydrostatic pressure:
\begin{align*}
	P(\xyzt) = \rho g ( \etaxt - z ) \quad \Longrightarrow   \pdiff{x}{P(\xyzt)} = \rho g \pdifft{x}{\etaxt}.
\end{align*}
%
Therefore, the FSEE \eqref{fseecont} -\eqref{moddyn}, reduce to the following system:
\begin{align}
&	  \pdiff{x}{u(\xyzt)} + \pdiff{y}{v(\xyzt)} + \pdiff{z}{w(\xyzt)}  = 0. \label{svmcont}
	\\
&	    \pdiff{t}{u(\xyzt)} 
	  +  \pdiffb{x}{u^2(\xyzt)} 
	  +  \pdiffb{y}{u(\xyzt)v(\xyzt)} 
	  +  \pdiffb{z}{u(\xyzt)w(\xyzt)} 
       =  -g\pdiff{x}{\etaxt}.\label{svmxmom}
\\
&	\left( 
		   u(\xyzt) \pdifft{x}{z_b(\xy)} 	+ v(\xyzt) \pdifft{y}{z_b(\xy)}	-  w(\xyzt) 
	\right)\bigg|_{z=z_b(\xy)}	
	= 0.  \label{modsvmkinbed}
    \\
&	\left( 
		\pdifft{t}{\etaxt} + u(\xyzt) \pdifft{x}{\etaxt} 	- w(\xyzt)
	\right)\bigg|_{z=\etaxt}	 
	= 0.   \label{modsvmkinsurf}
\end{align}

Integrating \eqref{svmcont} vertically (over $\zb \leq z \leq \bar{\eta}$),
and laterally (over $\yle \leq y \leq \yre$), applying the Leibnitz rule and using the 
kinematic boundary conditions, \eqref{modsvmkinbed}, \eqref{modsvmkinsurf}, we have the following mass 
equation \eqref{horAImportant} below see \cite{monnierMarin2009super}. Repeating the 
same process for \eqref{svmxmom}, gives the discharge equation \eqref{horQImportant}
below:
\begin{align}
   \pdifft{t}{A(\xt)}		+   \pdifft{x}{ Q(\xt)   } &= \Phi^A(\xt)   \label{horAImportant},
   \\
   \pdifft{t}	{ Q(\xt) 	}		+ \pdiffb{x}{ \frac{Q^2(\xt)}{A(\xt)}	}	
	   &= -gA(\xt)\pdifft{x}{\etaxt} + \Phi^Q(\xt)   \label{horQImportant},
\end{align}
see  \cite{chineduthesis, monnierMarin2009super} for details;   $\Phi^A(\xt) $ and $\Phi^Q(\xt)$ are the coupling terms defined as
\begin{align}
\begin{split}\label{eqn-couplA1}
  \Phi^A(\xt) &= 
  	      \pdifft{x}{\yre}
  	      \left[ 
            \sint{z}{\zb}{\etaxt}{u(\xyzt)} 
         \right]_{ y=\yre }       
         -
         \left[ 
            \sint{z}{\zb}{\etaxt}{v(\xyzt)} 
         \right]_{ y=\yre }     
         \\   
         &
         -
         \pdifft{x}{\yle }
         \left[ 
            \sint{z}{\zb}{\etaxt}{u(\xyzt)} 
         \right]_{ y=\yle }
         +
         \left[ 
            \sint{z}{\zb}{\etaxt}{v(\xyzt)} 
         \right]_{ y=\yle }.         
\end{split}
\\
%
\begin{split}\label{eqn-couplQ1}
  \Phi^Q(\xt) &=  
         \pdifft{x}{\yre }
  		\left[ 
            \sint{z}{\zb}{\etaxt}{u^2(\xyzt)} 
         \right]_{ y=\yre }
         -
         \left[ 
            \sint{z}{\zb}{\etaxt}{u(\xyzt)v(\xyzt)} 
         \right]_{ y=\yre }
        \\
        &
        -
        \pdifft{x}{\yle }
         \left[ 
            \sint{z}{\zb}{\etaxt}{u^2(\xyzt)} 
         \right]_{ y=\yle }
         +
         \left[ 
            \sint{z}{\zb}{\etaxt}{u(\xyzt)v(\xyzt)} 
         \right]_{ y=\yle }. 
\end{split}
\end{align}

Note that if the channel is not full, then we have $\Phi^A(\xt) = 0 $ and $\Phi^Q(\xt) = 0$ because 
non-full channel means $\etabar \le \zwall \Longrightarrow z_b(x, y_{l,r}(x, \etabar)) = \etabar$ by \eqref{hor-eqn-bank-1b}.
Hence, both limits in all the integrals in \eqref{eqn-couplA1}-\eqref{eqn-couplQ1}, so the coupling terms vanish.
In this case, the model, \eqref{horAImportant}-\eqref{horQImportant}
reduces to the standard 1D Saint-Venant Models.

It is straight forward to show that 
$$\partial_x y_{l,r}(x, \etabar) \bigg[ \int_{\zb}^{\etabar} \theta(\xyzt) dz \bigg]_{y=y_{l,r}(x,\etabar)}
=
\partial_x y_{l,r}^w(x) \bigg[ \int_{\zb}^{\etabar} \theta(\xyzt) dz \bigg]_{y=y_{l,r}^w(x)}
$$
where $\theta$ is any of the integrands appearing in \eqref{eqn-couplA1}-\eqref{eqn-couplQ1}.
Therefore, we can conveniently write the coupling terms as follows:
\begin{align}
\begin{split}
  \Phi^A(\xt) &= 
  	      \pdifft{x}{\yrwall}
  	      \left[ 
            \sint{z}{\zb}{\etaxt}{u(\xyzt)} 
         \right]_{ y=\yrwall }       
         -
         \left[ 
            \sint{z}{\zb}{\etaxt}{v(\xyzt)} 
         \right]_{ y=\yrwall }     
         \\   
         &
         -
         \pdifft{x}{\ylwall }
         \left[ 
            \sint{z}{\zb}{\etaxt}{u(\xyzt)} 
         \right]_{ y=\ylwall }
         +
         \left[ 
            \sint{z}{\zb}{\etaxt}{v(\xyzt)} 
         \right]_{ y=\ylwall }.         
\end{split}
\\
%
\begin{split}
  \Phi^Q(\xt) &=  
         \pdifft{x}{\yrwall }
  		\left[ 
            \sint{z}{\zb}{\etaxt}{u^2(\xyzt)} 
         \right]_{ y=\yrwall }
         -
         \left[ 
            \sint{z}{\zb}{\etaxt}{u(\xyzt)v(\xyzt)} 
         \right]_{ y=\yrwall }
        \\
        &
        -
        \pdifft{x}{\ylwall }
         \left[ 
            \sint{z}{\zb}{\etaxt}{u^2(\xyzt)} 
         \right]_{ y=\ylwall }
         +
         \left[ 
            \sint{z}{\zb}{\etaxt}{u(\xyzt)v(\xyzt)} 
         \right]_{ y=\ylwall }. 
\end{split}
\end{align}

To proceed, let us define the following quantities:
\begin{align}\label{horqxqy}
\begin{split}
&	q_x(\xyt) = \sint{z}{\zb}{\eTa}{ u(\xyzt)   }, \quad 	q_y(\xyt) = \sint{z}{\zb}{\eTa}{ v(\xyzt)   }. 		
    \\
\mbox{So that }   \quad \quad &	 \frac{q^2_x}{H} \approx \sint{z}{\zb}{\eTa}{ u^2(\xyzt)   },   \quad
	 \frac{q_xq_y}{H} \approx \sint{z}{\zb}{\eTa}{ u(\xyzt) v(\xyzt)  }  \quad \mbox{(see \cite{chineduthesis, toroshock})}.
%
\end{split}
\end{align}
Using equations \eqref{horqxqy} the coupling terms become
\begin{align}
&	\Phi^A(\xt) =     q_x|_{y=y_r^w(x)} \pdifft{x}{y_r^w(x)} - q_y|_{y=y_r^w(x)}
		       - q_x|_{y=y_l^w(x)} \pdifft{x}{y_l^w(x)} + q_y|_{y=y_l^w(x)}. 
	\\	
&	\Phi^Q(\xt) =      \frac{q_x^2}{H}  \bigg|_{y=y_r^w(x)} \pdifft{x}{y_r^w(x)}
                - \frac{q_xq_y}{H} \bigg|_{y=y_r^w(x)}
				- \frac{q_x^2}{H}  \bigg|_{y=y_l^w(x)} \pdifft{x}{y_l^w(x)}
                + \frac{q_xq_y}{H} \bigg|_{y=y_l^w(x)}.                
\end{align}
%

\subsubsection*{Notational  Simplification}
We now express the coupling terms as functions of the fluxes at the channel lateral boundaries which are easier to compute.
\begin{figure}[ht!]
 \subfigure[Flow cross section for non-full channel.]
          {\label{fig-non-full-channel}\includegraphics[width=0.66\textwidth]{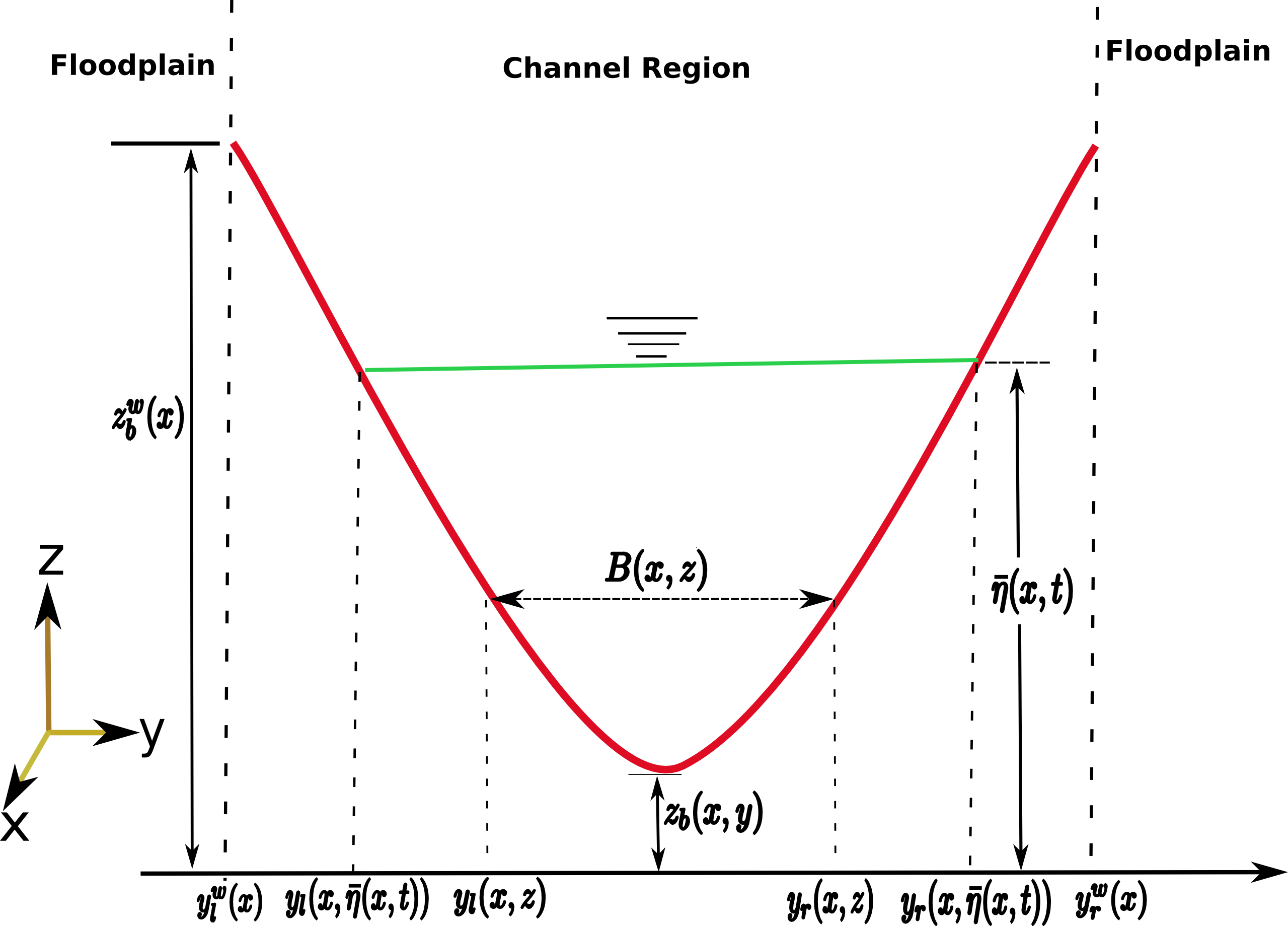}}
\hspace{0.03\textwidth}
 \subfigure[Top view of Lateral Boundaries (at elevation, $z=\zwall$]
          {\label{horfigLateralboundaries}    \includegraphics[width=0.3\textwidth]{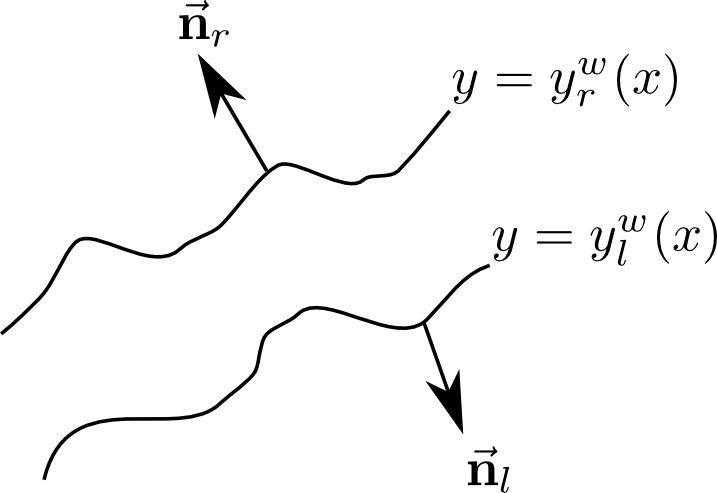}}
\caption{Non-full channel cross section (left) and channel top view (right).}
\end{figure}
Let $\myvec{n}_l= ( n_l^x, n_l^y )^T$ and $\myvec{n}_r= ( n_r^x, n_r^y )^T $
be the outward unit normal vectors to the lateral boundaries at $y=y^w_l(x)$ and $y=y_r^w(x)$
respectively (see figure \ref{horfigLateralboundaries}). Since these normal vectors are perpendicular
to the tangent lines to their respective lateral boundaries, we have
\begin{align*}
	\frac{n_l^y}{n_l^x} \pdifft{x}{y_l^w(x)}  = -1,
	\quad
	\frac{n_r^y}{n_r^x} \pdifft{x}{y_r^w(x)}  = -1
\quad \Longrightarrow \quad
	\pdifft{x}{y_l^w(x)}  = -\frac{n_l^x}{n_l^y},
	\quad
	\pdifft{x}{y_r^w(x)}  = -\frac{n_r^x}{n_r^y},  \quad n_l^y, n_r^y \neq 0.
\end{align*}
Define the vector $\myvec{q}=(q_x,q_y)^T$, then we write the coupling terms as 
\begin{equation}
	\Phi^A(x,t) = \Phi_L^A(x,t) + \Phi_R^A(x,t) \quad \mbox{ and } \quad
	\Phi^Q(x,t) = \Phi_L^Q(x,t) + \Phi_R^Q(x,t).
\end{equation}
where
\begin{align}
&  \Phi_L^A(x,t)  = \frac{1}{n_l^y}\bigg( \myvec{q}.\myvec{n}_l \bigg)\bigg|_{y=y_l^w(x)},   
  \quad
  \Phi_R^A(x,t)  =  -\frac{1}{n_r^y}\bigg( \myvec{q}.\myvec{n}_r \bigg)\bigg|_{y=y_r^w(x)}, \label{horcouplal}
  \\
&    \Phi_L^Q(x,t) = \frac{1}{n_l^y} \bigg( 
                     n_l^x \bigg [ \frac{q_x^2(\xyt)}{H(\xyt)}+ \frac{g}{2}H^2(\xyt) \bigg ]
                       + n_l^y\frac{q_y(\xyt)q_y(\xyt)}{H(\xyt)} 
                    \bigg ) \bigg|_{y=y_l^w(x)}
					- 
                  \frac{n_l^x}{n_l^y}\frac{g}{2}H^2(\xyt)\bigg|_{y=y_l^w(x)},\label{horcouplql}
  \\
 &   \Phi_R^Q(x,t) = -\frac{1}{n_r^y} \bigg( 
                     n_r^x \bigg [ \frac{q_x^2(\xyt)}{H(\xyt)}+ \frac{g}{2}H^2(\xyt) \bigg ]
                       + n_r^y\frac{q_y(\xyt)q_y(\xyt)}{H(\xyt)} 
                  \bigg ) \bigg|_{y=y_r^w(x)}
                  + 
                  \frac{n_r^x}{n_r^y}\frac{g}{2}H^2(\xyt)\bigg|_{y=y_r^w(x)}. \label{horcouplqr} 
\end{align}

Let $f^1_L(\xt)$ and $f^2_L(\xt)$ denote the first and second components, respectively, of the 
outgoing flux in the direction of $\myvec{n}_l$ at $y=\ylwall$, and $f^1_R(\xt)$ and
$f^2_R(\xt)$ be those in the direction of $\myvec{n}_r$, then
the coupling terms \eqref{horcouplal}-\eqref{horcouplqr}  can be written in the flowing  
forms:
\begin{align}
	&\Phi_L^A(\xt) = \frac{1}{n_l^y}f^1_L(\xt). \quad \Phi_R^A(\xt) = -\frac{1}{n_r^y}f^1_R(\xt).
	\\
	&\Phi_L^Q(\xt) = \frac{1}{n_l^y}f^2_L(\xt) - \frac{n_l^x}{n_l^y}\frac{g}{2} H^2(\xyt) \bigg|_{y=y_l^w(x)},
	\quad
	&\Phi_R^Q(\xt) = -\frac{1}{n_r^y}f^2_R(\xt) + \frac{n_r^x}{n_r^y}\frac{g}{2} H^2(\xyt) \bigg|_{y=y_r^w(x)}.
\end{align}

Hence, the 1D channel models with coupling term, in the presence of friction is
\begin{align}
\begin{split}\label{horAmodel}
 \pdifft{t}{A(\xt)}   +   \pdifft{x}{Q(\xt)}
  =& \frac{1}{n_l^y}f^1_L(\xt)  -\frac{1}{n_r^y}f^1_R(\xt),
 \end{split}
\\
\begin{split}\label{horQmodel}
 \pdifft{t}	{ Q(\xt) 	}	+ \pdiffb{x}{ \frac{Q^2(\xt)}{A(\xt)}	}	
	    =&   
      -    gA(\xt)\pdifft{x}{\etaxt}  +
     \frac{1}{n_l^y}f^2_L(\xt) 
      - \frac{n_l^x}{n_l^y}\frac{g}{2} H^2(\xyt) \bigg|_{y=y_l^w(x)}
      \\&
       -\frac{1}{n_r^y}f^2_R(\xt) 
	   + \frac{n_r^x}{n_r^y}\frac{g}{2} H^2(\xyt) \bigg|_{y=y_r^w(x)} + g A(\xt) S_f.
\end{split}              
\end{align}
where  $S_f = \frac{Q|Q|}{K^2}$ is the channel friction slope, 
$K = \frac{A^{k_1}}{nP^{k_2}}$ is the conveyance, 
$P$ is the wetted perimeter of channel cross-section,
$k_1=5/3, k_2=2/3$ and $n$ is the Manning coefficient, see \cite{cungeetal80, macdonaldthesis}.

\begin{myremark}
We obtained the above coupling terms without using or imposing any restriction on the channel width variation
as done in \citep{monnierMarin2009super}. And our coupling term clearly differs from theirs.
\end{myremark}

\subsubsection{Channel Flow Lateral Discharge Model}
To compute the lateral discharges in the channel,  we use
 the following $y$-discharge equation in the 2D Shallow water equations :
\begin{equation}\label{horqymodel}
\begin{split}
  \pdiff{t}{q_y(\xyt)} + \pdifft{x}{f_x(\Pi)} + \pdifft{y}{f_y(\Pi)} = -gH(\xyt)\pdifft{y}{\zb},\\
   \Pi = (H, q_x, q_y)^T,   \quad  f_x(\Pi) = \frac{q_xq_y}{H}, \quad 
   f_y( \Pi ) = \frac{q_y^2}{H} + \frac{1}{2}gH^2.
\end{split}
\end{equation}

\subsection{Floodplain Flow Model}\label{horsec2dswes}
We describe the flow in the floodplains using the 2D Shallow water equations, namely
\begin{align}\label{fv2dswewtfriceqnmodel}
	\partial_t \Pi + \nabla \cdot F(\Pi) =  S(\Pi, z_b ) + S_b(\Pi),
\end{align}
where
\begin{align}
\begin{split}
&	\Pi = \begin{pmatrix}
		H \\ q_x \\ q_y
	\end{pmatrix}, \quad
   F(\Pi) = \bigg( F_1(\Pi), F_2(\Pi)  \bigg),\quad
	   F_1(\Pi) = \begin{pmatrix}
   		q_x  \\ \frac{q_x^2}{H} + \frac{1}{2}gH^2 \\  \frac{q_xq_y}{H}
   \end{pmatrix}, \quad
   F_2(\Pi) = \begin{pmatrix}
   		q_y \\  \frac{q_xq_y}{H} \\ \frac{q_y^2}{H} + \frac{1}{2}gH^2
   \end{pmatrix},
\\
&	
   S_b(\Pi) = \begin{pmatrix}
   		                    0  \\  -g\frac{n^2}{H^{7/3}} \vec{q}|\vec{q}|
            \end{pmatrix}, \quad
  S(\Pi,z_b) = \begin{pmatrix}
   		                    0  \\  -gH\partial_x z_b(\xy) \\  -gh\partial_y z_b(\xy)
            \end{pmatrix},  
\end{split}                       
\end{align}
where $n$ is the manning coefficient, $S_b$ is the friction term and $S(\Pi,z_b)$
is the source term due to bottom topography term.
%

\section{Numerical Schemes}\label{hsecschemes}
\begin{figure}[ht!]
 \subfigure[Grid of the entire domain consisting of the 1D grid $\Omega_h^{1D}$ at the middle and the 2D grids $\Omega_h^{2D}$ for the floodplains.]
          {\label{fig2d1dgrid}\includegraphics[width=0.5\textwidth, height=0.35\textwidth]{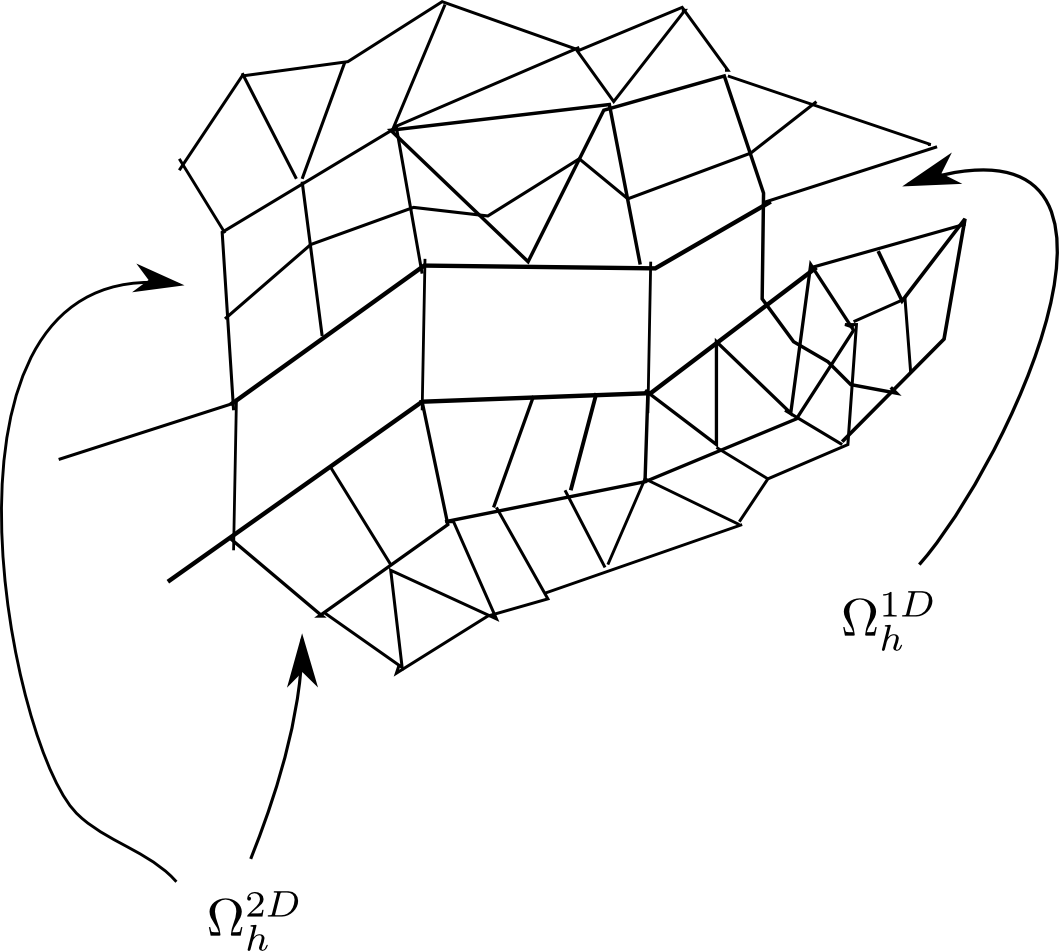}}
\hspace{0.02\textwidth}
 \subfigure[2D mesh showing two neighbour cells, $T_j$ and $T_k$, the edge $e_{jk}$ between them and the normal vector $\vec{n}_{jk}$.]
          {\label{fvm-fig-single2dmesh}\includegraphics[width=0.4\textwidth, height=0.35\textwidth]{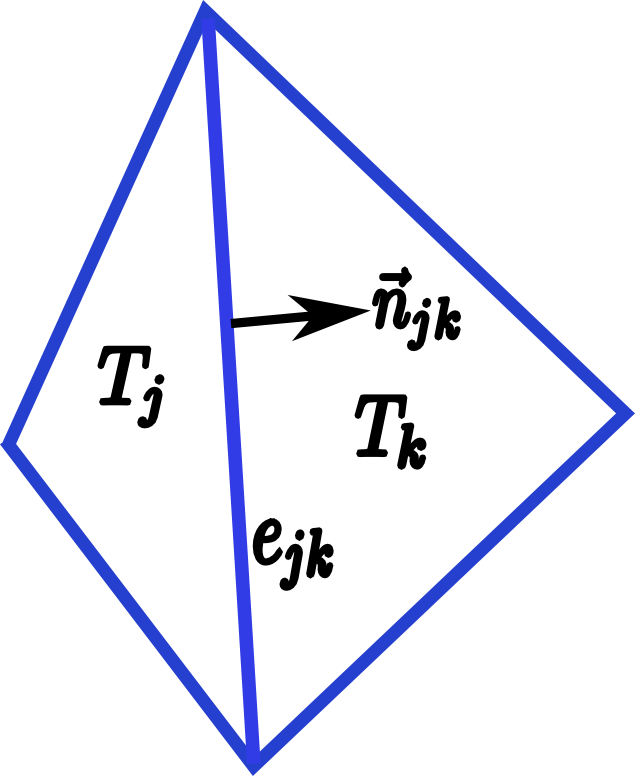}}
\caption{Grids}          
\end{figure}

In this section, we detail the numerical schemes for the models presented in previous sections.
To begin, we partition the channel into a 1D grid, 
$\Omega_h^{1D}$ made of cross sections and the floodplains, into a 2D grid $\Omega_h^{2D}$,
see figure \ref{fig2d1dgrid}. We first present the scheme for the 2D flood model, then
the schemes for the channel flow model is presented.

\subsection{Scheme For 2D Model}\label{hpsec2dscheme}
In this section, we present the scheme for the flood flow model, \eqref{fv2dswewtfriceqnmodel}.
Let $T_j\in \Omega_h^{2D}$ be an element of the 2D mesh $\Omega_h^{2D}$ in figure \ref{fig2d1dgrid}, 
and $T_k \in \Omega_h^{2D}$ be its neighbour cell, see figure \ref{fvm-fig-single2dmesh}.
Let $e_{jk}$ be the edge between $T_j$ and $T_k$, 
while $\vec{n}_{jk}$  is a unit vector normal to edge $e_{jk}$ and outward to $T_j$. 
Furthermore, let $|T_j|$ and $|e_{jk}|$ be the area of $T_j$ and length of $e_{jk}$ respectively and
let $\mathcal{E}_j$ be the set of all edges of $T_j$.  And let $\Pi_j^n=(H_j^n, q_{x,j}^n, q_{y,j}^n)^T$
be the approximate cell averages of the true solution in $T_j$, namely
\begin{align}
	\Pi_j^n = \frac{1}{|T_j|} \int_{T_j} \Pi(\xy,t^n) d\vec{x}.
\end{align}
Similarly, let $\Pi_k^n=(H_k^n, q_{x,k}^n, q_{y,k}^n)^T$ be the cell average vector $T_k$, while $z_{b,j}, z_{b,k}$
are the cell averages in $T_j, T_k $ respectively.

Then, we consider the following 2D hydrostatic reconstruction finite volume scheme \cite{audusseetal2004}:
\begin{align}\label{fvm2dswewbedfrictiongeneral}
\begin{split}
\Pi_j^{n+1} =& \Pi_j^n - \frac{\Delta t}{|T_j|}  \sum_{e_{jk}\in \mathcal{E}_j} |e_{jk}|  \bigg(
								T_{\vec{n}_{jk}}^{-1} \phi(  \widetilde{T_{\vec{n}_{jk}}\Pi_j^n}, \widetilde{T_{\vec{n}_{jk}}\Pi_k^n} )
							+  
								T_{\vec{n}_{jk}}^{-1}S^{hrm}( H_j^n, \tilde{H}_j^n) 
						\bigg)
\\						
  & + \Delta t S_b(\Pi_j^n),
\end{split}  
\end{align}	
where
\begin{align}
&  \tilde{H}_p^n := \max( H_p^n + z_{b,p} - \max(z_{b,j},z_{b,k})  ), \quad 
 \widetilde{T_{\vec{n}_{jk}} \Pi_p^n} := \frac{\tilde{H}^n_p}{H_p^n} T_{\vec{n}_{jk}} \Pi_p^n, \quad p = j,k.
  \\
&  S^{hrm}( H_j^n, \tilde{H}_j^n)  :=  \begin{pmatrix}   0 \\ \frac{g}{2}( (H_j^n)^2 - (\tilde{H}_j^n)^2   )  \\ 0  \end{pmatrix}.
\end{align}
The function,$\phi$ is any numerical flux function consistent with the 1D component, $F_1(\Pi)$ of the 2D flux, $F(\Pi)$. Here, we consider
the HLL scheme \cite{hll}, namely
\begin{align}\label{fvmhllceqnhllflux}
	\phi(\Pi_L, \Pi_R ) = \begin{cases} F_1(\Pi_L),   & \mbox{ if } s_L \geq 0,\\
	                                      F_1^* := \frac{s_RF_1(\Pi_L) - s_LF_1(\Pi_R) + s_Ls_R(\Pi_R-\Pi_L) }{s_R-s_L}, & \mbox{ if } s_L \leq 0 \leq s_R, \\
	                                      F_1(\Pi_R),  & \mbox{ if }  s_R \leq 0,
	\end{cases}
\end{align}
where $s_L$ and $s_R$ are estimates of the smallest and largest wave speeds in the solution of the associated 1D Riemann problem
\cite{toroshock}. There are several choices for $s_L, s_R$ \cite{toroshock, tororiemann}.
We use the ones given in \cite{rotating} namely
\begin{align}
 s_L = \min_{k}\{ \lambda_k(\Pi_L), \lambda_k(\Pi_R) \}, \quad s_R = \max_{k}\{ \lambda_k(\Pi_L), \lambda_k(\Pi_R) \},
\end{align}
where 
$\lambda_k, k=1,...,M$ , are the eigenvalues of the Jacobian matrix of the system.

And $T_{\vec{n}}$ is a rotation matrix which depends on the normal vector, $\vec{n}=(n_x,n_y)^T$ and $T_{\vec{n}}^{-1}$ is its inverse; they are given by
\begin{align}
	T_{\vec{n}} = \begin{pmatrix}
	                     1    &   0     &      0 \\
	                     0    &   n_x   &      n_y\\
	                     0    &   -n_y   &      n_x
	\end{pmatrix}, \quad
	T_{\vec{n}}^{-1} = \begin{pmatrix}
	                     1    &   0     &       0 \\
	                     0    &   n_x   &      -n_y\\
	                     0    &   n_y   &       n_x
	\end{pmatrix},
\end{align}
see \citep{toroshock}.

\subsection{Schemes for the Channel Models}
Here, we describe the finite volume method to discretize 
the channel models with coupling terms, equations \eqref{horAmodel}, \eqref{horQmodel} and 
\eqref{horqymodel}.
To design a method which reuses existing 1D channel solvers, we discretize the purely 1D channel
models, \eqref{horAmodel} and \eqref{horQmodel} separately from the lateral discharge
model, \eqref{horqymodel}.
Let  $\{x_{i+1/2}\}_{i=1}^{N_{1Dcell}}$ be points
in the 1D grid, $\Omega_h^{1D}$  and $K_i=(x_{i-1/2}, x_{i+1/2}), i=1,2,..,N_{1Dcell}$  be a cell centred
at $x_i=(x_{i-1/2}+ x_{i+1/2})/2$ in  $\Omega_h^{1D}$ . Where $N_{1Dcell}$ is the number of cells
in the 1D grid. Let $\w(\xt)=(A(\xt), Q(\xt))^T$ be a vector of conserved quantities at point, $x$ and time, 
$t$, then the cell average vector, $\w_i^n = (A_i^n, Q_i^n)^T$ in cell $K_i$ is defined as
\begin{equation} \label{horcellaverage}
  \w^n_i := \frac{1}{\Delta x_i} \sint{x}{ K_i}{}{ \w( x, t^n) }
\end{equation}
where $\Delta x_i = x_{i+1/2}-x_{i-1/2},  t^n = t^{n-1} + \Delta t$, and $\Delta t $ is the time step.
\begin{figure}
	\begin{center}
		\includegraphics[scale=0.3]{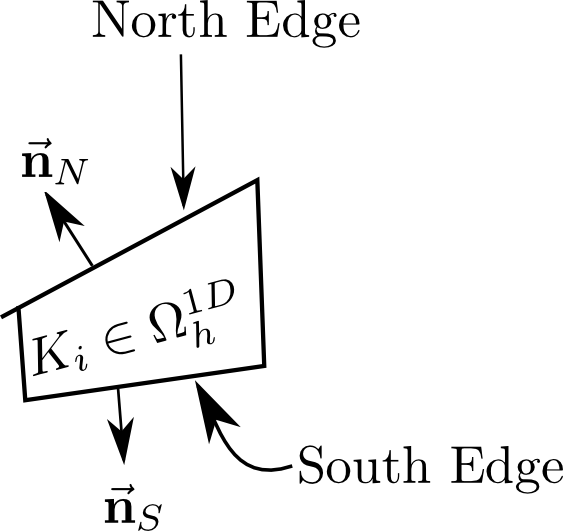}
	\end{center}
	\caption{A single cell, $K_i$ in the 1D channel mesh showing its lateral edges; South edge $e_i^S$ is on the negative $y$-direction
			 while the North edge $e_i^N$ is on the positive $y$-direction. These edges  are the interfaces between the 1D cell and the
			  adjacent 2D floodplain cells.}
	\label{horfig1dcellwithedges}
\end{figure}

For each 1D cell, $K_i \in \Omega_h^{1D}$, the channel lateral boundaries, $y=y_l^w(x)$ and $y=y_r^w(x)$ 
are approximated with straight edges which we call South (S) and North (N) edges (or faces) respectively with
unit normals $\myvec{n}_S=(n_S^x,n_S^y)^T$ and $\myvec{n}_N = (n_N^x, n_N^y)^T$ 
(see figure \ref{horfig1dcellwithedges}). This means that the channel normals, $\myvec{n}_l$ and $\myvec{n}_r$
are approximated with the edge normals $\myvec{n}_S$ and $\myvec{n}_N$ respectively, that is
\begin{equation}\label{horeqnedgenormal}
	\myvec{n}_l \approx \myvec{n}_S, \quad \myvec{n}_r \approx \myvec{n}_N.
\end{equation}
We start by presenting the scheme of the 1D model without coupling terms, next the scheme for channel lateral
discharge model is presented and finally, the discrete coupling term is derived.

\subsubsection{Scheme for 1D Model without Coupling Terms}
We now focus on the 1D channel model \eqref{horAmodel}-\eqref{horQmodel} but without the coupling terms, namely,
\begin{align}\label{hor1dodelnocoupTerm}
\begin{split}
 \partial_t A(\xt)   &+   \partial_x Q(\xt) = 0. 
\\
 \partial_t Q(\xt) 	&+ \partial_x \bigg( \frac{Q^2(\xt)}{A(\xt)} \bigg)   
    =      -    gA(\xt)\pdifft{x}{\etaxt} + gA(\xt) S_f.
\end{split} 
\end{align}
We consider the scheme of \cite{moralesetallargedt} as summarised in \cite{moralesetal2013}.
The scheme is based on the formulation of the St Venant model as presented in \cite{cungeetal80}
and rewrites the model in the quasi-linear form.
\begin{align}
	\partial_t \w + J(\w,B)\partial_x \w = s'(x,\w),
\end{align}
where $\w=(A,Q)^T$, the Jacobian matrix,$J$ is given by
\begin{flushleft}
\begin{align}\label{fvmmorchansolvereqnjacobian}
	& J(\w,B) = \begin{pmatrix}
				0 & 1 \\ c^2 - \uoned^2 & 2 \uoned
	\end{pmatrix},
	\quad \uoned = \frac{Q}{A}, c = \sqrt{ g\frac{A}{B}  }, \quad
    s'(x,\w) = \begin{pmatrix}  0 \\ gA \bigg[ S_o - S_f - \frac{d\honed}{dx} + \frac{1}{B}\frac{dA}{dx}  \bigg]   \end{pmatrix},
\end{align}
\end{flushleft}	
$B$ is the top width at the free-surface, $\honed$ is the water depth from the 1D bottom elevation, $Z_b(x)$ to
the flat free-surface, $\bar{\eta}$ and $S_o = -\frac{dZ_b}{dx}$ is the negative of channel bed slope.
Details about this formulation can be found in  \citep{moralesetal2013}.
The eigenvalues and eigenvectors of $J(\w,B)$ are 
\begin{align*}
\lambda_1(\w,B) = \uoned - c, \,  \lambda_2 (\w,B)= \uoned + c
\mbox{ and }
\textbf{e}_1(\w,B) = (1, \lambda_1(\w,B) )^T, \,  \textbf{e}_2(\w,B) = ( 1, \lambda_2(\w,B) )^T
\end{align*}
respectively.

Define the  Roe averages:
\begin{align}
& \hat{A}_{i+1/2} = \frac{1}{2}( A_i + A_{i+1} )  \quad
  \hat{\uoned}_{i+1/2} = \frac{ \sqrt{A_i} \uoned_i + \sqrt{A_{i+1}}\uoned_{i+1} }{ \sqrt{A_i} + \sqrt{A_{i+1}}  }, \quad
  \hat{\w}_{i+1/2} = \begin{pmatrix} \hat{A}_{i+1/2}    \\  \hat{A}_{i+1/2} \hat{\uoned}_{i+1/2} \end{pmatrix}  
   \label{fvmmorchansolvereqnroeavgu}.
\\&
	\hat{B}_{i+1/2} = \frac{1}{2}( B_{i}+B_{i+1})  \quad
	\hat{\honed}_{i+1/2} = \bigg(\frac{\hat{A}}{\hat{B}} \bigg)_{i+1/2},  \quad
	\hat{c}_{i+1/2} =  \sqrt{  g \hat{\honed}_{i+1/2} }.   \label{fvmmorchansolvereqnroeavgh}
\\
&  (\hat{S_o})_{i+1/2} = \frac{ Z_{b,i+1} - Z_{b,i}   }{ x_{i+1}-x_i } , \quad
  (\hat{S_f})_{i+1/2} =  S_f(\hat{w}_{i+1/2}).
\end{align}

Define $(\Delta p)_{i+1/2} = p_{i+1} - p_i$ for any quantity, $p$. Then, define
\begin{align}
&	( \hat{\alpha}_1 )_{i+1/2} = \bigg[ \frac{\hat{\lambda}_2 \Delta A -    \Delta Q }{2\hat{c}}  \bigg]_{i+/2}, \quad
	( \hat{\alpha}_2 )_{i+1/2} = \bigg[ \frac{-\hat{\lambda}_1 \Delta A +   \Delta Q }{2\hat{c}} \bigg]_{i+1/2}. 
\\
&	(\hat{\beta}_1)_{i+1/2} =  \bigg( -g\frac{\hat{A}} {2 \hat{c}} \left[  ( \hat{S}_0 - \hat{S}_f )\Delta x   - \Delta \honed + \frac{1}{\hat{B}}\Delta A \right] \bigg)_{i+1/2}, \quad
   (\hat{\beta}_2)_{i+1/2} = -( \hat{\beta}_1 )_{i+1/2}.
\end{align}

The Roe averaged eigenvalues and eigenvectors are 
\begin{align}
	(\hat{\lambda}_m )_{i+1/2} := \lambda_m( \hat{\w}_{i+2}, \hat{B}_{i+1/2} ), \quad (\hat{\textbf{e}}_m)_{i+1/2} := \textbf{e}_m(\hat{\w}_{i+1/2}, \hat{B}_{i+1/2} ), \quad m=1,2.
\end{align}
The artificial viscosity (entropy fix), $\hat{\nu}$ to correct the entropy problem associated with the Roe method \citep{moralesetal2013} is given by
\begin{align}
&	( \hat{\nu}_m)_{i+1/2} = \begin{cases}   \frac{1}{4} \bigg[ (\lambda_m)_{i+1} - (\lambda_m)_i    \bigg], & \mbox{ if } (\lambda_m)_i < 0 < (\lambda_m)_{i+1} \\
							0 ,  &  \mbox{ else }
	\end{cases}, \quad	m=1,2. 
\end{align}

Hence, the numerical scheme of \citep{moralesetal2013, moralesetallargedt} 
for the 1D channel model without coupling term, \eqref{hor1dodelnocoupTerm}
is given by 
\begin{align}\label{horeqn1dscheme}
	\w_i^{n+1*} = \w_i^n - \frac{\Delta t}{\Delta x} \bigg[
						\sum_{m=1}^2  \bigg( \hat{\gamma}_m^+ \hat{\textbf{e}}_m  \bigg)_{i-1/2}
					 +  \sum_{m=1}^2  \bigg( \hat{\gamma}_m^- \hat{\textbf{e}}_m  \bigg)_{i+1/2} 
					 \bigg]^n,
\end{align}
where 
\begin{align}
	\bigg( \hat{\gamma}_m^{\pm} \bigg)_{i+1/2} = \bigg[ \frac{1}{2}[ 1 \pm sgn(\hat{\lambda})  ]\hat{\gamma}  \pm \hat{\nu}\hat{\alpha}   \bigg]_{m,i+1/2},  \quad 
	\bigg( \hat{\gamma}_m \bigg)_{i+1/2}  = \bigg(  \hat{\lambda}\hat{\alpha} - \hat{\beta}      \bigg)_{m,i+1/2},  \quad	m=1,2.
\end{align}

\subsubsection{Approximating Channel Lateral Discharge}\label{horsecqyscheme}
\begin{figure}[ht!]
  \begin{center}
   \subfigure[A 1D cell $K_i$ in the channel grid. $e_i^N$ and $e_i^S$ are the North and South edges respectively.
             $T_{ij}^N, j=1,2...,N_n$ are its adjacent/neighbour 2D floodplain cells on the North edge, while
      		  $T_{ij}^S, j=1,2,...,N_s$ are the adjacent 2D floodplain cells on the South edge. 
              $K_{i-1}, K_{i+1}$ are the left and right neighbours of $K_i$ in 1D channel grid, $\Omega_h^{1D}$.]
              {\label{fig1d2dcells} \includegraphics[width=0.46\textwidth]{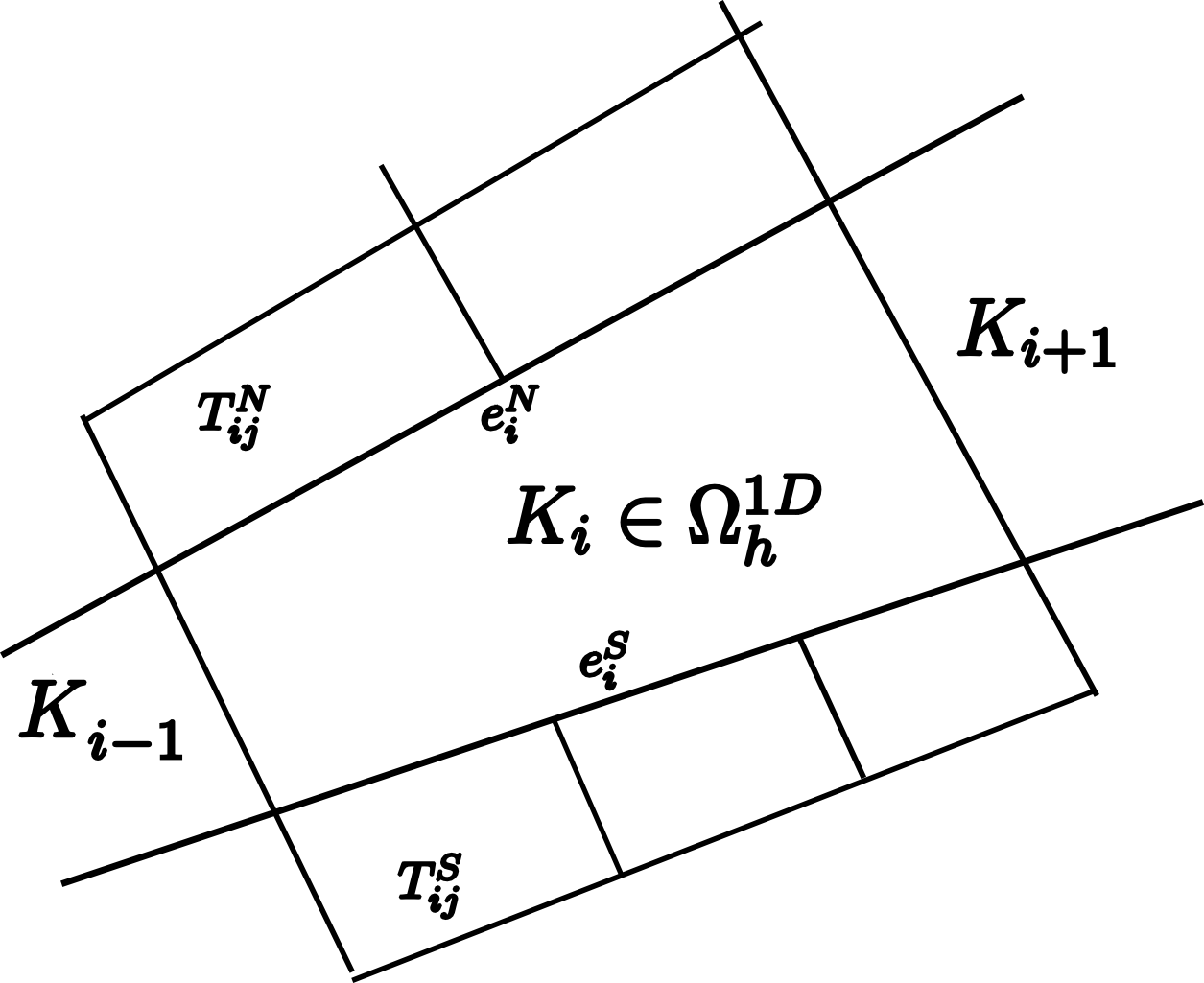}}
               \hspace{5mm}
   \subfigure[The single 1D cell subdivided into two subcells $K^N_i$ and $K^S_i$ which are then viewed as 2D cells.
             $e_{xf}, e_{NS}, e_{xb}$ and $e_{ij}^N, j=1,2,...N_n$ are the edges of $K_i^N$
             with their outward unit normal vectors as indicated. Similarly, the edges of $K_i^S$
              and their normal vectors are indicated.]
   			{ \label{fig1d2dsubcells} \includegraphics[width=0.46\textwidth]{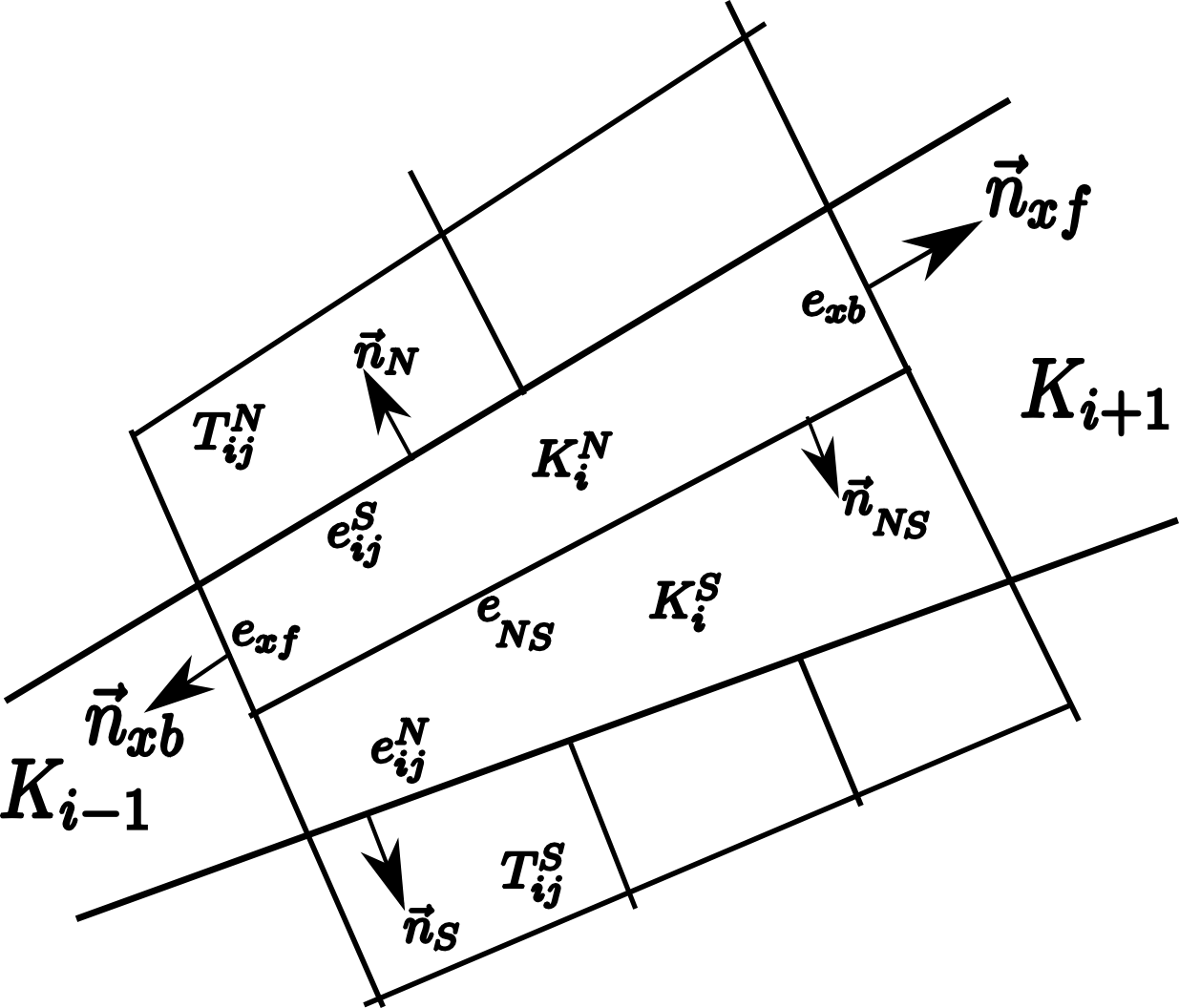}}
  \end{center}
  \caption{ To the left is a 1D channel cell and its adjacent 2D floodplain cells while to the right is the 1D cell subdivided into
            two subcells viewed as 2D cells}
\end{figure}
Here, the goal to solve the lateral discharge model, \eqref{horqymodel} along the channel.
Consider the 1D channel cell, $K_i \in \Omega_h^{1D}$, with cell average vector $\w_i^n=(A_i^n,Q_i^n)^T$ and denote by 
$T_{ij}^N \in \Omega_h^{2D}, j=1,2,...,N_n$ and $T_{ij}^S \in \Omega_h^{2D}, j=1,2,...,N_s$, 
the $j$-th 2D floodplain cells adjacent to $K_i$ on its North edge $e_i^N$ and South edge $e_i^S$ respectively, see figure \ref{fig1d2dcells}.
Let the cell averages in the adjacent 2D cells $T_{ij}^N$ and $T_{ij}^S$ be 
\begin{equation}
  \begin{split}
   (\Pi^N)_{i,j}^n  = ( (H^N)_{i,j}^n, (q_x^N)_{i,j}^n, (q_y^N)_{i,j}^n )^T
  \mbox{ and  }  
  (\Pi^S)_{i,j}^n  = ( (H^S)_{i,j}^n, (q_x^S)_{i,j}^n, (q_y^S)_{i,j}^n )^T,
\end{split}
\end{equation} 
respectively.
$N_n$ and $N_s$ are the number of the adjacent 2D cells on the North and South edges, respectively of $K_i$.
$(H^N)_{ij}^n, (q_x^N)_{ij}^n$ and $ (q_y^N)_{ij}^n $ are the average water depth,  average discharge along $x$-direction and
average discharge along $y$-direction respectively, in 2D cell $T_{ij}^N$ while
$ (H^S)_{ij}^n, (q_x^S)_{ij}^n$ and $ (q_y^S)_{ij}^n $ are those of cell $T_{ij}^S$.

To discretize the lateral discharge model \eqref{horqymodel} in $K_i$, we subdivide  $K_i$ into two subcells,
$K_i^N$ and $K^S_i$ and view them as 2D cells within the channel, see figure \ref{fig1d2dsubcells}. 
Let $(\w^N)_i^n$ and $(\w^S)_i^n$ be the 2D cell average vectors in the
subcells, $K_i^N$ and $K_i^S$ respectively.
Then, we define them as
\begin{equation}\label{horeqnsubcellsavg}
	(\w^N)_i^n = ( \honed_i^n, \honed_i^n \uoned_i^n, (q_y^N)_i^n )^T \mbox{  and  }
	(\w^S)_i^n = ( \honed_i^n, \honed_i^n \uoned_i^n, (q_y^S)_i^n )^T, 
\end{equation}
where $\honed_i^n$ and  $\uoned_i^n=\frac{Q^n_i}{A_i^n}$  are the 1D cell average water depth  and section-averaged velocity
in the channel, and  $(q_y^S)_i^n, (q_y^N)_i^n$ are computed at every time step
using the scheme presented below.

The motivation to compute $q_y^{N/S}$ is to apply a well-balanced scheme to the model,\eqref{horqymodel} 
in the subcells $K_i^{N/S}$ by taking the bottom to be flat across all the edges within the
channel, $e_{xb}, e_{NS}, e_{xf}$ (see figure \ref{fig1d2dsubcells}).
To this end, we define the following:
\begin{align}\label{horeqn-qn-data} 
	h_{2ij}^N &= max(0, \bar{\eta}_i^n - z_{b,ij}^N), 
\quad
	(\tilde{\w}^N)_i^n = (h_{2ij}^N, h_{2ij}^N \uoned_i^n, h_{2ij}^N (v^N)_i^n)^T, 
	\quad
	(v^N)_i^n = (q_y^N)_i^n/\honed_i^n,  
\end{align}
where $\bar{\eta}_i^n = \honed_i^n + Z_{b,i}$ is the discrete flat free surface elevation in 1D cell, $K_i$
and $z_{b,ij}^N$ is the bed elevation of the adjacent 2D cell $T_{ij}^N$, see figure \ref{fig1d2dsubcells}. 
Therefore, we propose the following hydrostatic reconstruction scheme \citep{audusseetal2004} for
the lateral discharge in subcell $K_i^N$:
\begin{align}\label{horeqnqn}
\begin{split}
	(q_y^N)_i^{n+1} =& (q_y^N)_i^n - \frac{\Delta t}{|K_i^N|} 
					   \bigg[					 
						     |e_{xb}|\phi_3^{2D}( (\w^N)_i^n, (\w^N)_{i-1}^n, \vec{n}_{xb}  )			    					    
				+
					    	|e_{xf}|\phi_3^{2D}( (\w^N)_i^n, (\w^N)_{i+1}^n, \vec{n}_{xf}  )
                \\
				&
					    	+				    	
					    	|e_{NS}|\phi_3^{2D}( (\w^N)_i^n, (\w^S)_i^n, \vec{n}_{NS}  )  
			    	 \bigg]	
			       	- \frac{\Delta t}{|K_i^N|} 
					    	\sum_{j=1}^{Nn} |e_{ij}^N|
					    	\bigg[ 
					    		   \phi_3^{2D}( (\tilde{\w}^N)_i^n , (\Pi^N)^n_{ij},  \vec{n}_{N} )
                    \\
                &
                     +
				   	       \frac{g}{2} \vec{n}_{N} \cdot \begin{pmatrix}  0 \\  (\honed_i^n)^2 - (h_{2ij}^N)^2  \end{pmatrix}
					  		\bigg],
\end{split}					  					
\end{align}
$\phi_{3}^{2D}( w_L, w_R, \vec{n} )$ denotes the 3rd component of numerical flux, 
$\phi^{2D}( w_L, w_R,\vec{n}) := T_{\vec{n}}^{-1}\phi( T_{\vec{n}}w_L, T_{\vec{n}}w_R)$.
The quantities; $ |e_{xb}|, |e_{xf}|, |e_{NS}| $ and $|e_{ij}^N|$ are the lengths of the corresponding edges of $K_i^N$ and
$\vec{n}_N$ is the outward unit normal of $K_i^N$ towards $T_{ij}^N$
(see figure \ref{fig1d2dsubcells}).

Similarly, by defining
\begin{align}\label{horeqn-qs-data} 
	h_{2ij}^S &= max(0, \bar{\eta}_i^n -  z_{b,ij}^S), 
	\quad
	(\tilde{\w}^S)_i^n := (h_{2ij}^S, h_{2ij}^S \uoned_i^n, h_{2ij}^S (v^S)_i^n)^T , 
	\quad
	(v^S)_i^n := (q_y^S)_i^n/\honed_i^n,  
\end{align}
we propose the following scheme for the lateral discharge in $K_i^S$:
\begin{align}\label{horeqnqs}
\begin{split}
	(q_y^S)_i^{n+1} =& (q_y^S)_i^n - \frac{\Delta t}{|K_i^S|} 
					     \bigg[					 
						    |e_{xb}|\phi_3^{2D}(  (\w^S)_i^n, (\w^S)_{i-1}^n , \vec{n}_{xb} )			    					    						    
						    +
					    	|e_{xf}|\phi_3^{2D}(  (\w^S)_i^n, (\w^S)_{i+1}^n, \vec{n}_{xf}  ) 
                \\
                &
					    	+				    	
					    	|e_{SN}|\phi_3^{2D}( (\w^S)_i^n, (\w^N)_i^n, \vec{n}_{SN}  )  				
				     \bigg]		
				   	-      	\frac{\Delta t}{|K_i^S|} 
					    	\sum_{j=1}^{Ns} |e_{ij}^S|
					    	\bigg[ 
					    		   \phi_3^{2D}(  (\tilde{\w}^S)_i^n , (\Pi^S)_{ij}^n, \vec{n}_{S} ) 
                     \\
                &
					    	        +
					    	       \frac{g}{2} \vec{n}_{N} \cdot\begin{pmatrix}  0 \\  (\honed_i^n)^2 - (h_{2ij}^S)^2  \end{pmatrix}
					  		\bigg],
\end{split}					  			  								
\end{align}
where
$|e_{SN}|$ is length of edge, $e_{SN}$ between $K_i^N$ and $K_i^S$, and $z_{b,ij}^S$ is the bed elevation of  2D cell, $T_{ij}^S$.
The initial values of $(q_y^S)_i^n, (q_y^N)_i^n$ are obtained as explained in the following remark.
\begin{myremark}
 At initial time ($n=0$), only the lateral discharge, $(q_y)_i^0$ for the full cell, $K_i$ is given.
Then we initialize $(q_y^N)_i^0$ and $(q_y^S)_i^0$ to be equal to it, namely
\begin{equation}\label{horinitialqySN}
	 (q_y^N)_i^0  = (q_y^S)_i^0  = (q_y)_i^0.
\end{equation}
For other times, ($n > 0$), we compute $(q_y^{S/N})_i^n$ using the scheme described above.
\end{myremark}

\subsubsection{Discrete Coupling Terms}\label{horseccopscheme}
To discretize the coupling term, let us denote by $f_{i,j}^S$ the approximation of a function, $f$
at the edge, $e^S_{ij}$ between 1D cell $K_i$ and 2D cell $T_{ij}^S$ (see figure \ref{fig1d2dsubcells}).   
Then the approximation, $f|_L$ of $f$ over the
entire South edge, $e_i^S$ of $K_i$, is given by averaging over all edges on the South edge, namely
\begin{equation}
 ( f|_L )_i = \sum_{j=1}^{Ns} \bigg( f_{i,j}^S \frac{|e_{ij}^S|}{|e_i^S|}  \bigg),
\end{equation}
where
$|e^S_{ij}|$ is the length of edge, $e_{ij}^S$;  $|e^S_i|$ is sum of all edges of $K_i$ on the 
South side.
Similarly $ (f|_R)_i  =  \sum_{j=1}^{Nn} \bigg( f_{i,j}^N \frac{|e_{ij}^N|}{|e_i^N|} \bigg)    $ for North edge,
where $|e^N_{ij}|$ is the length of edge between $K_i$ and $T_{ij}^N$, and $|e_i^N|$ is the sum of all edges on North 
side of $K_i$.
Hence we can approximate the coupling term as
\begin{align}\label{hordiscretecouplinterm}
  \Phi_i^n   =   \sum_{j=1}^{Ns} \Psi_{i,j}^S +     \sum_{j=1}^{Nn} \Psi_{i,j}^N,
\end{align}
where
\begin{align}
   \Psi^S_{i,j}    
                =   \begin{pmatrix}  \frac{1}{n^y_S}  f^{1,S}_{i,j}  \\
                           \frac{1}{n^y_S}   f^{2,S}_{i,j} -  \frac{n^x_S}{n^y_S} \frac{g}{2}  (H^{S*}_{i,j})^2 
                    \end{pmatrix}
                    \frac{|e_{ij}^S|}{|e_i^S|}, 
\quad
  \Psi^N_{i,j}      
                =   \begin{pmatrix}  - \frac{1}{n^y_N}     f^{1,N}_{i,j}  \\
                          - \frac{1}{n^y_N}   f^{2,N}_{i,j} +  \frac{n^x_N}{n^y_N}\frac{g}{2} (H^{N*}_{i,j})^2
                    \end{pmatrix}
                    \frac{|e_{ij}^N|}{|e_i^N|} , 
\end{align}
are the discrete coupling terms at the edges  $e_{ij}^S$ and $e_{ij}^N$ respectively (see figure \ref{fig1d2dsubcells}). 
$H^{S*}_{i,j}, f^{1,S}_{i,j}$ and $f^{2,S}_{i,j}$ are respectively, the discrete water depth, first and second components of 
2D numerical flux at edge $e_{ij}^S$. While  $H^{N*}_{i,j}, f^{1,N}_{i,j}$ and $f^{2,N}_{i,j}$ are respectively, the 
water depth, first and second components of 2D numerical flux at edge, $e_{ij}^N$.
We now focus on how to compute them.

Given the 1D cell average, $\w_i^n$ in $K_i$ from which we obtain the cell average
$(\w^S)_i^n$ in the subcell, $K_i^S$ (using equation \eqref{horeqnsubcellsavg}).
Then we directly approximate $f^{1,S}_{i,j}$ and $f^{2,S}_{i,j}$ by computing the 2D numerical flux,
$\phi^{2D}( (\tilde{\w}^S)^n_i, (\Pi^S)_{ij}^n, \vec{n}_S ) $,
at edge, $e_{ij}^S$  (see figure \ref{fig1d2dsubcells}), namely
\begin{equation}\label{horsch1dflux}
\begin{split}
                 f^{1,S}_{i,j} = \phi_1^{2D}( (\tilde{\w}^S)^n_i, (\Pi^S)_{ij}^n , \vec{n}_S ),  \quad
                 f^{2,S}_{i,j} = \phi_2^{2D}( (\tilde{\w}^S)^n_i, (\Pi^S)_{ij}^n , \vec{n}_S ).         
\end{split}             
\end{equation}
Similarly, by using $(\tilde{\w}^N)^n_i$ and $(\Pi^N)_{ij}^n$  we  approximate  $f^{1,N}_{i,j}$ and $f^{2,N}_{i,j}$
using
\begin{equation}
                \begin{split} f^{1,N}_{i,j} = \phi_1^{2D}((\tilde{\w}^N)^n_i, (\Pi^N)_{ij}^n, \vec{n}_N ), \quad
                              f^{2,N}_{i,j} = \phi_2^{2D}((\tilde{\w}^N)^n_i, (\Pi^N)_{ij}^n, \vec{n}_N ), 
                \end{split}
\end{equation}
The hydrostatically reconstructed quantities, $(\tilde{\w}^N)^n_i$ and $(\tilde{\w}^S)^n_i$
are defined in equations \eqref{horeqn-qn-data} and \eqref{horeqn-qs-data} respectively.

To approximate $H_{i,j}^{S*}$ and  $H_{i,j}^{N*}$, we propose to adapt the hydrostatic reconstruction 
approach \citep{audusseetal2004}, namely 
\begin{align}
	H_{i,j}^{S*} = \max(h_{2ij}^S, (H^S)_{i,j}^n),
\quad
	H_{i,j}^{N*} = \max(h_{2ij}^N, (H^N)_{i,j}^n),
\end{align}
We therefore summarise the discrete coupling term as
%
%
\begin{align}\label{horeqnnumcouplingterm}
  \begin{split}
    \Phi_i^n  &=   
              \frac{1}{|e_i^S|} \sum_{j=1}^{Ns}  |e_{ij}^S|
              \begin{pmatrix} 
                 \frac{1}{n^y_S}  \phi_1^{2D}(  (\tilde{\w}^S)_i^n, (\Pi^S)_{ij}^n,\vec{n}_S ) 
                \\
                 \frac{1}{n^y_S}  \phi_2^{2D}(  (\tilde{\w}^S)_i^n, (\Pi^S)_{ij}^n, \vec{n}_S )  
                - \frac{g}{2}\frac{n_S^x}{n_S^y} \bigg[\max( h_{2ij}^S, (H^S)_{ij}^n) \bigg ]^2
              \end{pmatrix}
              \\
              &
              -
               \frac{1}{|e_i^N|} \sum_{j=1}^{N_n}  |e_{ij}^N|
               \begin{pmatrix} 
                 \frac{1}{n^y_N}  \phi_1^{2D}( (\tilde{\w}^N)_i^n, (\Pi^N)_{ij}^n, \vec{n}_N  ) 
                \\
                 \frac{1}{n^y_N}  \phi_2^{2D}(  (\tilde{\w}^N)_i^n, (\Pi^N)_{ij}^n, \vec{n}_N )  
                - \frac{g}{2} \frac{n_N^x}{n_N^y} \bigg[ \max( h_{2ij}^N, (H^N)_{ij}^n) \bigg]^2
              \end{pmatrix},
   \end{split} \\
    & n_N^y, n_S^y \neq 0. \notag
\end{align}

\subsection{Summary of the Channel Flow Solver}\label{summarychannelsolver}
The complete scheme for the channel flow model with coupling term,  \eqref{horAmodel}-\eqref{horQmodel} is
\begin{equation} \label{horeqnfully1dschemeFinal}
    \w_i^{n+1} = \w_i^{n+1*} + \Delta t \Phi_i^n,
\end{equation}
where $\w_i^{n+1*}$ is the solution of the purely channel model without the coupling terms,\eqref{hor1dodelnocoupTerm}
which is given in \eqref{horeqn1dscheme}. $\Phi_i^n$ is the discrete coupling term summarised in \eqref{horeqnnumcouplingterm}.
The channel lateral discharges are computed using the schemes in \eqref{horeqnqn}
and \eqref{horeqnqs}.
%

%

\subsection{Summary of the Horizontal Coupling Method}\label{summary-hcm}
\begin{figure}[ht!]
	\begin{center}
		\includegraphics[width=0.95\textwidth, height=190mm]{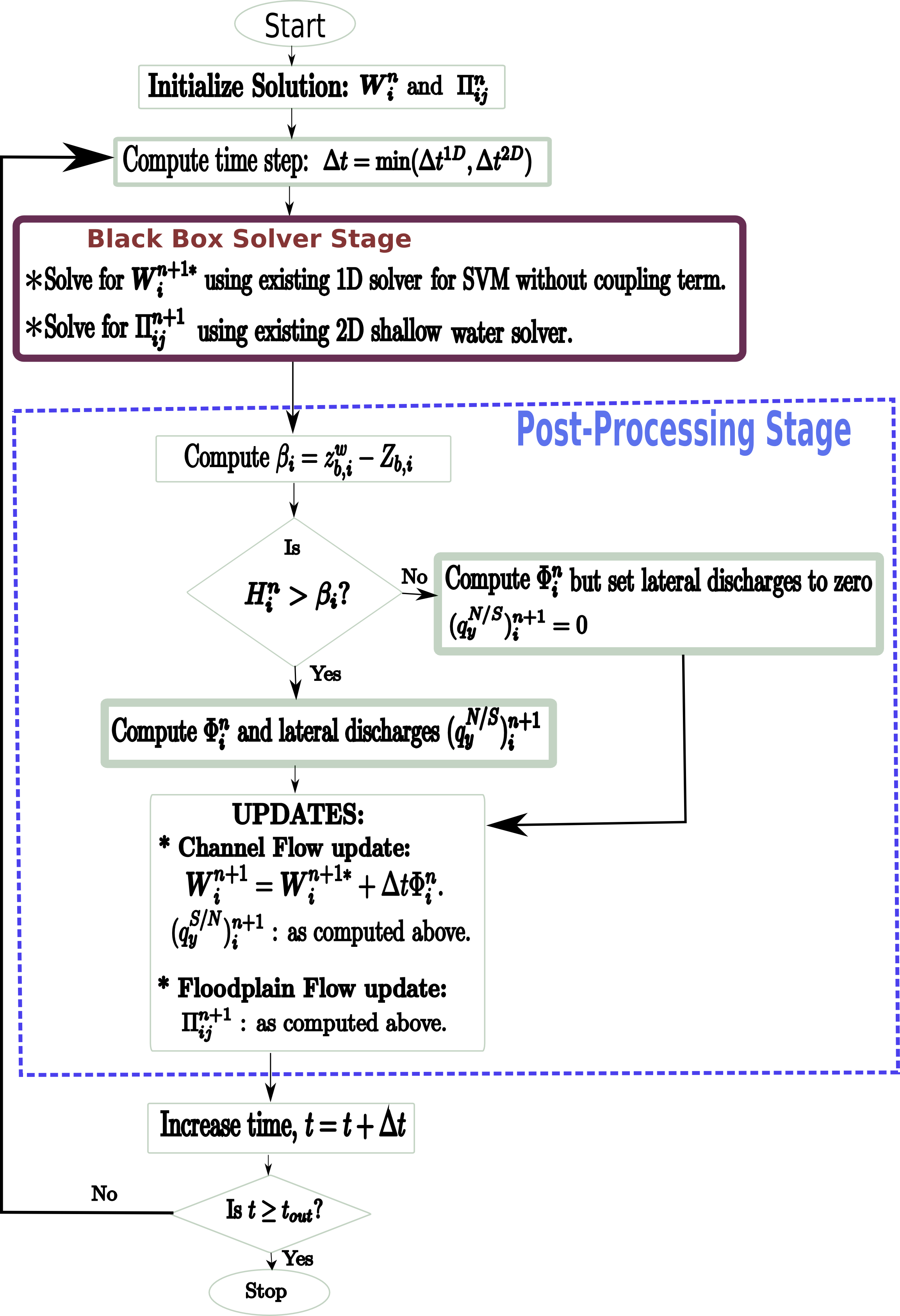}
	\end{center}
	\caption{Flow Chart for implementation of the Horizontal Coupling Method (HCM)}
	\label{figflowchart}
\end{figure}

The channel flow model, \eqref{horAmodel}, \eqref{horQmodel} and \eqref{horqymodel}
is simulated as summarised in section \ref{summarychannelsolver} while the flood flow model, 
\eqref{fv2dswewtfriceqnmodel} is solved with the 2D solver,  \eqref{fvm2dswewbedfrictiongeneral}.
At 2D/1D edge,  the 2D numerical flux is computed by using the 2D cell, $T_j \in \Omega_h^{2D}$ 
averages and the averages obtained from the adjacent channel subcell, ($K_i^{N}$ or $K_i^S$) 
as described in \eqref{horeqnsubcellsavg}. A flow chart for the implementation
of the HCM is given in figure \ref{figflowchart}.

\section{Properties of the HCM}\label{horiprops}
We discuss a few properties of the method in this section.
\begin{mydefinition}[Well Balance of Lake at rest]
	Assuming that the existing numerical schemes for the uncoupled 1D and 2D models are well balanced
	with respect to lake at rest, then the coupled scheme is said to be well balanced with respect to lake at rest
	if the coupling term vanishes whenever the lake at rest condition holds.
\end{mydefinition}
\begin{mytheorem}
	The coupling term  derived in equation \eqref{horeqnnumcouplingterm} leads to a fully well-balanced 
	scheme with respect to lake at rest.
\end{mytheorem}
\begin{myproof}
Assuming that the condition of water at rest holds, then
\begin{equation*}
	\bar{\eta}_i^n = (\eta^N)_{ij}^n = (\eta^S)_{ij}^n  \quad \forall j
\end{equation*}
where $(\eta^S)_{ij}^n $ and $(\eta^N)_{ij}^n $ are the free surface elevation in the adjacent 2D cells, $T_{ij}^S$ and $T_{ij}^N$
respectively.
Hence, 
\begin{align*}
\begin{split}
	&h_{2,ij}^S := \max(0, \bar{\eta}_i^n-z_{b,ij}^S) = \max(0, (\eta^S)_{ij}^n-z_{b,ij}^S)  = \max(0, (H^S)_{ij}^n ) = (H^S)_{ij}^n, 
	\\
	&h_{2,ij}^N := \max(0, \bar{\eta}_i^n-z_{b,ij}^N) = \max(0, (\eta^N)_{ij}^n-z_{b,ij}^N)  = (H^N)_{ij}^n.  
	\end{split} 
\end{align*}
Therefore, 
\begin{align}\label{horeqnhwellbalance}
	\max( h_{2ij}^S, (H^S)_{ij}^n) =  h_{2ij}^S \mbox{ and }
	\max( h_{2ij}^N, (H^N)_{ij}^n) =  h_{2ij}^N.
\end{align}
Since all velocities (and discharges) are zero, then
$$
	(\Pi^S)_{ij}^n = ( (H^S)_{ij}^n, 0, 0 )^T =  ( h_{2ij}^S, 0, 0  )^T  = (\tilde{\w}^S)_i^n.
$$
Hence, by the consistency  of the numerical flux, $\phi^{2D}$ with the physical flux, $F( \cdot )$, we have
$$
	\phi ^{2D}( (\tilde{\w}^S)_i^n, (\Pi^S)_{ij}^n , \myvec{n}_S  ) = F ( (\tilde{\w}^S)_i^n ) \cdot \myvec{n}_S
	=   ( 0,  n_S^x\frac{g}{2} \bigg[h_{2ij}^S \bigg]^2,  n_S^y\frac{g}{2} \bigg[h_{2ij}^S \bigg]^2  )^T.
$$
That is 
\begin{align}\label{horeqnsflux}
	\phi_1^{2D}( (\tilde{\w}^S)_i^n, (\Pi^S)_{ij}^n,  \myvec{n}_S  ) = 0, \quad 
	\phi_2^{2D}( (\tilde{\w}^S)_i^n, (\Pi^S)_{ij}^n,  \myvec{n}_S ) = n_S^x\frac{g}{2} \bigg[ h_{2ij}^S \bigg]^2.
\end{align}
Similarly, 
\begin{align*}
\begin{split}
	(\Pi^N)_j^n = ( (H^N)_{ij}^n, 0, 0  )^T = (\tilde{\w}^N)_i^n 
\mbox{ and } 
	\phi^{2D} ( (\tilde{\w}^N)_i^n, (\Pi^N)_{ij}^n,  \myvec{n}_N    ) 
		= ( 0,  n_N^x\frac{g}{2} \bigg[h_{2ij}^N\bigg]^2,  n_N^y\frac{g}{2} \bigg[h_{2ij}^N \bigg]^2  )^T.
\end{split}		
\end{align*}
So that
\begin{align}\label{horeqnnflux}
	\phi_1^{2D}( (\tilde{\w}^N)_i^n, (\Pi^N)_{ij}^n, \myvec{n}_N ) = 0, \quad 
	\phi_2^{2D}( (\tilde{\w}^N)_i^n, (\Pi^N)_{ij}^n, \myvec{n}_N ) = n_N^x\frac{g}{2} \bigg[h_{2ij}^N \bigg]^2.
\end{align}		

Therefore, using equations \eqref{horeqnhwellbalance}-\eqref{horeqnnflux}, then the discrete coupling term
in \eqref{horeqnnumcouplingterm} becomes
\begin{align*}
  \begin{split}
    \Phi_i^n   &=
              \frac{1}{|e_i^S|} \sum_{j=1}^{Ns}  |e_{ij}^S|
              \begin{pmatrix} 
                  0
                \\
                 \frac{1}{n^y_S}  n_S^x\frac{g}{2} \bigg[h_{2ij}^S \bigg]^2
                - \frac{g}{2}\frac{n_S^x}{n_S^y} \bigg[ h_{2ij}^S \bigg ]^2
              \end{pmatrix}
               -
               \frac{1}{|e_i^N|} \sum_{j=1}^{Nn}  |e_{ij}^N|
               \begin{pmatrix} 
             	0
                \\
                \frac{1}{n^y_N}  n_N^x\frac{g}{2}\bigg[h_{2ij}^N \bigg]^2
                - \frac{g}{2} \frac{n_N^x}{n_N^y} \bigg[h_{2ij}^N\bigg]^2
              \end{pmatrix}
              = 
              \begin{pmatrix}
              	0 \\0
              \end{pmatrix}.
\end{split}              
\end{align*}
as claimed.
\end{myproof}
We now introduce the concept of "No-Numerical Flooding".

\begin{mydefinition}[No Numerical Flooding Property]
  We shall say that a 2D/1D coupling scheme preseves the \textbf{No Numerical Flooding property}, if all
	its coupling terms vanish whenever there is no flooding or draining.
\end{mydefinition}
\begin{mytheorem} 
	The scheme \eqref{horeqnfully1dschemeFinal} preserves the no numerical
	flooding property.
\end{mytheorem}
\begin{myproof}
If no flooding, then 
\begin{align*}
   \bar{\eta}_i^n \leq z_{b,i}^w \leq z^S_{b,ij}  \mbox{ and } \bar{\eta}_i^n \leq z^N_{b,ij}  \quad  \forall j 
   \Longrightarrow h_{2ij}^N=h^S_{2ij} = 0, \quad  \mbox{ by definition (\eqref{horeqn-qn-data}, \eqref{horeqn-qs-data})}.
\end{align*}
Hence, 
\begin{align*}
   (\tilde{\w}^N)^n_i = (\tilde{\w}^S)^n_i = (0, 0, 0)^T \quad \mbox{ by definition (\eqref{horeqn-qn-data}, \eqref{horeqn-qs-data})}.
\end{align*}   
Again, since floodplain is dry, we have 
\begin{align*}
(H^N)^n_{ij}=(H^S)^n_{ij} = 0 \quad \forall j \Longrightarrow (\Pi^N)^n_{ij} = (\Pi^S)^n_{ij} = (0,0,0)^T \quad \forall  j. 
\end{align*}
These give the numerical fluxes:
\begin{align*}
 \phi^{2D}((\tilde{\w}^N)^n_i , (\Pi^N)^n_{ij},\vec{n}_N  ) = \phi^{2D}( (\tilde{\w}^S)^n_i , (\Pi^S)^n_{ij}, \vec{n}_S ) = (0,0,0)^T.
\end{align*}
So all the flux terms in $\Phi_i^n$ are zero.
Finally,
\begin{align*}
	 \max( h_{2ij}^N, (H^N)_{ij}^n) = \max( h_{2ij}^S, (H^S)_{ij}^n) = \max( 0,0) = 0.
\end{align*}	 
Therefore,
$ \Phi_i^n  = 0$.
Which means that no water is gained from or lost to the floodplain as required.
\end{myproof}

\section{Numerical Results}\label{hsecnumresults}
%
In this section, we present some numerical experiments to investigate the performance of the proposed method.
We use full 2D simulation results as the reference solution and compare these results with those of
the HCM  and of the flux-based method (FBM) of \cite{bladeetal2012}.  All the algorithms are 
implemented in a C++ code and the experiments are run on the Cluster of Workstations
(COW) of the Centre for Scientific Computing, University of Warwick, United Kingdom.

%
\subsection{Test Case 1 : Dam-Break Flow into a Flat Floodplain}
\begin{figure}[ht!] 
  \begin{center}
	\includegraphics[width=\textwidth]{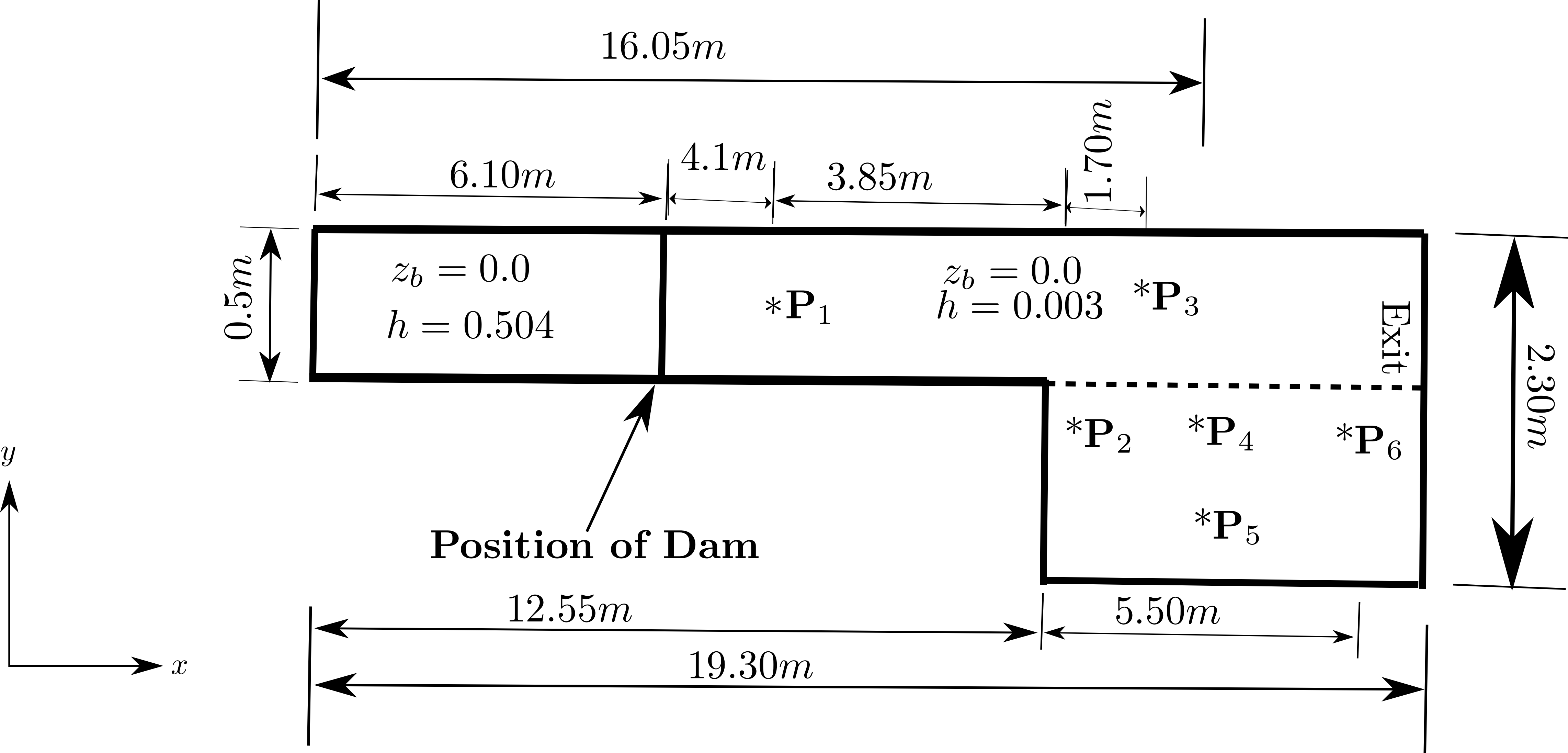}
  \end{center}
  \caption{Top view of Channel and Floodplain for river-flooding problem for test case 1}
  \label{fvmnumfigwhycouplinggeometry}
\end{figure}

The first test case is suggested in \citep{moralesetal2013}. 
The setup consists of a dam break flow in a 19.3 meter long, 0.5 meter constant width flat channel 
with adjacent flat floodplain, see figure \ref{fvmnumfigwhycouplinggeometry}. 
The National Laboratory of Civil Engineering in the IST in Portugal designed
and measured this test case \citep{viseu3numerical, moralesetal2013}.
A reservoir is located from the left end of the channel to 6.10 metres 
(position of dam in figure \ref{fvmnumfigwhycouplinggeometry}). The initial condition is 
\begin{align*}
&	H(x,y,0) = \begin{cases}
						0.504, & \mbox{ at the reservoir, that is }  0\leq x \leq 6.10 \mbox{ and } 1.8\leq y \leq 2.3 \\
						0.003, & \mbox{ elsewhere}
		       \end{cases},\\
&	u(x,y,0) = v(x,y,0)	= 0 \quad \mbox{ everywhere}.	       
\end{align*}
The manning coefficient, $n$ (see section \ref{sec-derive-models}) for both channel and floodplain is $0.009$s/m$^{1/3}$ and the boundaries are all closed walls except the right side
as indicated in figure \ref{fvmnumfigwhycouplinggeometry}. The labels $P_1, P_2, \dots, P_6$ are probe points in the flow domain.
More about this test case can be found in \citep{moralesetal2013,viseu3numerical}.

Here, a full 2D simulation was run with a grid of $193\times25$ cells in the channel and $68\times90$
cells in the floodplain, while a simulation with the proposed
method was run with a grid of $68\times 90$ cells for the floodplain and $193\times2$ cells in the
channel, and the simulation using the FBM was also run with a grid of $68\times90$ cells in floodplain and $193\times1$
in the channel, see table \ref{horinumtabletest1}. 

\begin{table}[ht!] 
\begin{center}
\begin{tabular}{ | l |  c | c |  c | c |}
\hline
                 &      Channel Grid   &  Floodplain Grid &  No. of time steps   &  Processor time (in seconds)           \\
\hline
Full 2D          &     $193\times25$   &  $68\times90$   &    3,669              &     3,110.31  \\  
\hline
 HCM            &      $193\times2$   & $68\times90$   & 2,616              &     1,420.4    \\ 
\hline   
FBM             &      $193\times1$  & $68\times90$    & 2,592              &     1,311.26   \\ 
\hline 
\end{tabular}
\end{center}
\caption{Grid cells, simulation times and number of time steps}
\label{horinumtabletest1}
\end{table}

\begin{figure}[ht!] 
	\includegraphics[width=\linewidth]{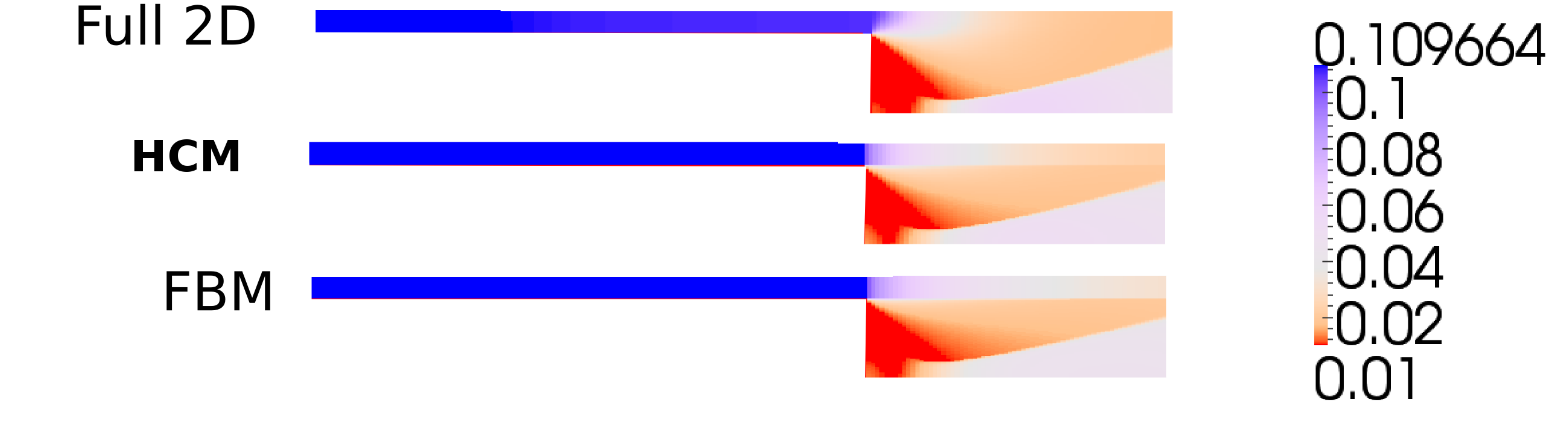}
  \caption{Comparison of the final free surface elevation after ten seconds for test 1.} 
  \label{horinumfigtest1eta}
\end{figure}

\begin{figure}[ht!] 
  \begin{center}
	\subfigure{\includegraphics[width=0.48\textwidth]{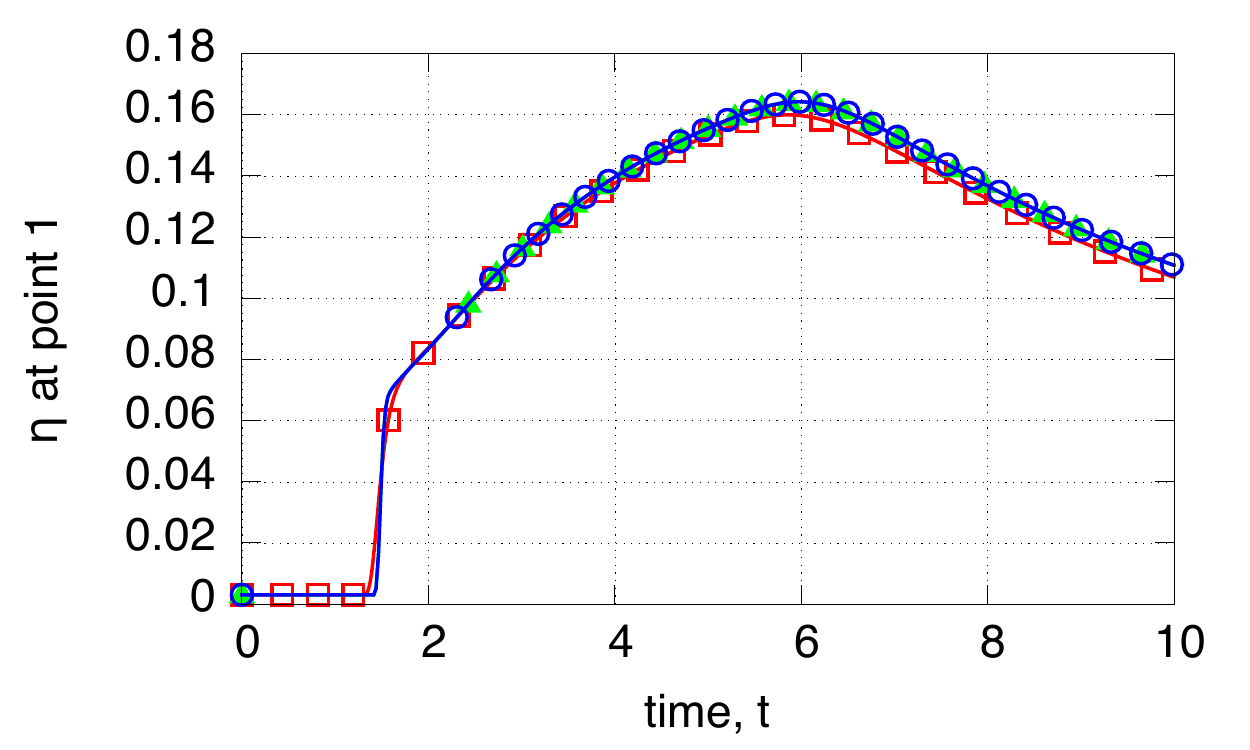}}
	\subfigure{\includegraphics[width=0.48\textwidth]{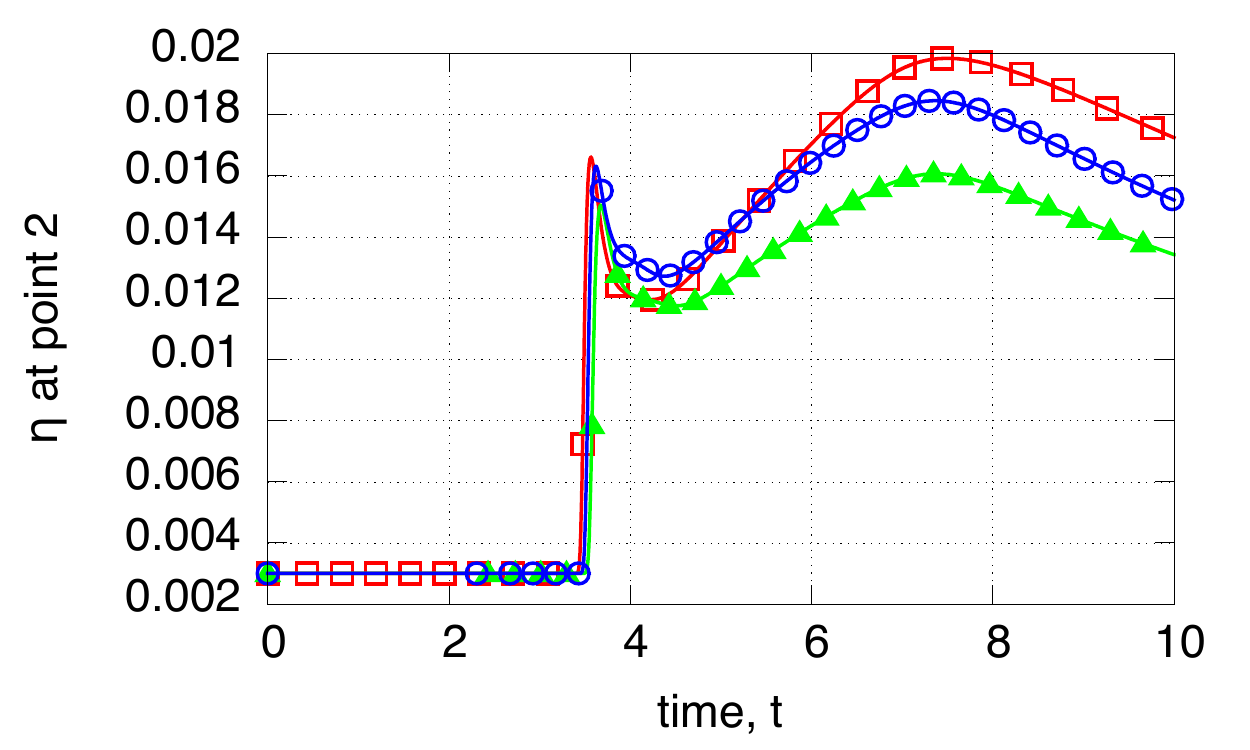}}
	\subfigure{\includegraphics[width=0.48\textwidth]{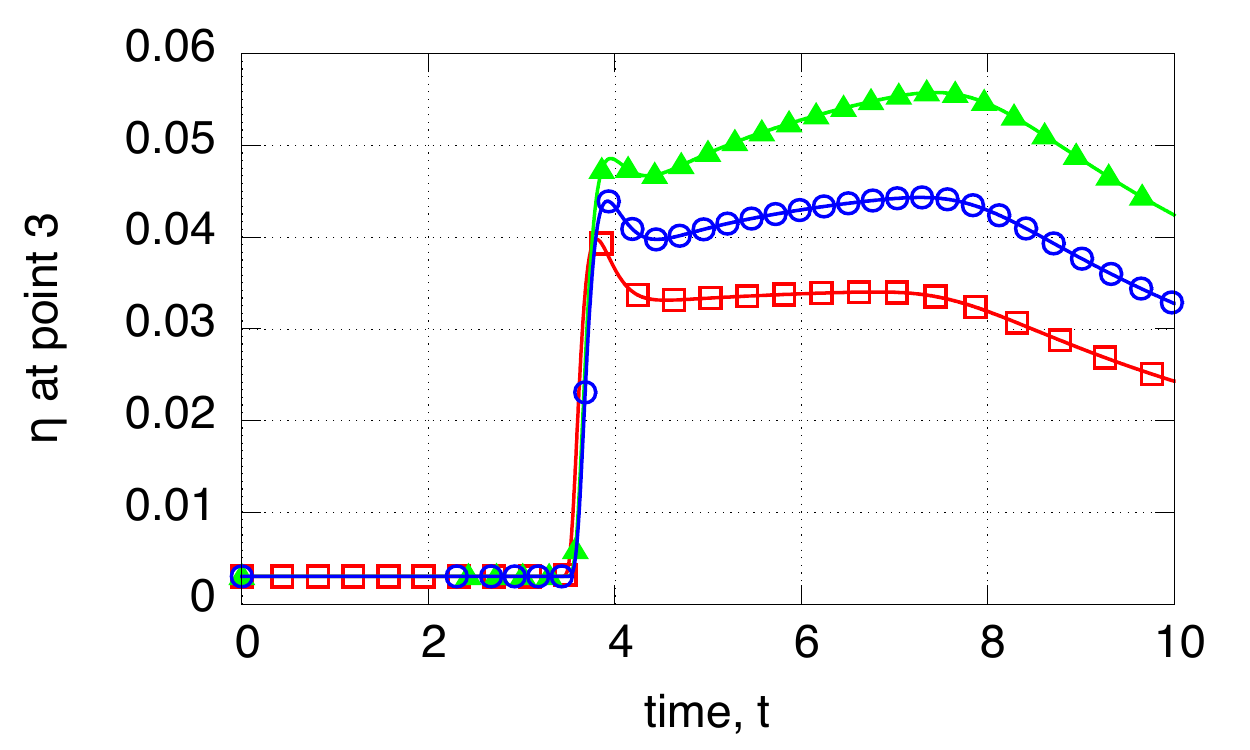}}
    \subfigure{\includegraphics[width=0.48\textwidth]{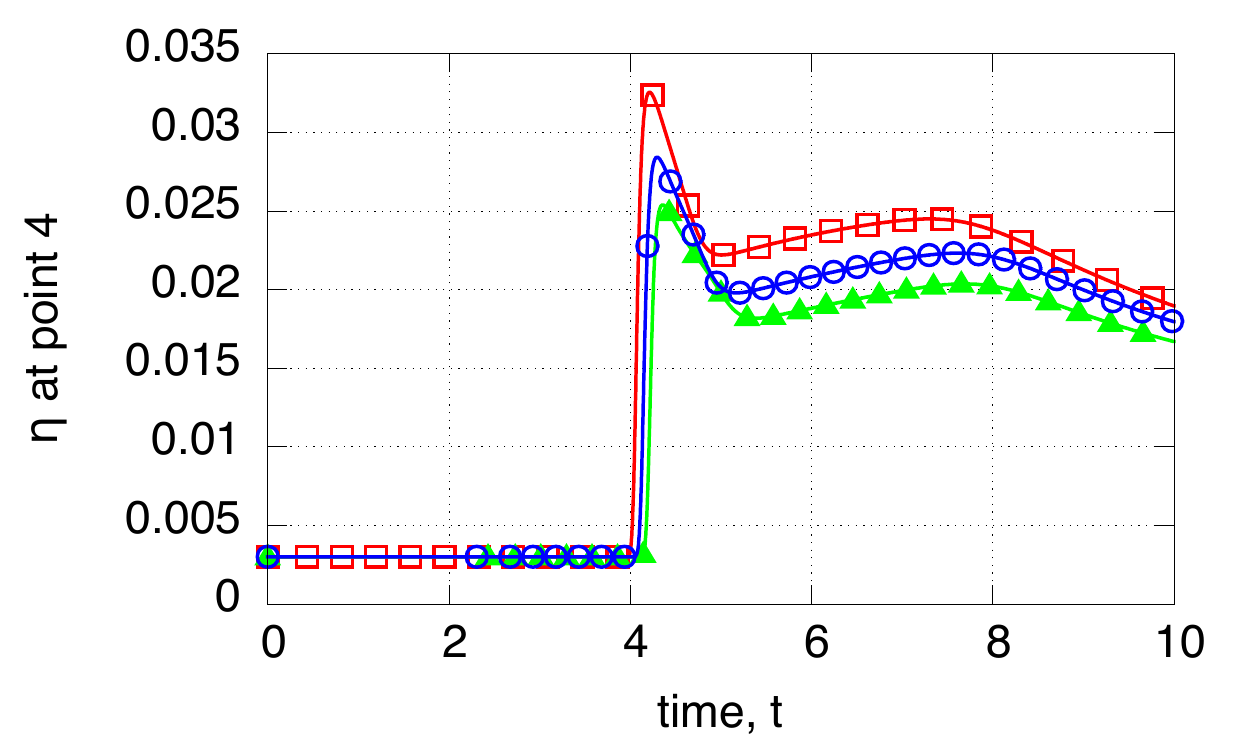}}
	\subfigure{\includegraphics[width=0.48\textwidth]{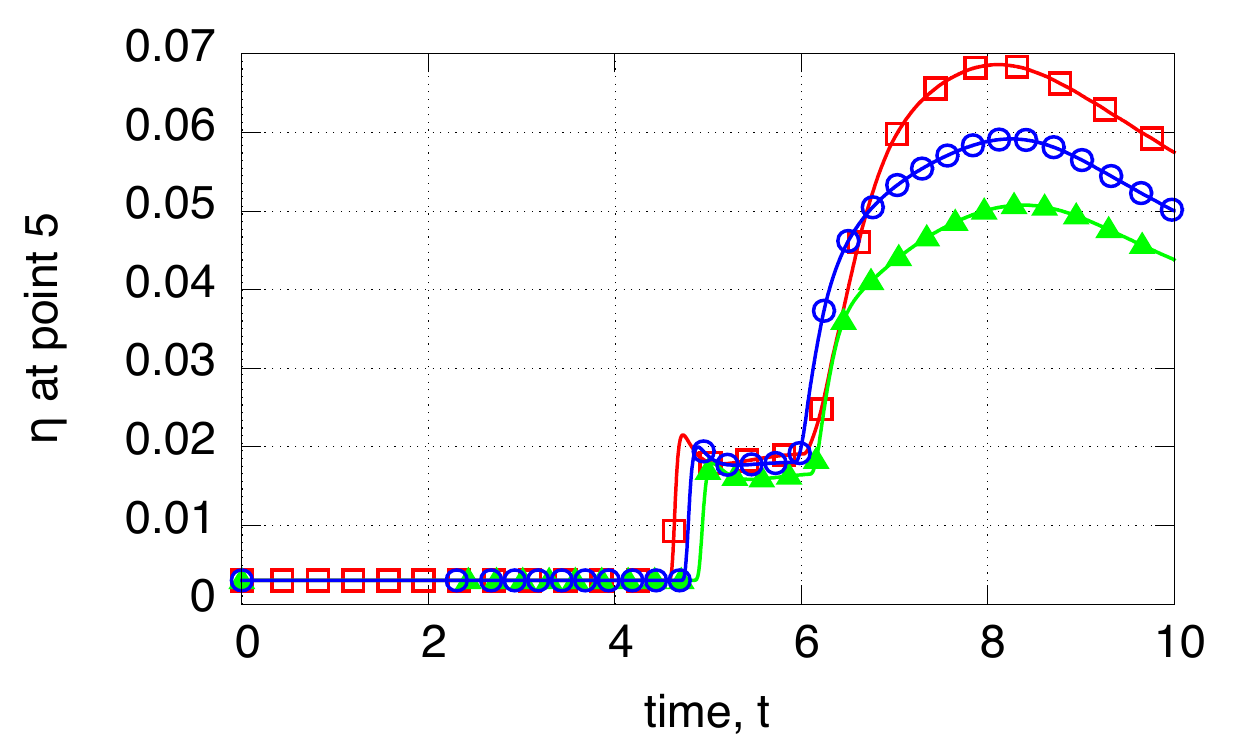}}
	\subfigure{\includegraphics[width=0.48\textwidth]{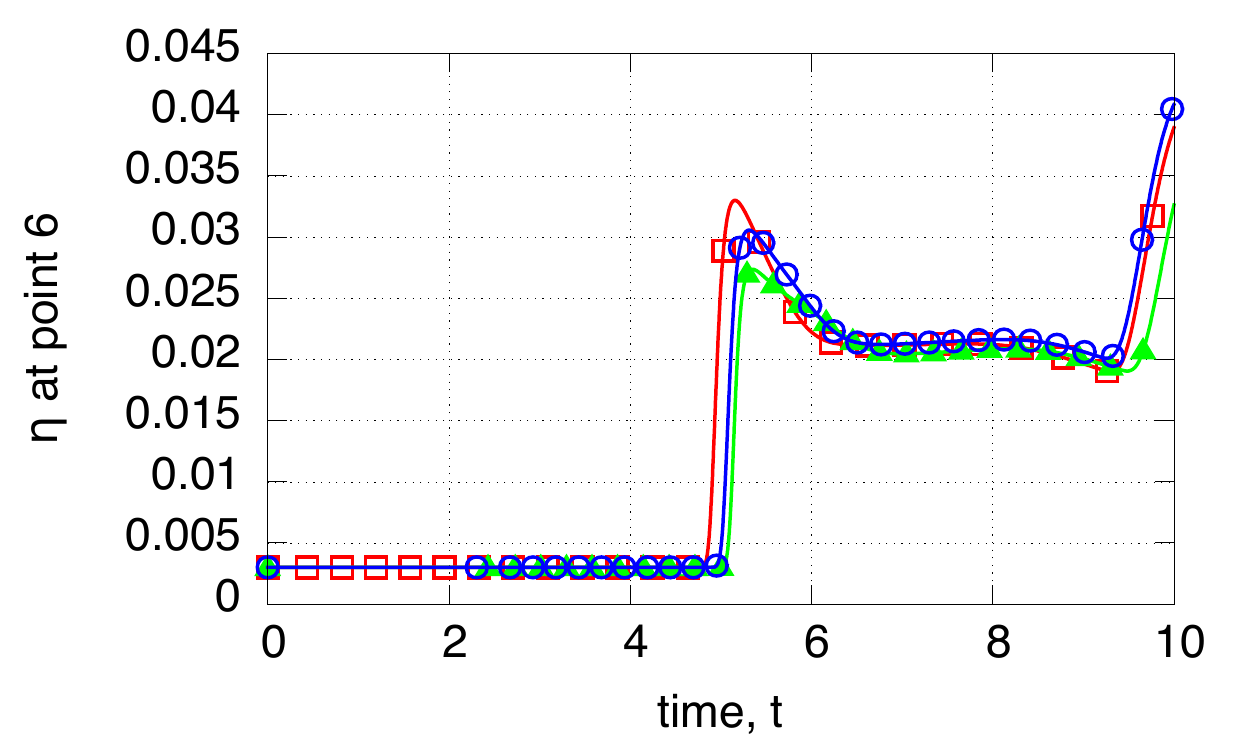}}
	\subfigure{\includegraphics[width=0.5\textwidth]{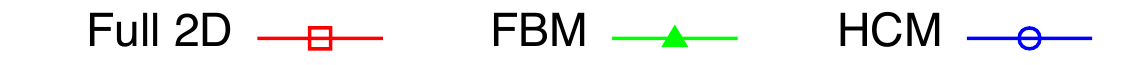}}
  \end{center}
  \caption{Test 1 : Comparison of time evolution of the free surface elevation, $\eta$ at the probe points indicated in figure \ref{fvmnumfigwhycouplinggeometry}.} 
  \label{horinumfigtest1probepoints}
\end{figure}

\begin{figure}[ht!] 
 \subfigure{\includegraphics[width=0.48\textwidth]{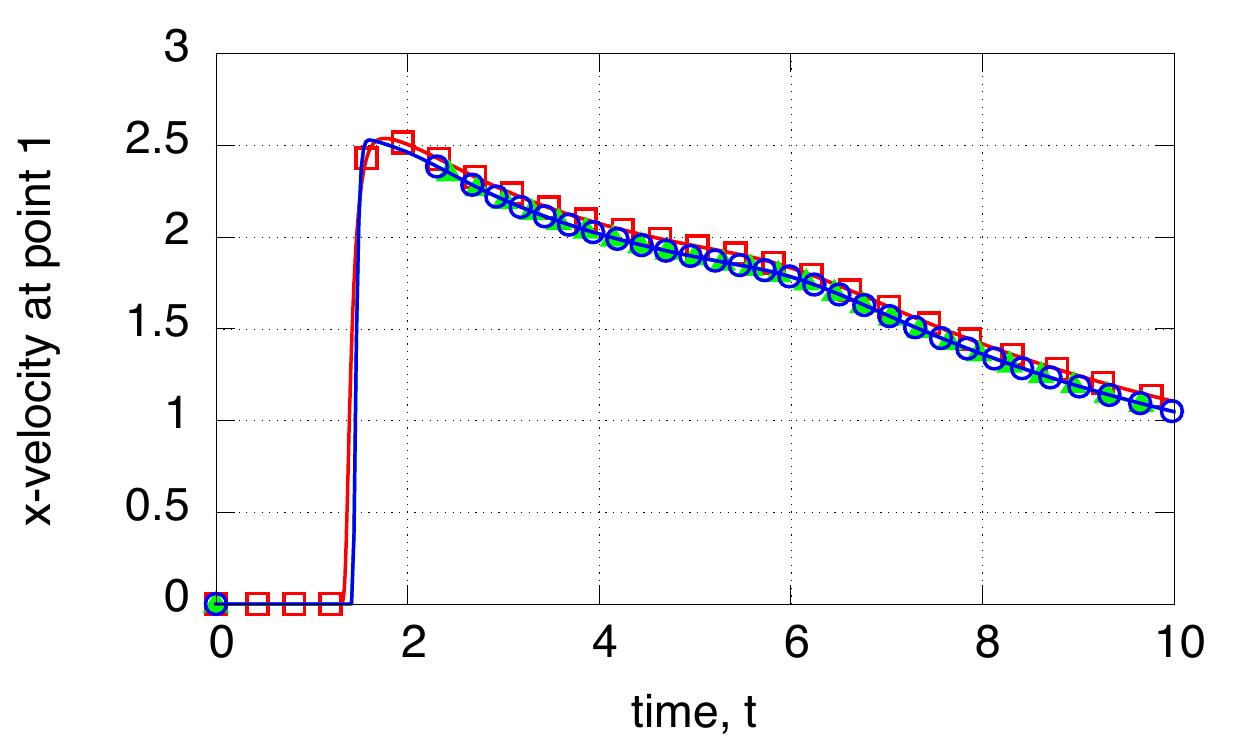}}
 \subfigure{\includegraphics[width=0.48\textwidth]{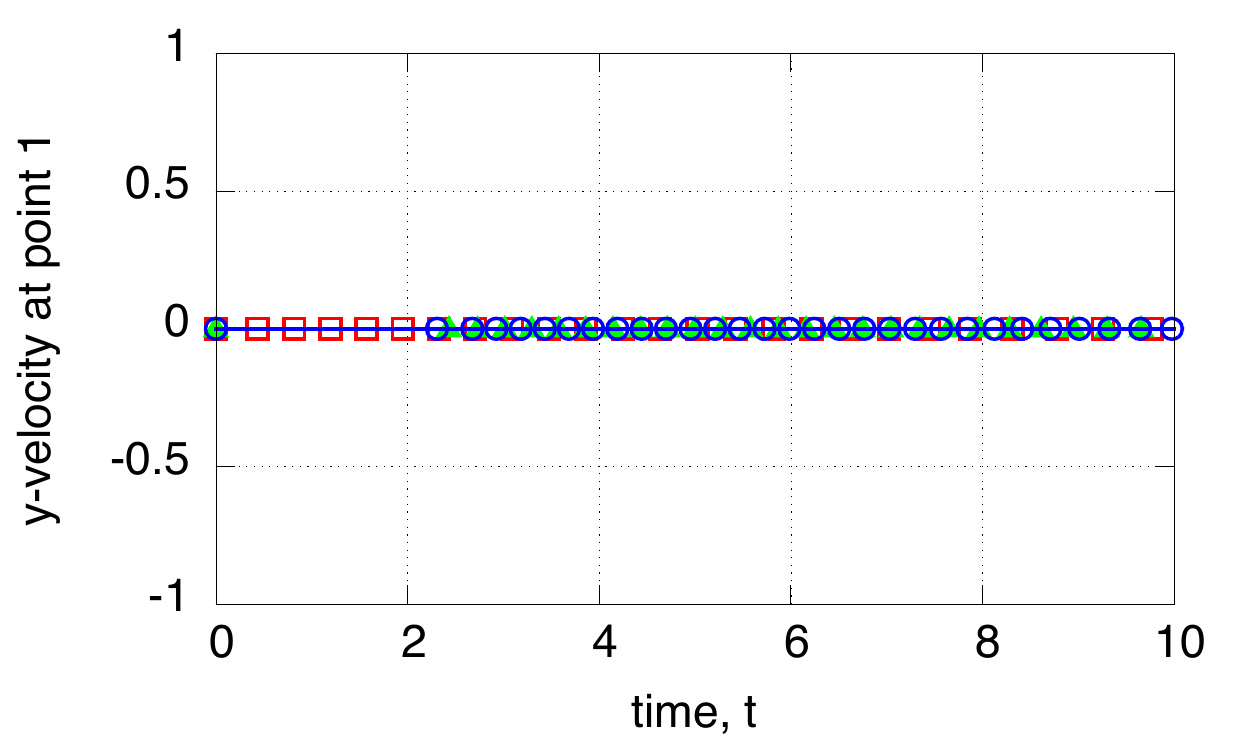}}
 \\
 \subfigure{\includegraphics[width=0.48\textwidth]{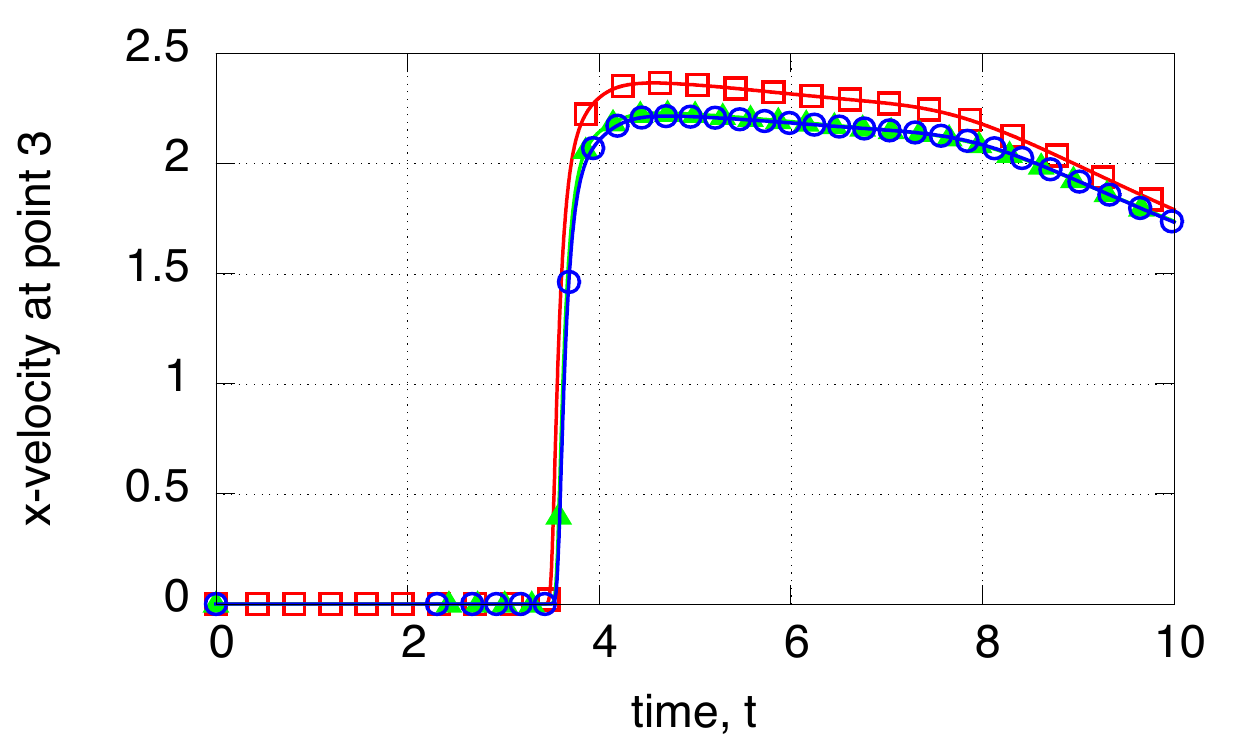}}
 \subfigure{\includegraphics[width=0.48\textwidth]{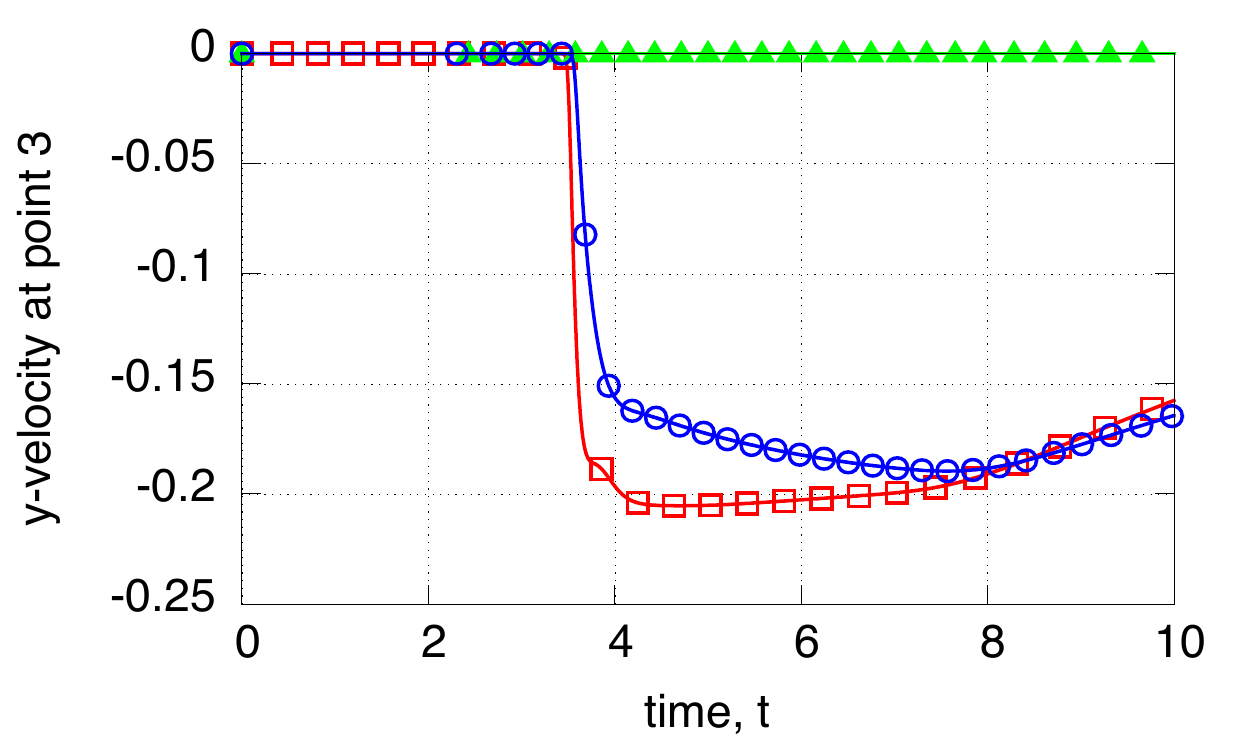}}
 \\
\subfigure{\includegraphics[width=0.48\textwidth]{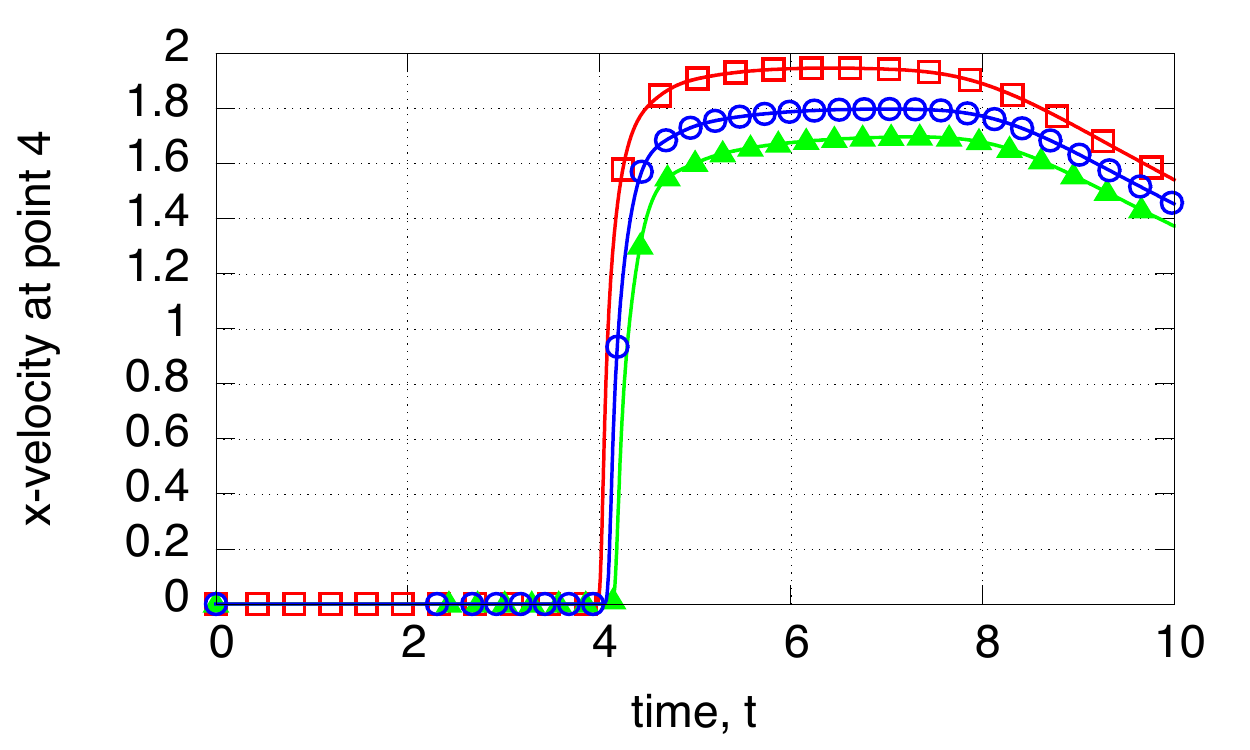}}
\subfigure{\includegraphics[width=0.48\textwidth]{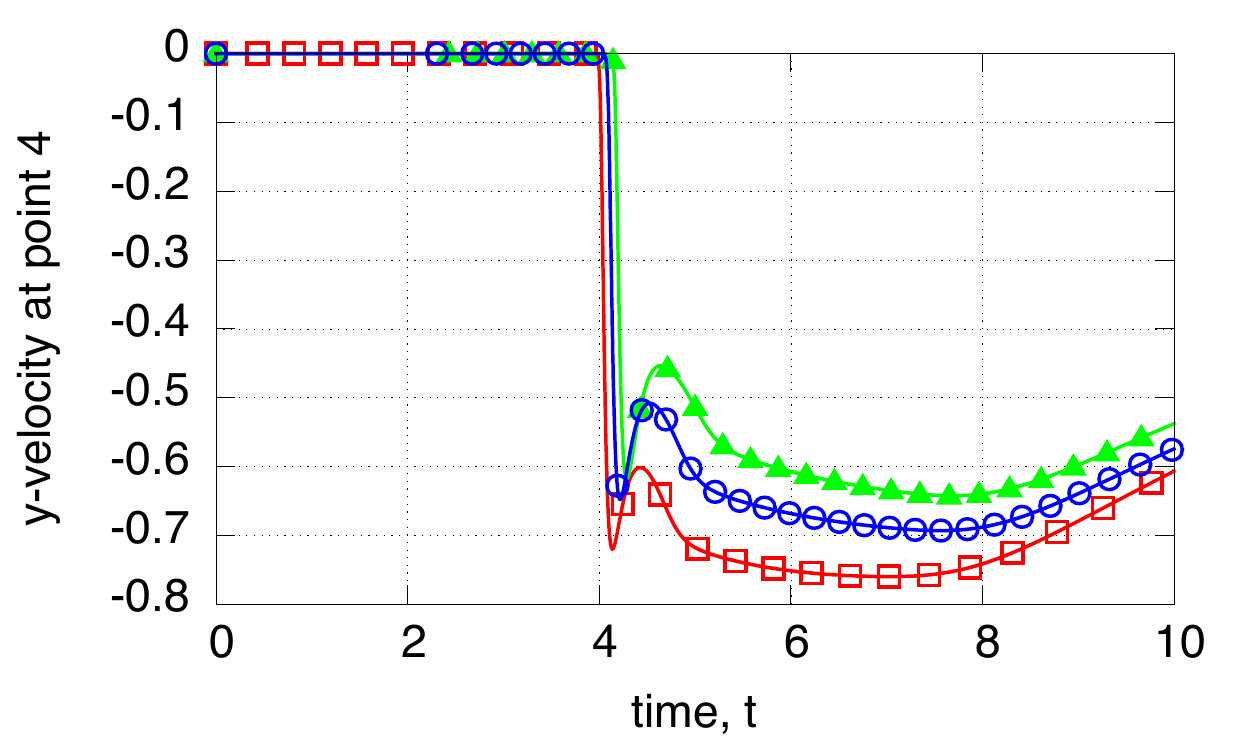}}
\\
\subfigure{\includegraphics[width=0.48\textwidth]{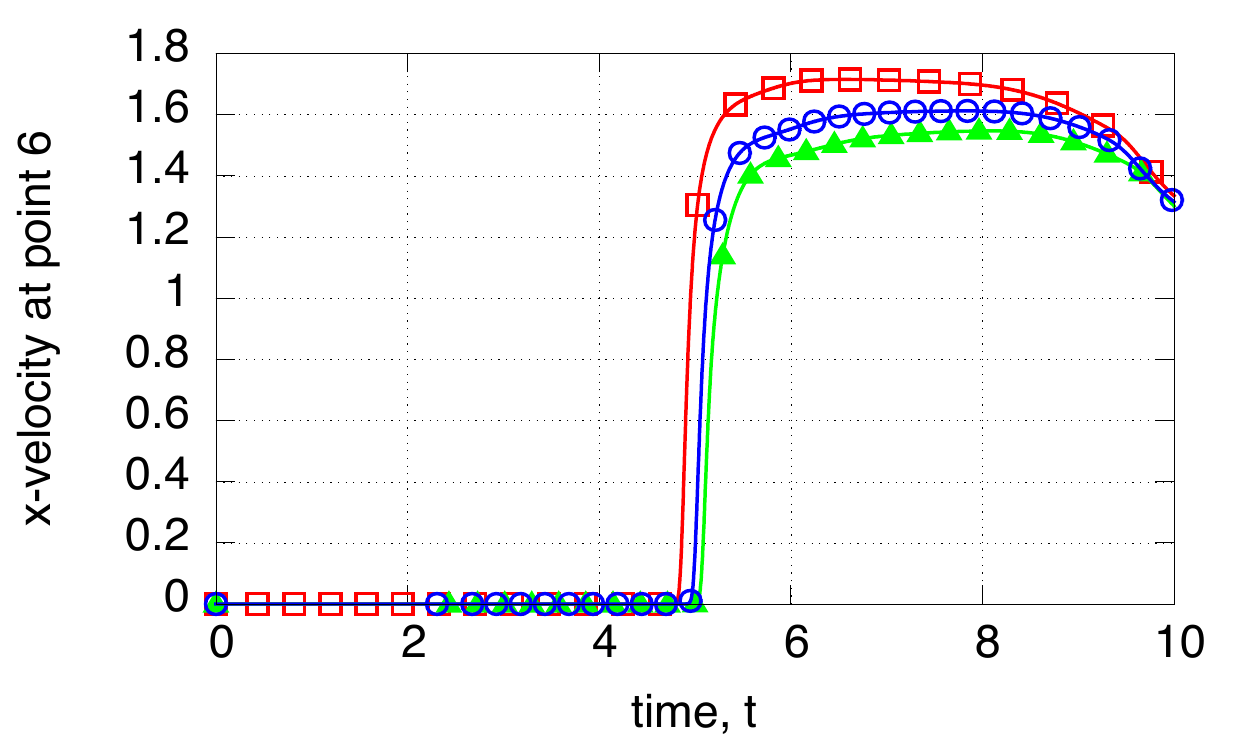}}
\subfigure{\includegraphics[width=0.48\textwidth]{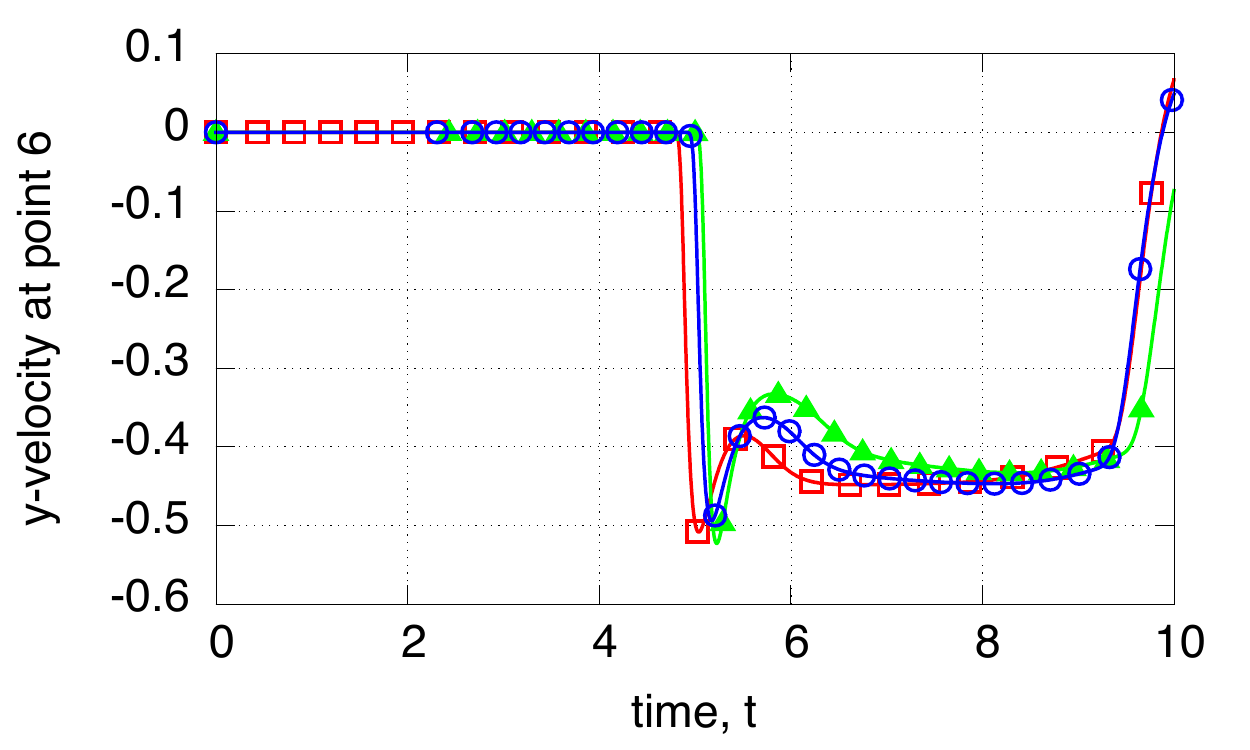}}
\\
  \subfigure{\includegraphics[width=0.5\textwidth]{legend-use}}
  \caption{Time evolution  of the $x$-velocity (left column) and  $y$-velocity (right column) at the indicated selected probe points for test       case 1.}
  \label{numfigtest1-xyvel-floodplainprobepoints}
\end{figure}

Figure \ref{horinumfigtest1eta} displays the free surface elevation for the three simulation methods. We can see that both the proposed
method and the FBM capture the behaviour of the full 2D simulation, and from the right end of the channel, one can also
see  that the proposed method approximates the full 2D result better than the FBM. This is more obvious in figure \ref{horinumfigtest1probepoints},
in which the time evolution of the free surface elevation at the probe points (indicated in figure \ref{fvmnumfigwhycouplinggeometry}), are plotted for all the simulation methods.
Furthermore, figure \ref{numfigtest1-xyvel-floodplainprobepoints} displays the plots of the time evolution of the
$x$- and $y$- velocity components at selected probe points. One can see that the HCM performs better than the FBM
at the indicated points. Of particular interest is the $y$-velocity component at the probe point, $P_3$ which
is located within the channel.
We can see that while the FBM wrongly computed a zero $y$-velocity all the time,
the HCM computed the correct none-zero values with very good accuracy.
It is thus very clear that the proposed method, HCM significantly outperforms the FBM for this test case.

In terms of efficiency, as shown in table
\ref{horinumtabletest1}, the full 2D simulation took 3,669 time steps and 3,100.31 seconds to complete this simulation while the coupling methods took less number of time steps
and more than 50\% reduced time to complete the same simulation. With the above observations, we conclude that it is 
possible to efficiently
use coupling methods instead of full 2D simulation, and that the proposed method is capable of reproducing the full 2D solutions with
greater accuracy than the FBM.

%

\subsection{Test case 2 : Channel Flow into Elevated 2D Floodplain}\label{numsectest2}
\begin{figure}[ht!] 
  \begin{center}
	\includegraphics[width=0.98\textwidth]{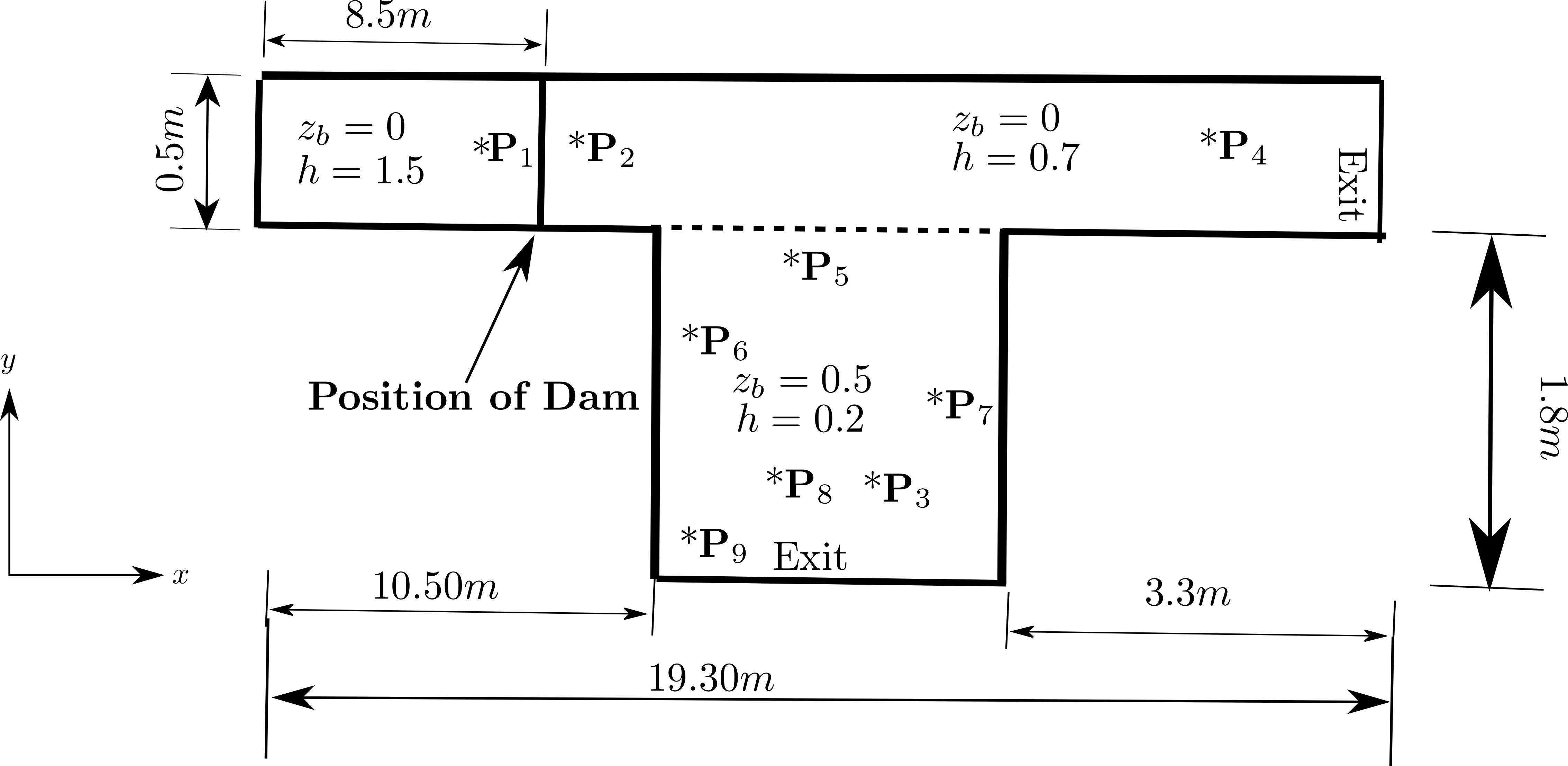}
  \end{center}
   \caption{ Top view of Channel and Floodplain for test case 2 showing the floodplain region in $(x,y)\in[10.5,16.0]\times[0,1.8]$   and the channel region in $(x,y)\in[0,19.3]\times[1.8,2.3]$.}
   \label{numfiggeometrytest2}
\end{figure}

This test case involves the same channel as in the previous example but connected to an elevated floodplain located in the
region $10.5\leq x \leq 16.0$ (see figure \ref{numfiggeometrytest2}). The channel bed is flat and the floodplain bed is 0.5 meters high.
The initial condition is the following.
\begin{align}
 & H(x,y,0)  = \begin{cases} 1.5, & \mbox{ if } x \leq 8.5, \quad y \geq 1.8, \\
                            0.7, & \mbox{ if } x > 8.5, \quad y \geq 1.8, \\
                            0.2, & \mbox{ if } 10.5\leq x \leq 16.0, \quad 0 \leq y \leq 1.8,
               \end{cases} \\
& u(x,y,0) = v(x,y,0) = 0.               
\end{align}
The manning coefficient for both channel and floodplain is taken as 0.009$s/m^{1/3}$ the boundaries
are only open at the sides indicated "exit" in figure \ref{numfiggeometrytest2}, others are closed.
Just like the previous test case, here nine probe points, $P_1-P_9$ are identified, see figure
\ref{numfiggeometrytest2}.
%

\begin{figure}[ht!] 
	\includegraphics[width=0.98\textwidth]{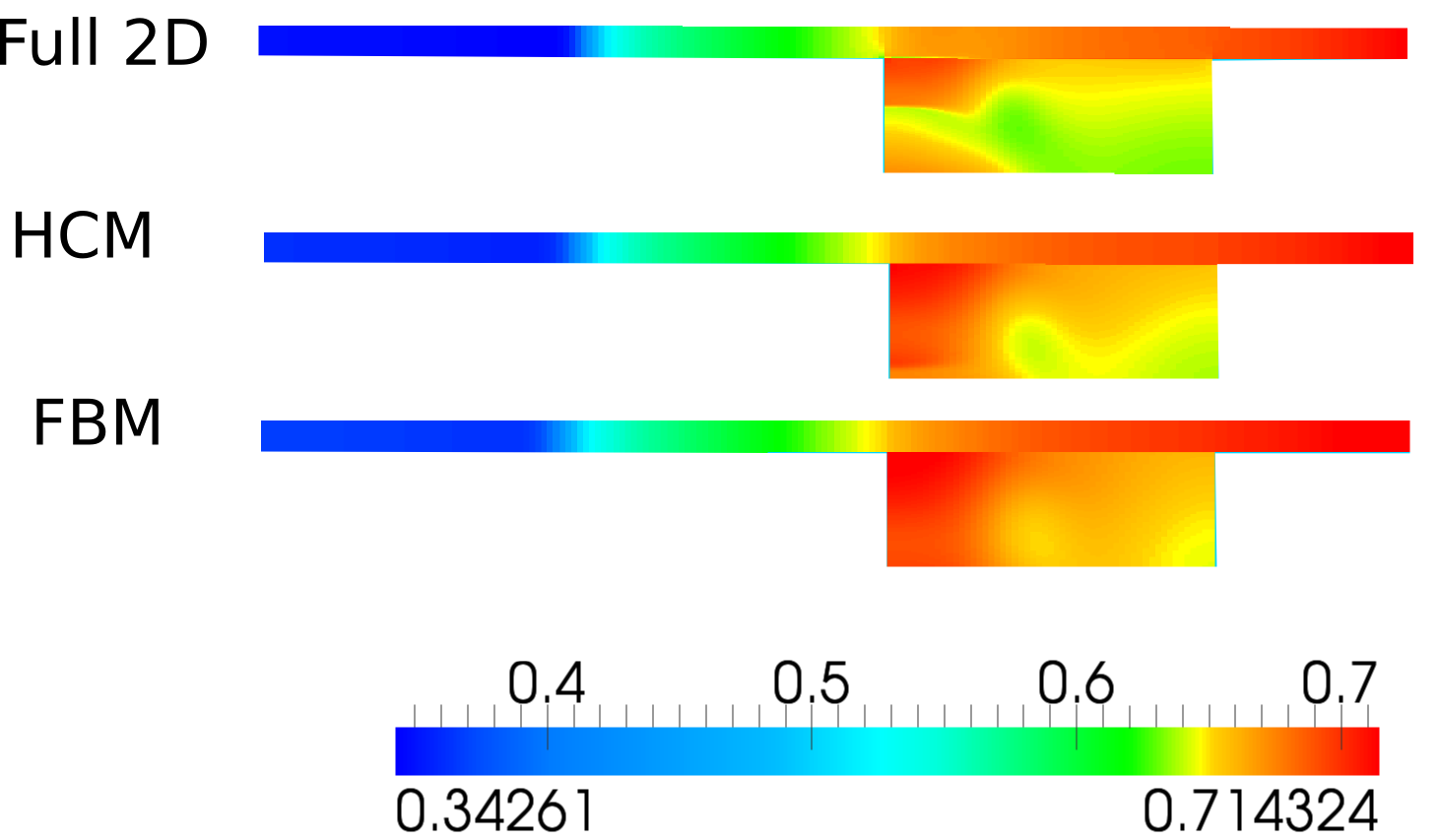}
  \caption{Comparison of free surface elevation for the different methods after the last time step: Test 2 } 
  \label{numfigtest1etaview}
\end{figure}

\begin{figure}[ht!] 
	\includegraphics[width=0.98\textwidth]{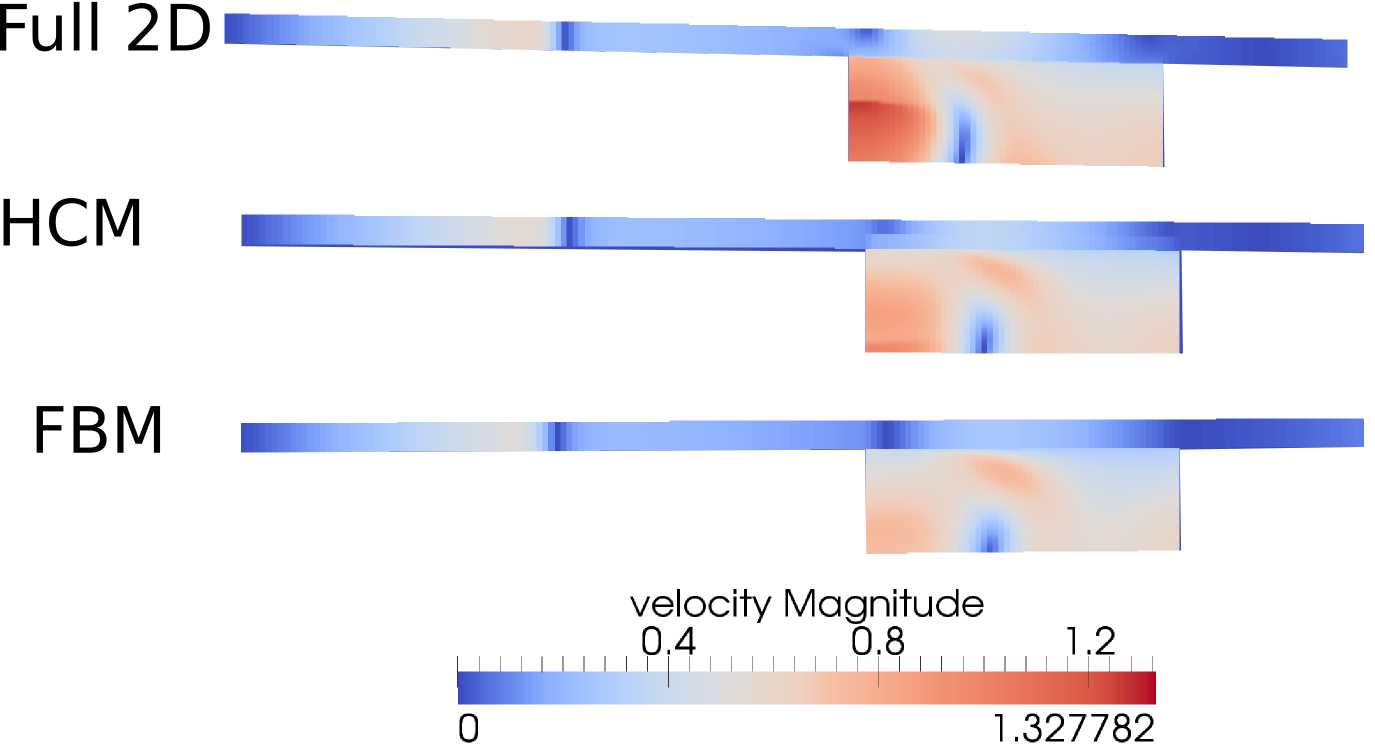}
  \caption{Comparison of velocity magnitude for the different methods after the last time step: Test 2} 
  \label{numfigtest2velomagview}
\end{figure}

\begin{figure}[ht!] 
	\subfigure{\includegraphics[width=0.48\textwidth]{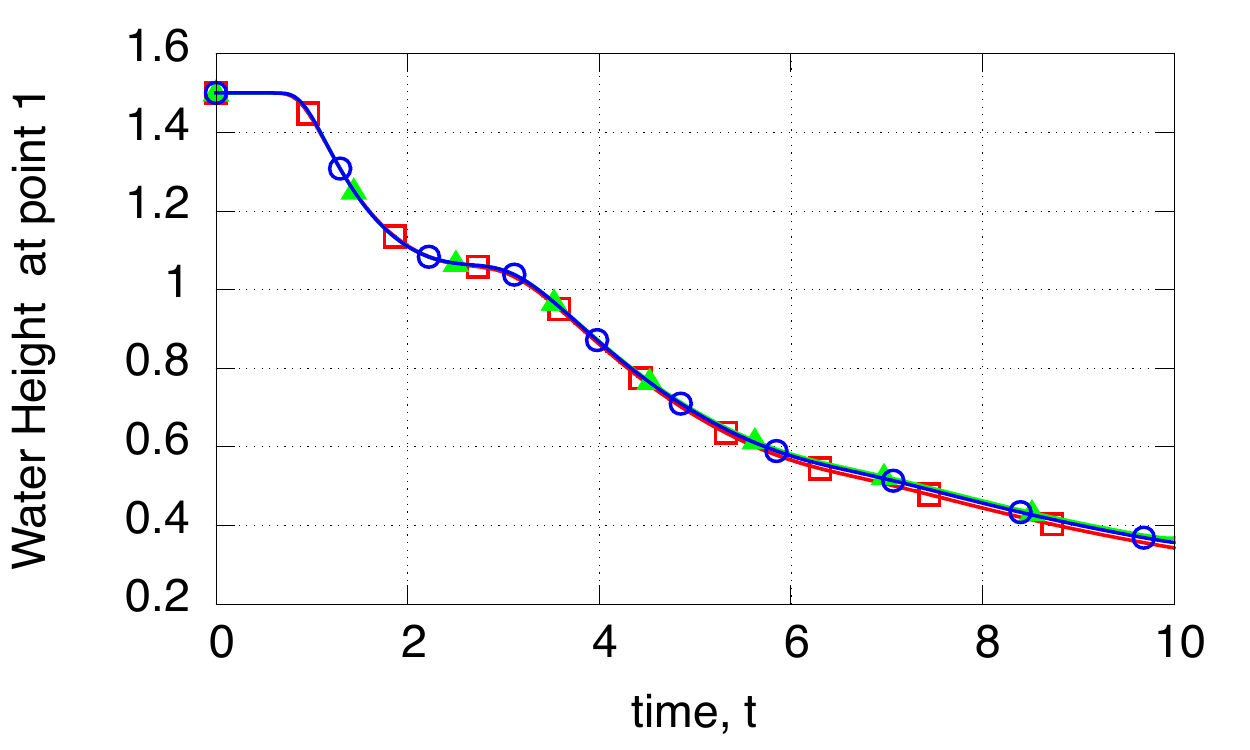}}
	\subfigure{\includegraphics[width=0.48\textwidth]{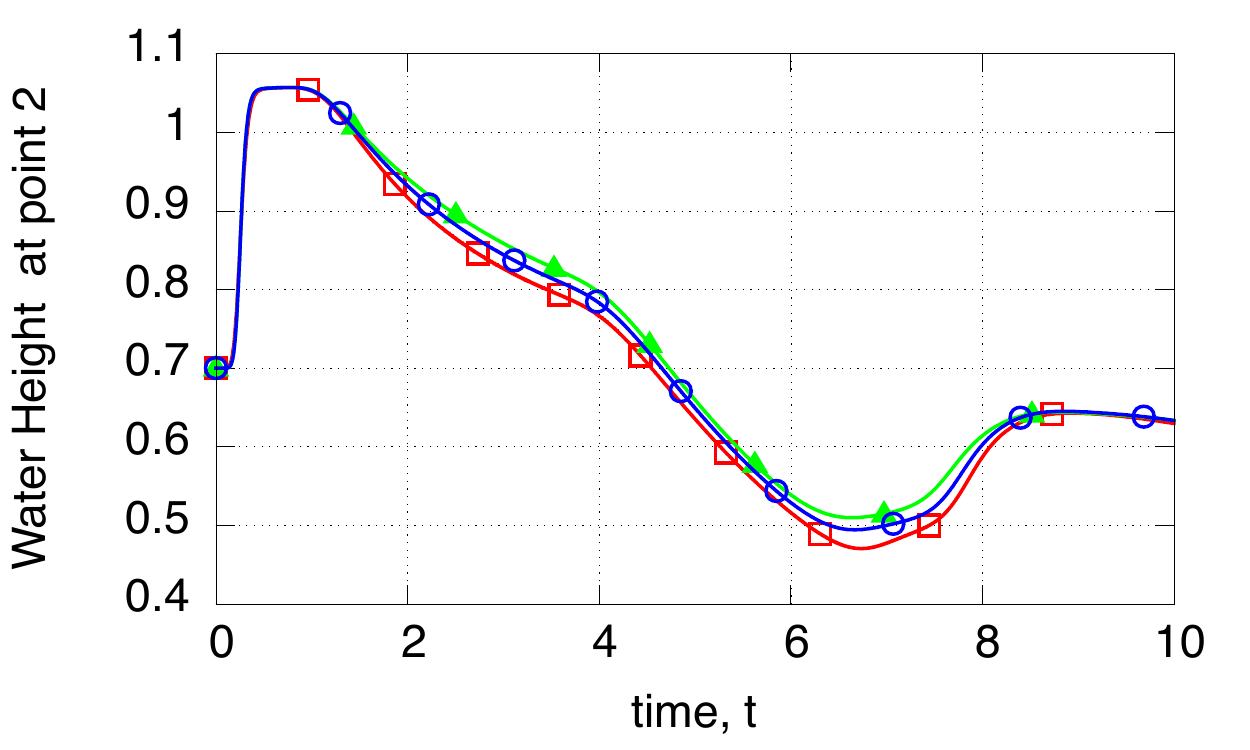}}
	\subfigure{\includegraphics[width=0.48\textwidth]{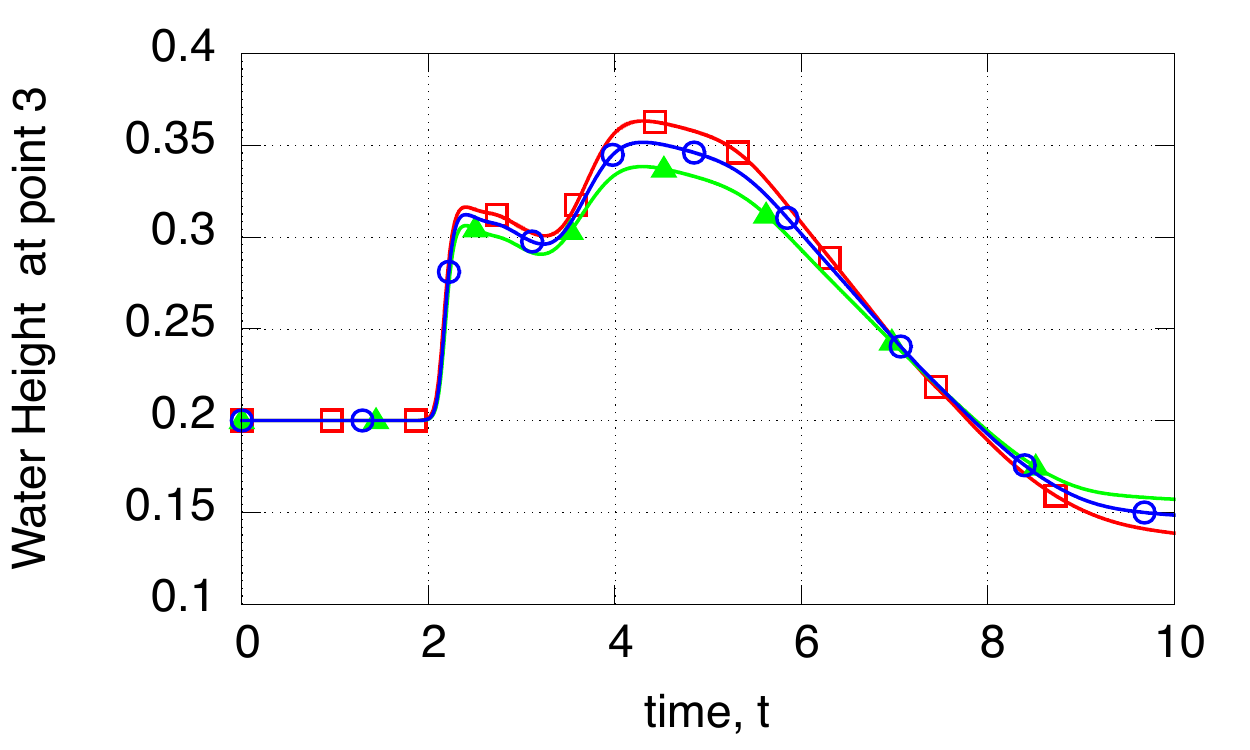}}
	\subfigure{\includegraphics[width=0.48\textwidth]{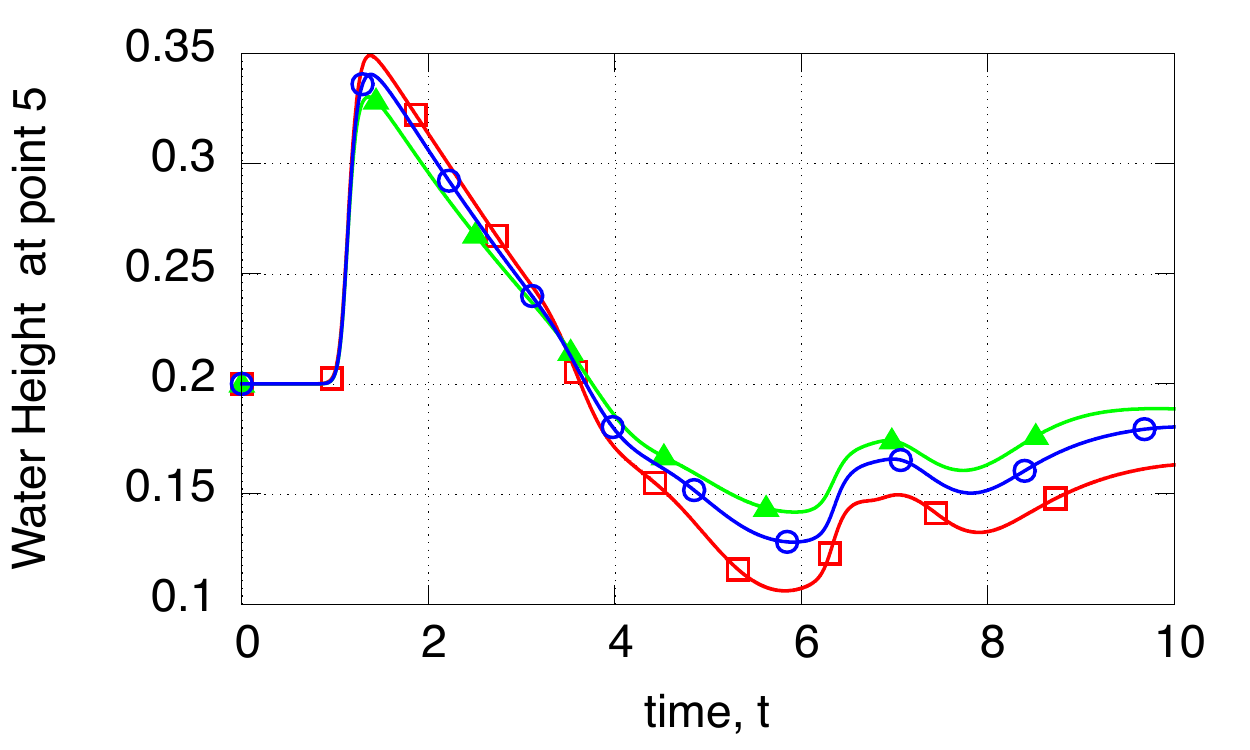}}
	\subfigure{\includegraphics[width=0.48\textwidth]{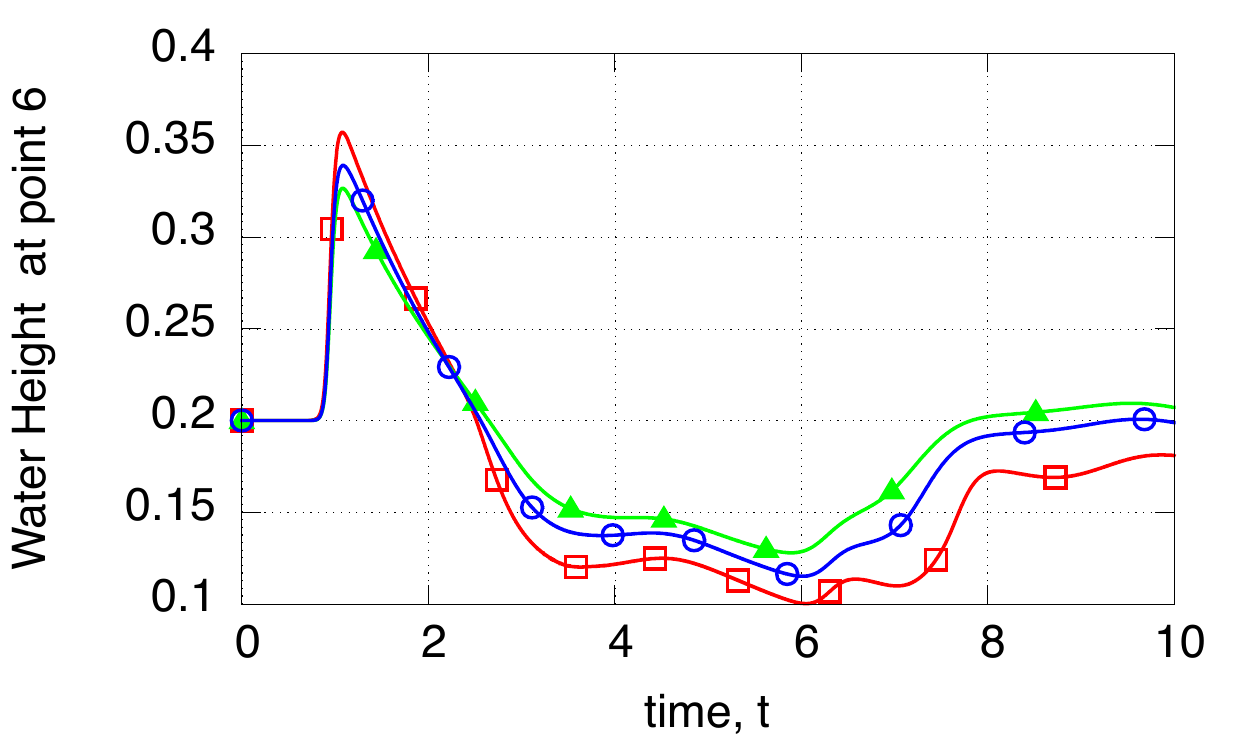}}
	\subfigure{\includegraphics[width=0.48\textwidth]{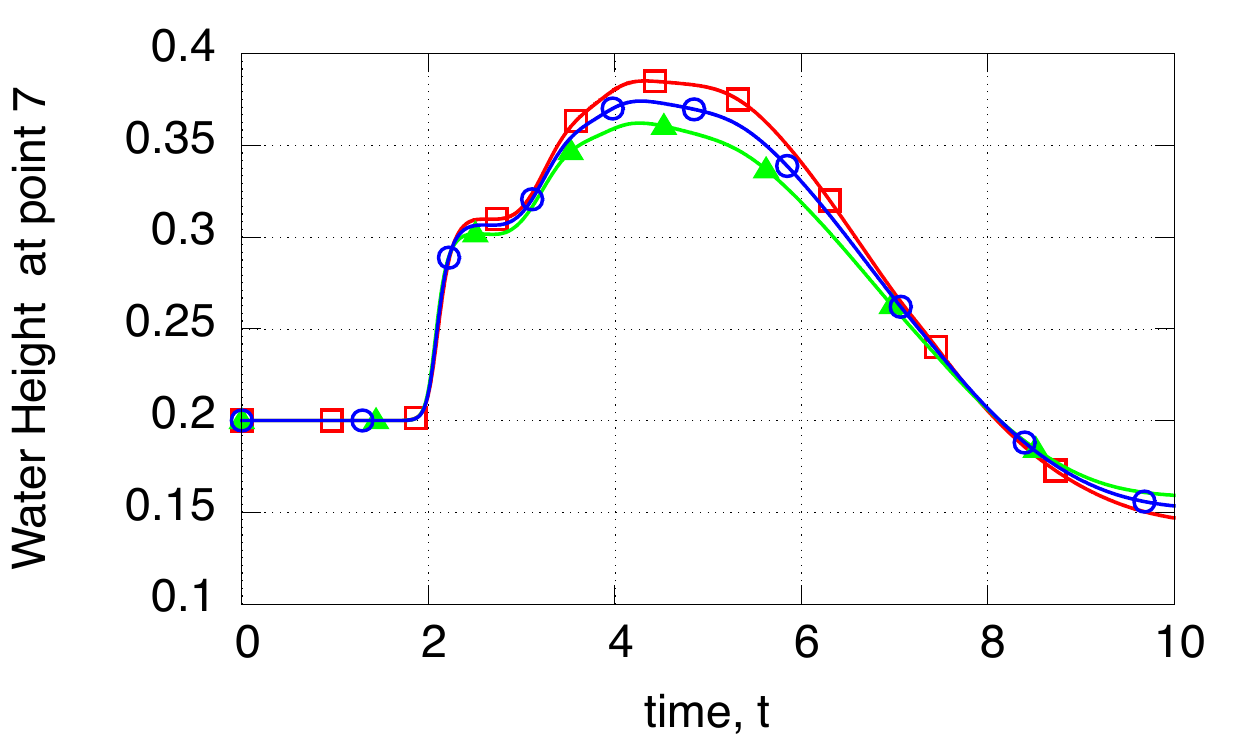}}
	\subfigure{\includegraphics[width=0.48\textwidth]{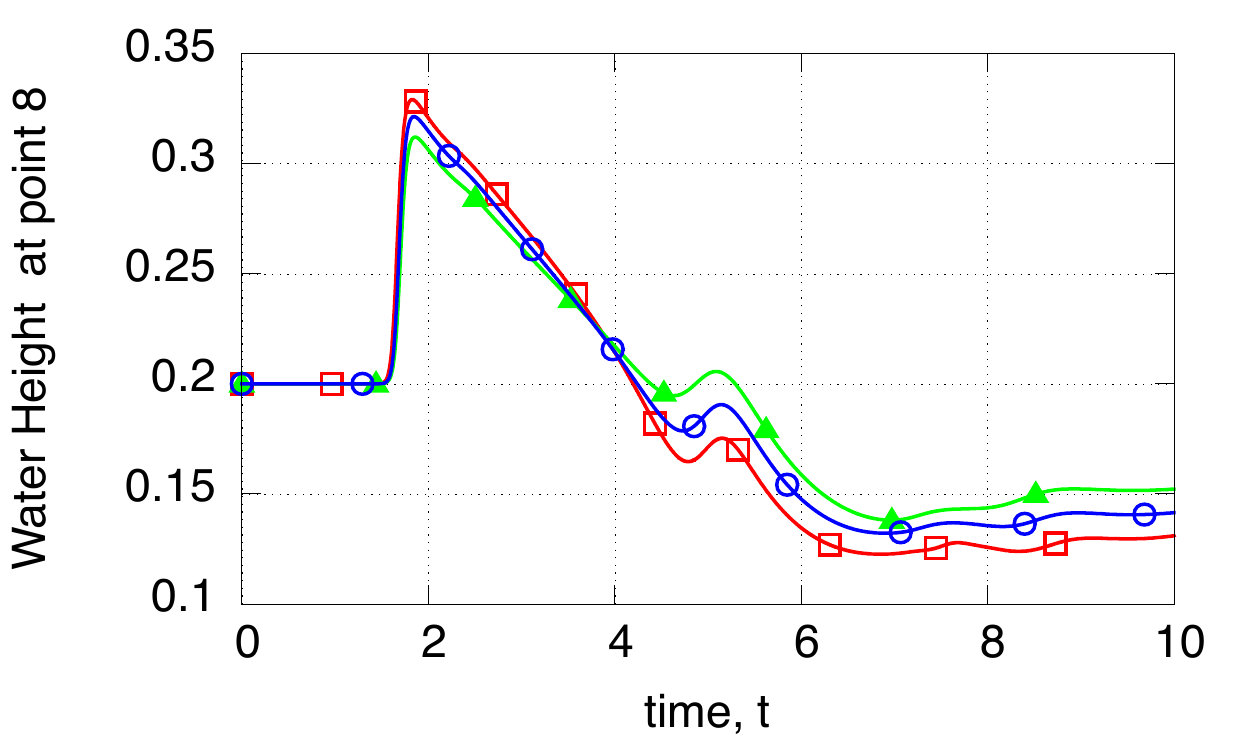}}
	\subfigure{\includegraphics[width=0.48\textwidth]{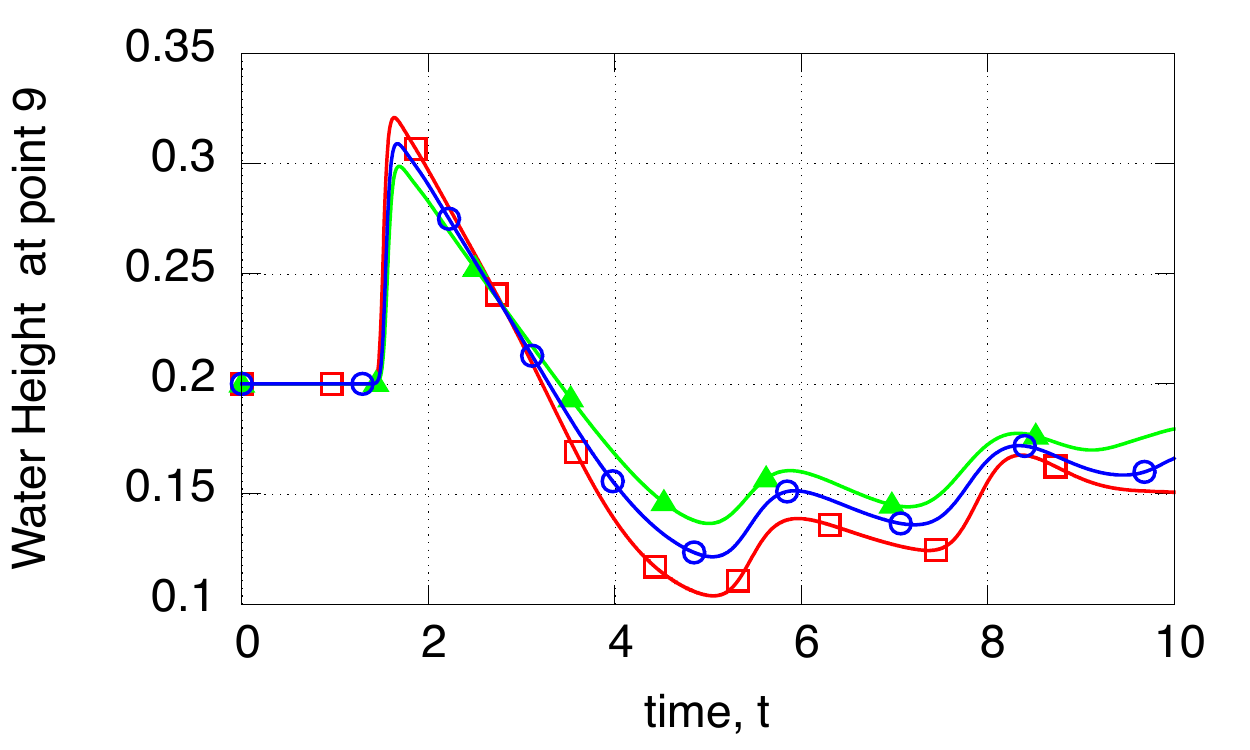}}\\
	\subfigure{\includegraphics[width=0.50\textwidth]{legend-use}}
  \caption{Comparison of time evolution of water height at probe points : Test 2} 
  \label{numfigtest2etaprobepoints}
\end{figure}

\begin{figure}[ht!] 
 \subfigure{\includegraphics[width=0.48\textwidth]{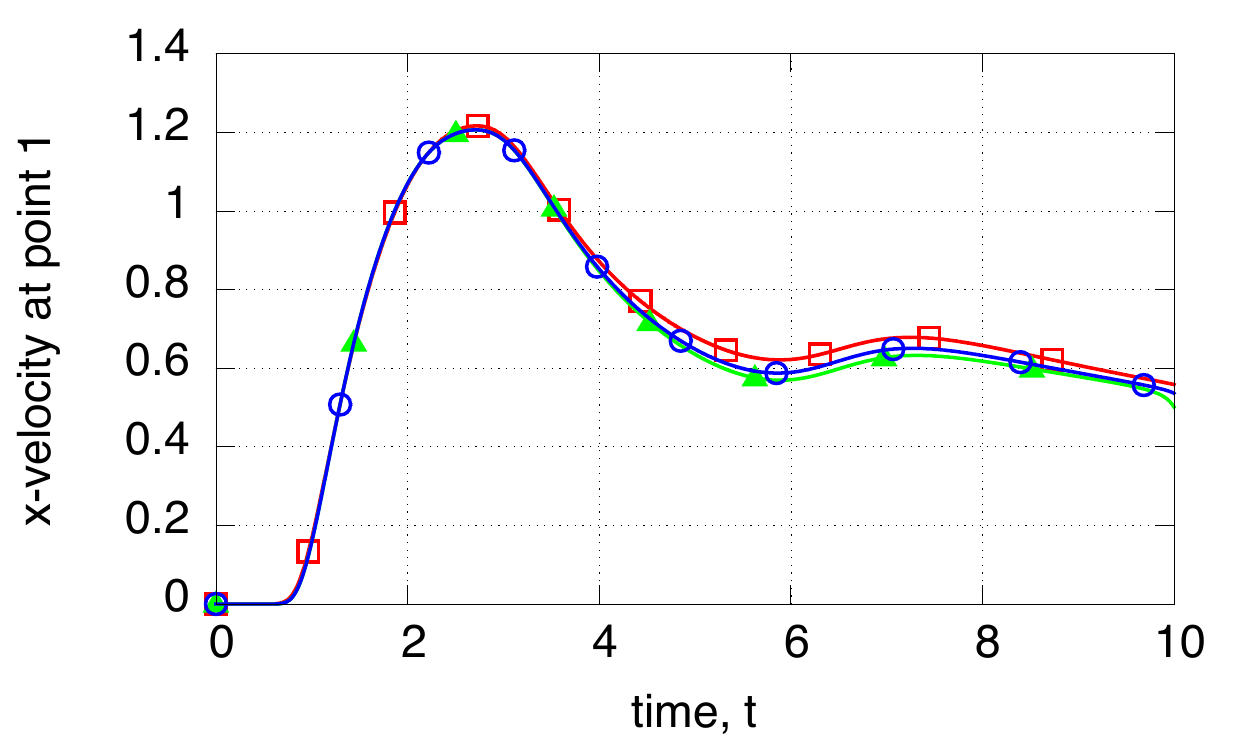}}
 \subfigure{\includegraphics[width=0.48\textwidth]{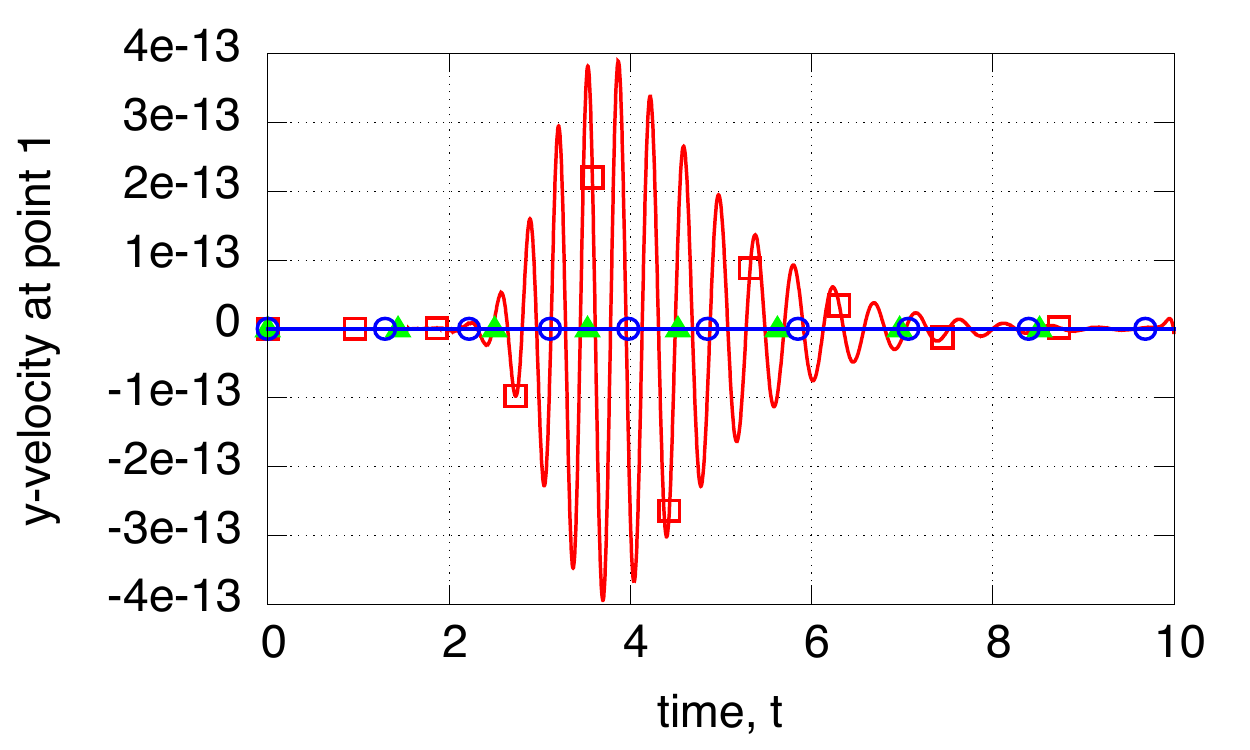}}
 \\
 \subfigure{\includegraphics[width=0.48\textwidth]{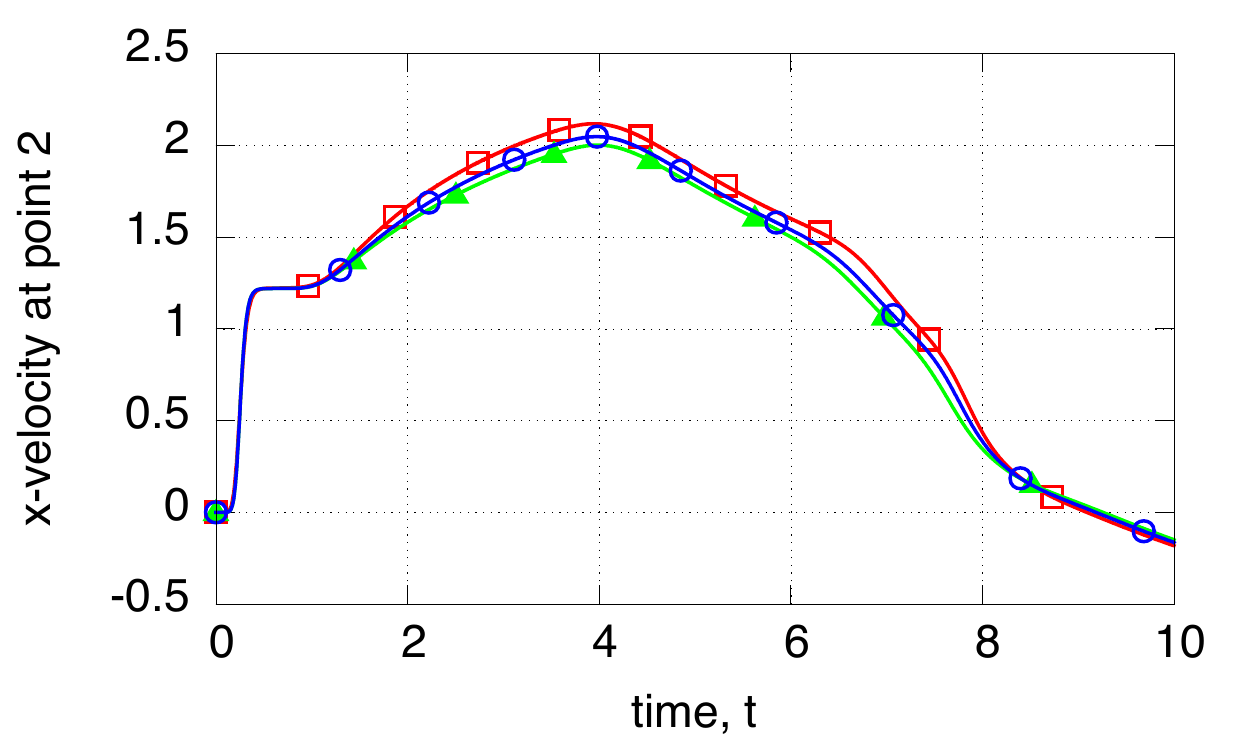}}
 \subfigure{\includegraphics[width=0.48\textwidth]{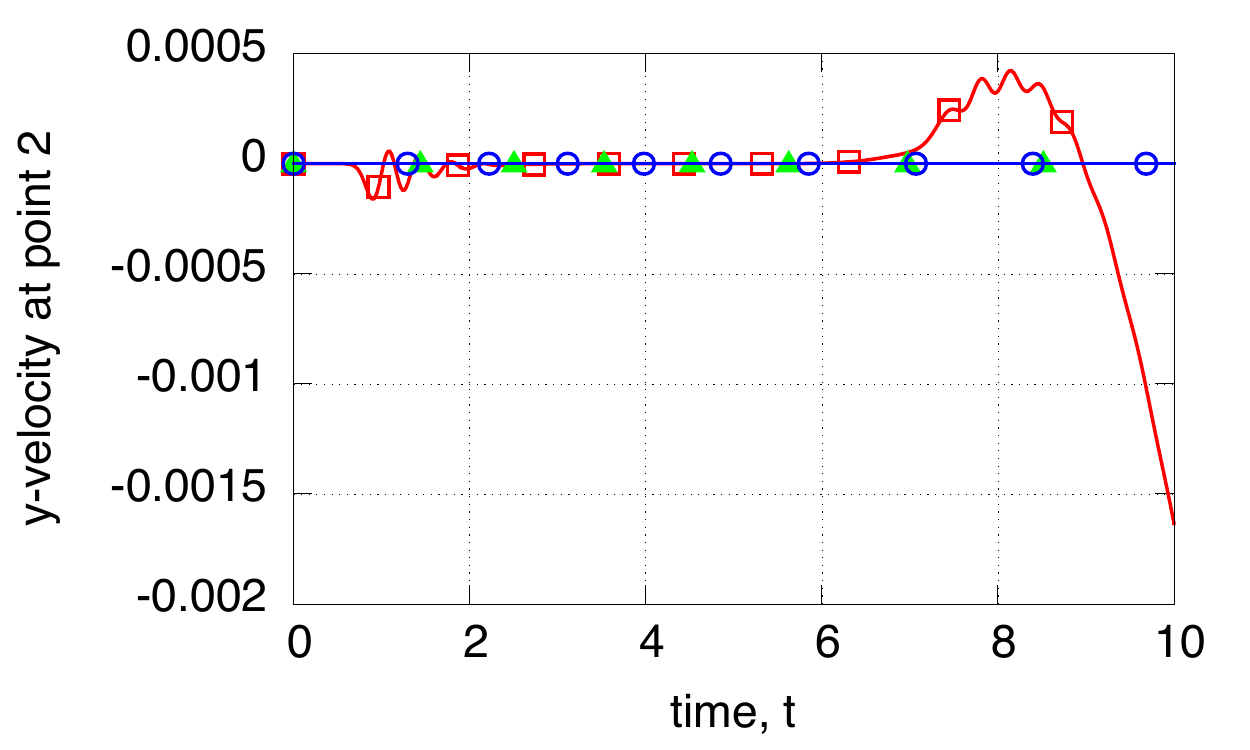}}
 \\
\subfigure{\includegraphics[width=0.48\textwidth]{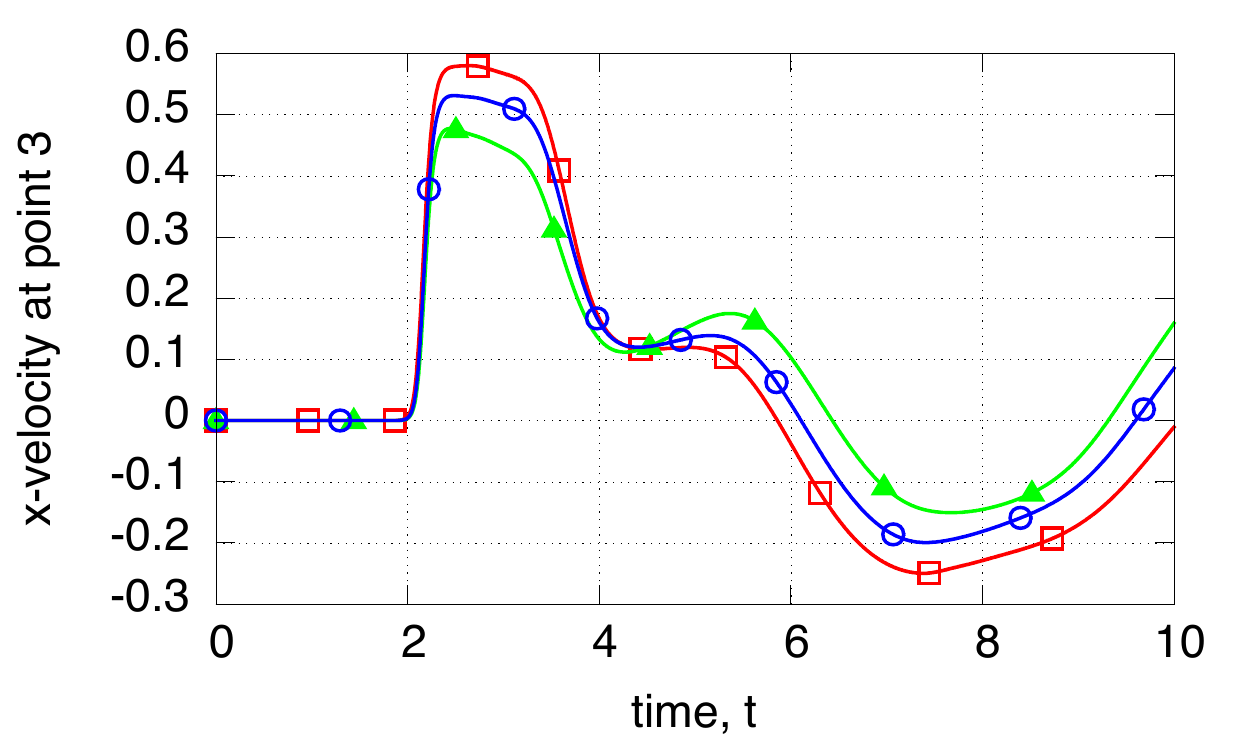}}
\subfigure{\includegraphics[width=0.48\textwidth]{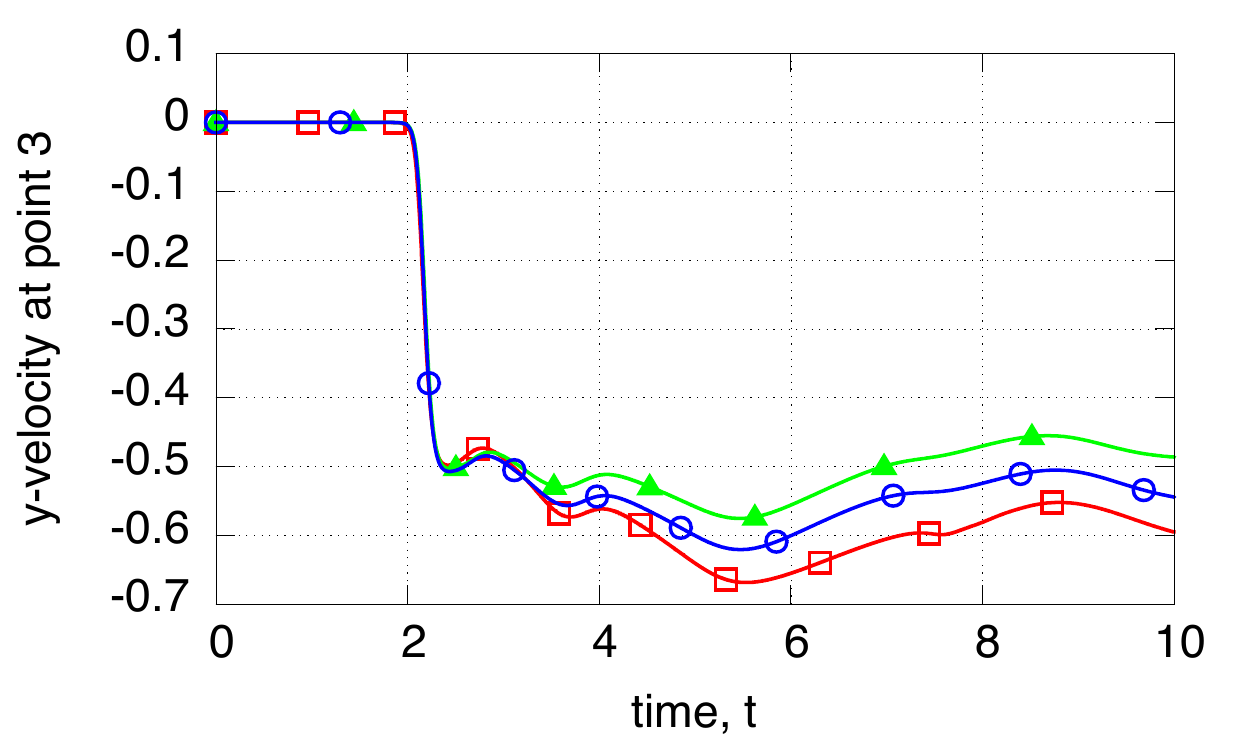}}
\\
\subfigure{\includegraphics[width=0.48\textwidth]{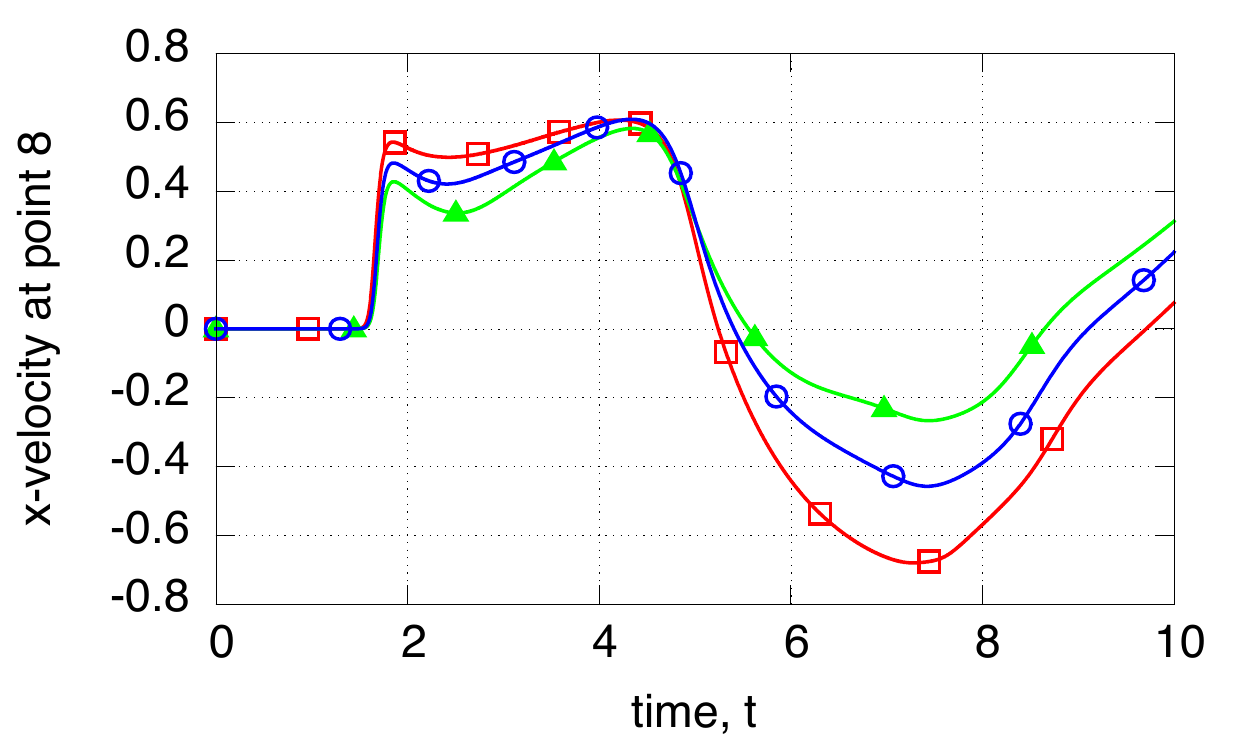}}
\subfigure{\includegraphics[width=0.48\textwidth]{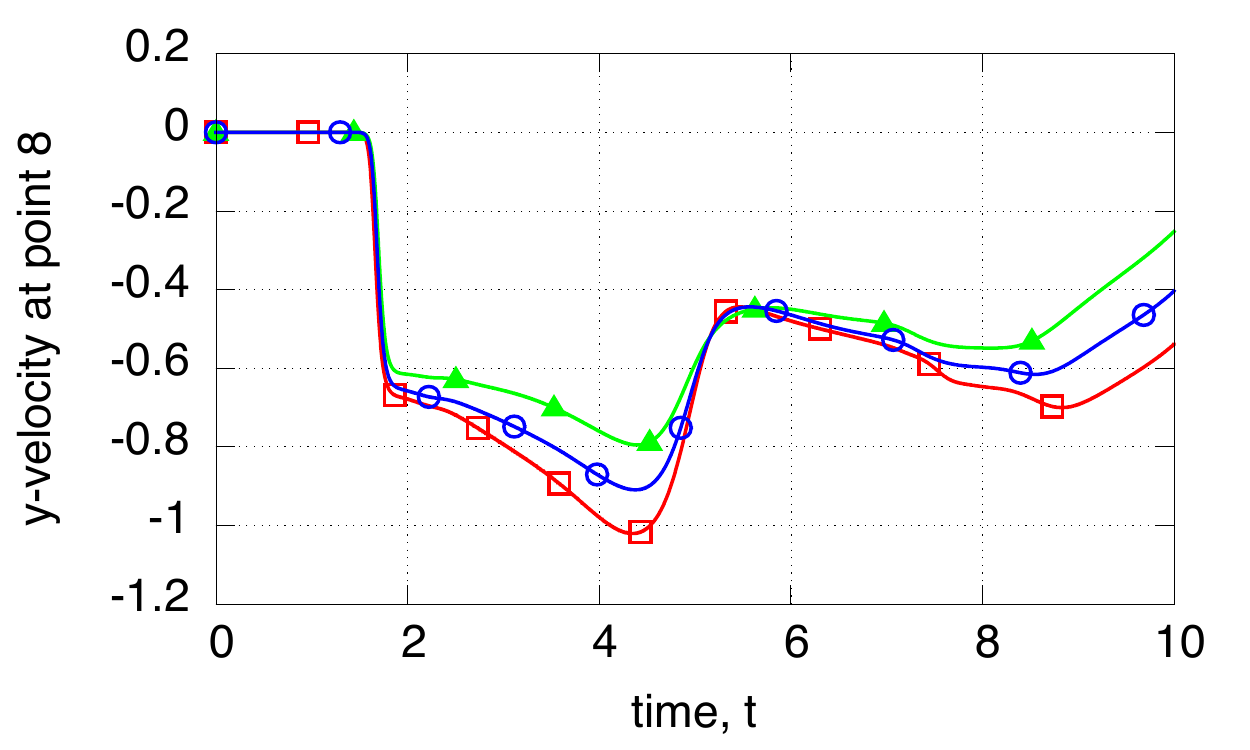}}
\\
  \subfigure{\includegraphics[width=0.5\textwidth]{legend-use}}
  \caption{Time evolution of  $x$-velocity (left column) and  $y$-velocity (right column) at the indicated selected probe points for test       case 2.}
  \label{numfigtest2-xyvel-floodplainprobepoints}
\end{figure}

All the methods solved this problem with a grid of $55\times90$ cells in the floodplain while the channel 
consists of $193\times25$,  $193\times2$ and $193\times1$ cells for the full 2D, the HCM and the FBM
respectively, see table \ref{numtabletest2}. 
The simulation was run for ten seconds.
Figures \ref{numfigtest1etaview} and \ref{numfigtest2velomagview} show the free surface elevation
and velocity magnitude after the last time step for each method. It can be seen that the HCM provides a
better approximation of the full 2D results than the FBM.
As a further validation of this claim, the time evolution of the free surface elevation
is plotted in figures \ref{numfigtest2etaprobepoints}, while those of the 
$x$-velocity and $y$-velocity components are plotted in figure \ref{numfigtest2-xyvel-floodplainprobepoints}
for selected probe points. It can be seen that the horizontal coupling
method is more accurate than the FBM at the points for all flow quantities and almost
all the time. Again, the HCM really captures the flow structure of the full 2D simulation.
This proves the accuracy of the proposed methods over the FBM for this test case.

\begin{table}[ht!] 

\begin{center}
\begin{tabular}{ | l |  c | c |  c | c |}
\hline
                 &      Channel Grid     &   Floodplain Grid    & No. of time steps   &  Processor time (in seconds)           \\
\hline
Full 2D          &     $193\times25$    &    $55\times90$   &       4,963             &     4,100.47  \\ 
\hline
HCM              &      $193\times2$    &    $55\times90$   &  3,235              &    1,710.36    \\ 
\hline   
FBM             &       $193\times1$    &    $55\times90$  & 3,178              &    1,555.08   \\ 
\hline 
\end{tabular}
\end{center}
\caption{Grid cells, simulation times and number of time steps : Test 2}
\label{numtabletest2}
\end{table}

For efficiency, we see from the processing time in table \ref{numtabletest2}
that the FBM is very efficient but not very accurate while the horizontal coupling 
method is both very efficient and also has good accuracy. 

%

\subsection{Test case 3 : Flooding of an initially dry floodplain}\label{numsectest3}
\begin{figure}[ht!] 
	\subfigure[Bottom elevation, $z_b(x,y)$ in $\Omega^2$]{\label{bedtest3}\includegraphics[width=0.58\textwidth]{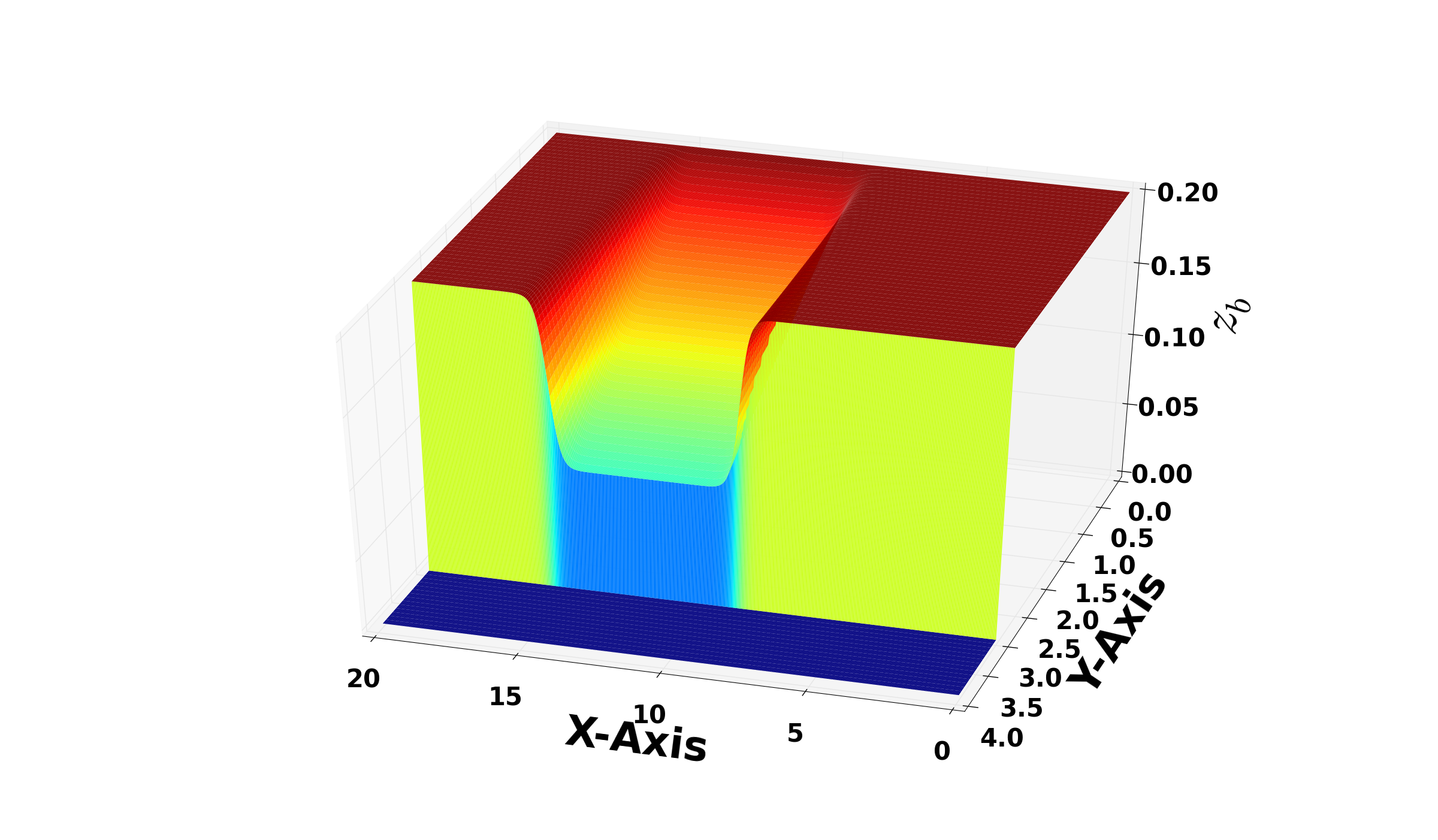}}
	\subfigure[Channel Wall Elevation, $\zwall$]{\label{zwalltest3}\includegraphics[width=0.4\textwidth]{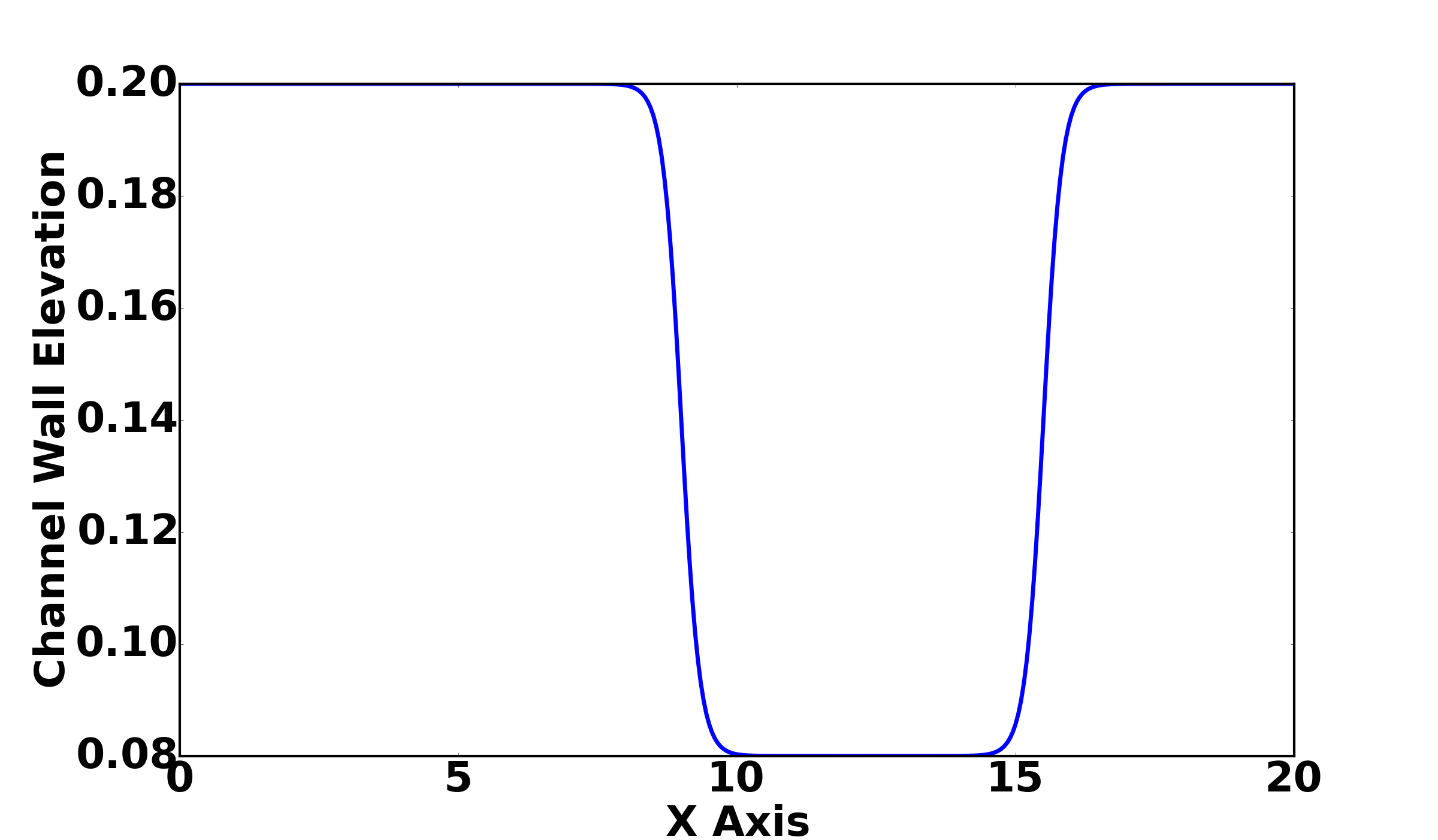}}
  \caption{The 2D bed elevation and channel wall elevation for test 3}
  \label{numfigtest3bedandchannelwall}
\end{figure} 
The final test case involves the overflowing of a channel onto an initially dry floodplain.
Both the channel and the floodplain are located in the $2D$ domain, $[0, 20]\times [0, 4]$.
The channel occupies the region, $[0,20]\times[y_c,4]$ with flat bottom, $Z_b(x)=0$, while the floodplain occupies the rest
of the domain, $[0,20]\times[0,y_c]$, where $y_c=3$. The bottom topography of the entire domain is the following
\begin{align}
  z_b(x,y) = \begin{cases} Z_b(x) = 0, & \mbox{ if } y \ge y_c, \\
                          0.2 + \frac{ \zwall - 0.2}{y_c}y , & \mbox{ otherwise },
             \end{cases},
\end{align}
where 
\begin{align}
  \zwall = \begin{cases}
                                -0.06 \tanh ( 3(x-9   ) ) + 0.14, & \mbox{ if } x \leq 10.5, \\
                                 0.06 \tanh ( 3(x-15.5) ) + 0.14,  & \mbox{ otherwise }
                   \end{cases}
\end{align}
is the elevation of the channel wall, see figures \ref{bedtest3} and \ref{zwalltest3} for the
plots of $z_b(x,y)$ and $\zwall$.

The initial condition consists of stationary water of depth, 0.08 meters in the channel and dry floodplain.
The boundary conditions are time-dependent water depth at the left boundary of the channel and zero
velocity at the right channel boundary, namely
\begin{align}
  &H(0,  y, t ) = \begin{cases} h_b(t), & \mbox{ if } t \leq 4a, \\
                                h_b(4a), & \mbox{ if } t>4a,
                  \end{cases}  \\
                  \mbox{ for } y \ge y_c.  \notag  \\ 
  &u(20, y, t ) = 0.0, \mbox{ for all } t \ge 0, y \ge y_c,
\end{align}
where
\begin{equation}
  h_b(t) = \eta_0 + r + r\sin \bigg( \frac{(t-a)\pi}{2a}  \bigg),
\end{equation}
where $a=10$ and $r=0.025$. $\eta_0=0.08$ is a constant initial free-surface elevation inside the channel.
The remaining boundaries are closed and the manning coefficients are the same as used in the previous cases.
The following probe points are chosen, $P_1=(2.5, 3.5)$, $P_2=(4.0,3.8)$, $P_3=(7.0,3.3)$, $P_4=(10.0,3.4)$,
$P_5=(11, 3.5)$, $P_6=(12,3.3)$, $P_7=(14,3.4)$, $P_8=(16,3.5)$, $P_9=(17.3,3.5)$, $P_{10}=(19,3.5)$, $P_{11}=(12,2.8)$, $P_{12}=(13,2.8)$,
$P_{13}=(12, 2.5)$, $P_{14}=(12, 2.0)$ and  $P_{15}=( 13.0, 1.0  )$.

\begin{figure}[ht!] 
	\includegraphics[width=0.98\textwidth]{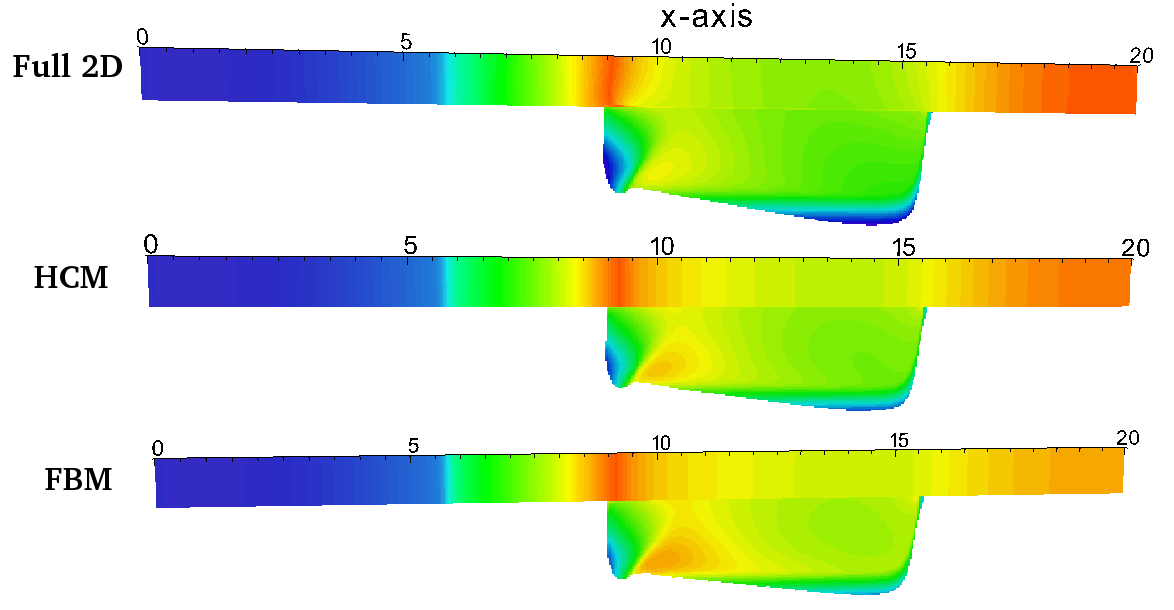}
  \caption{Visualisation of free surface elevation after $t=40$ for test case 3. The $x$-axis is from left to right,
  while the $y$-axis is from the bottom to the top.} 
  \label{numfigtest3etaview}
\end{figure}

\begin{figure}[ht!] 
	\includegraphics[width=0.98\textwidth]{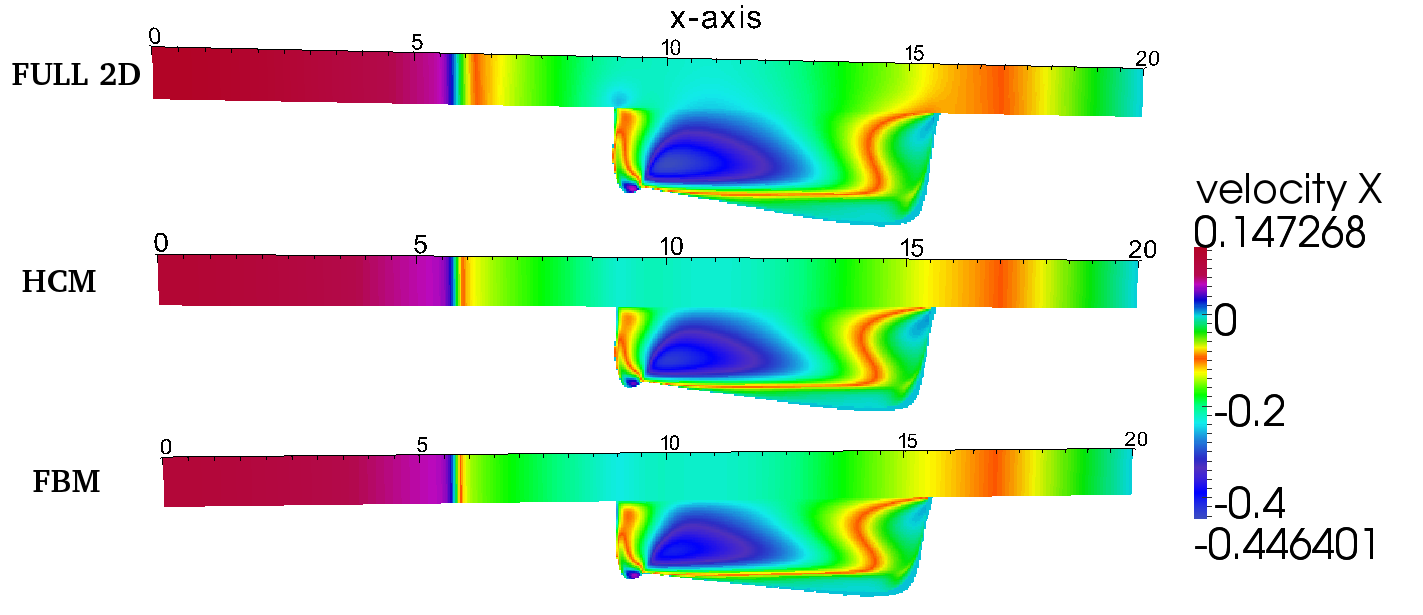}
  \caption{Visualisation of $x$-velocity after $t=40s$ for test case 3.
  The $x$-axis is from left to right,
  while the $y$-axis is from the bottom to the top.}
  \label{numfigtest3xveloview}
\end{figure}

\begin{figure}[ht!] 
	\includegraphics[width=0.98\textwidth]{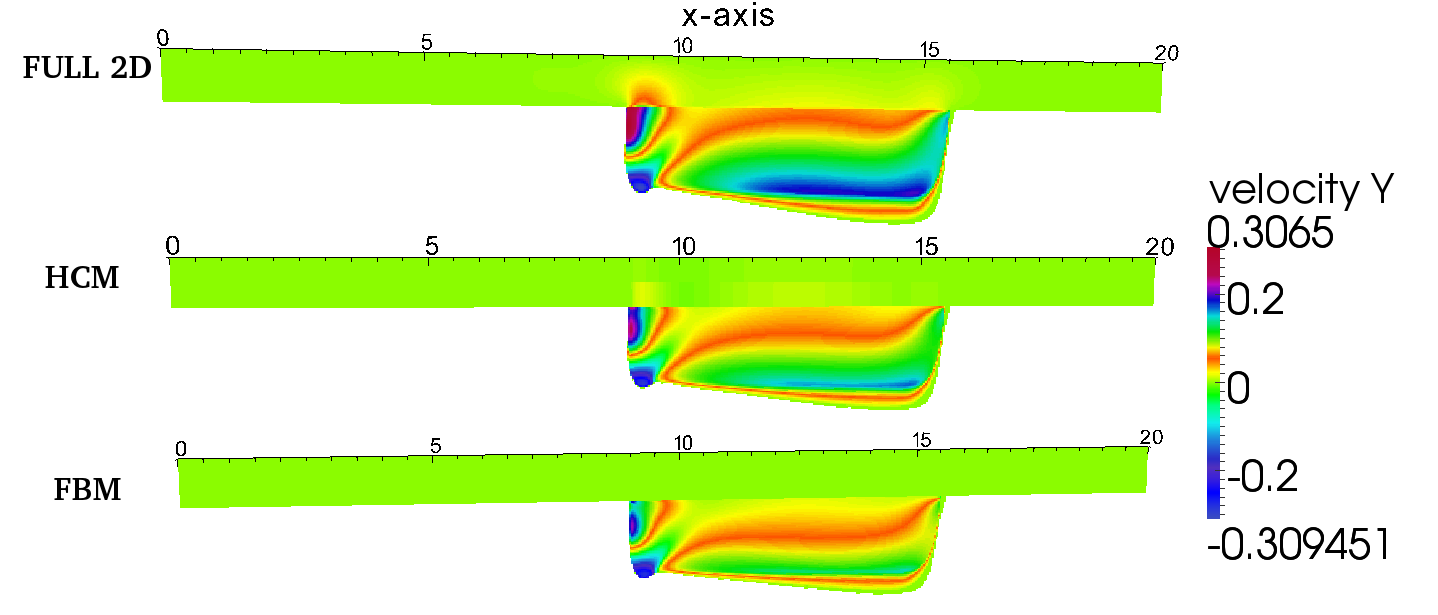}
  \caption{Visualisation of $y$-velocity after $t=40s$ for test case 3. The $x$-axis is from left to right,
  while the $y$-axis is from the bottom to the top.}
  \label{numfigtest3yveloview}
\end{figure}

\begin{figure}[ht] 
	\includegraphics[width=0.98\textwidth]{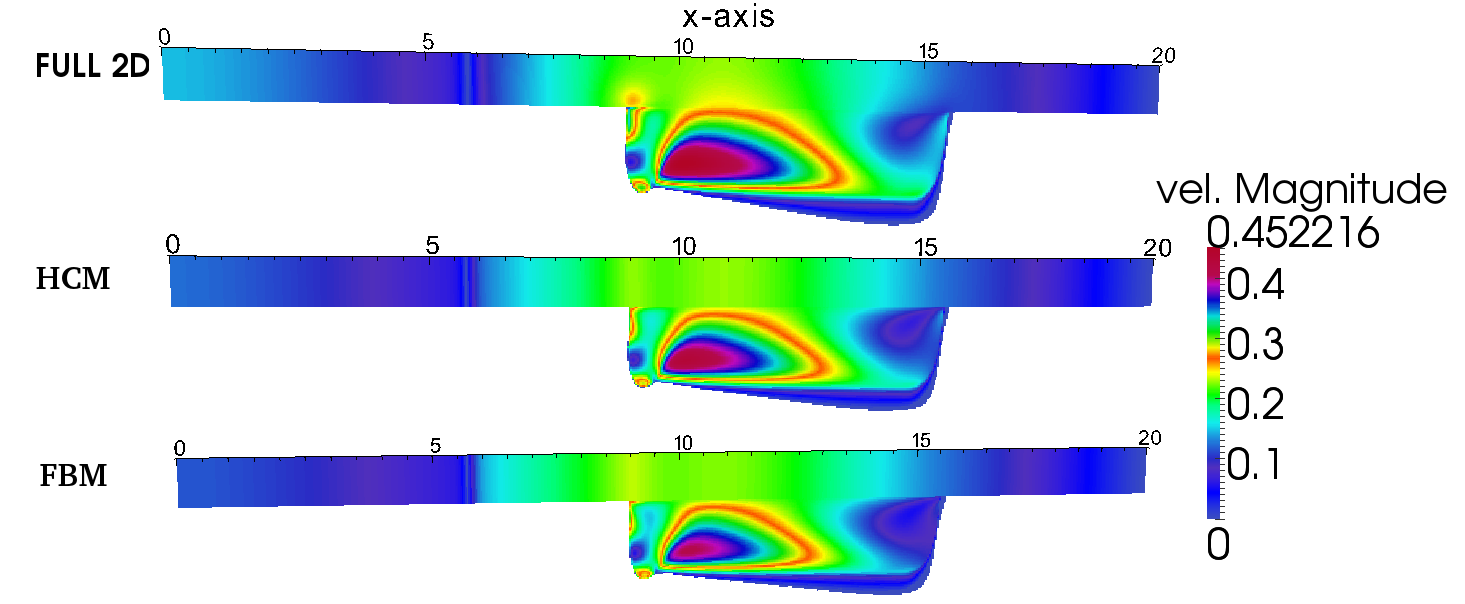}
  \caption{Visualisation of velocity magnitude after $t=40s$ for test case 3.
  The $x$-axis is from left to right,
  while the $y$-axis is from the bottom to the top.}
  \label{numfigtest3velomagview}
\end{figure}

\begin{figure}[ht!] 
 \subfigure{\includegraphics[width=0.33\textwidth, height=0.30\textwidth]{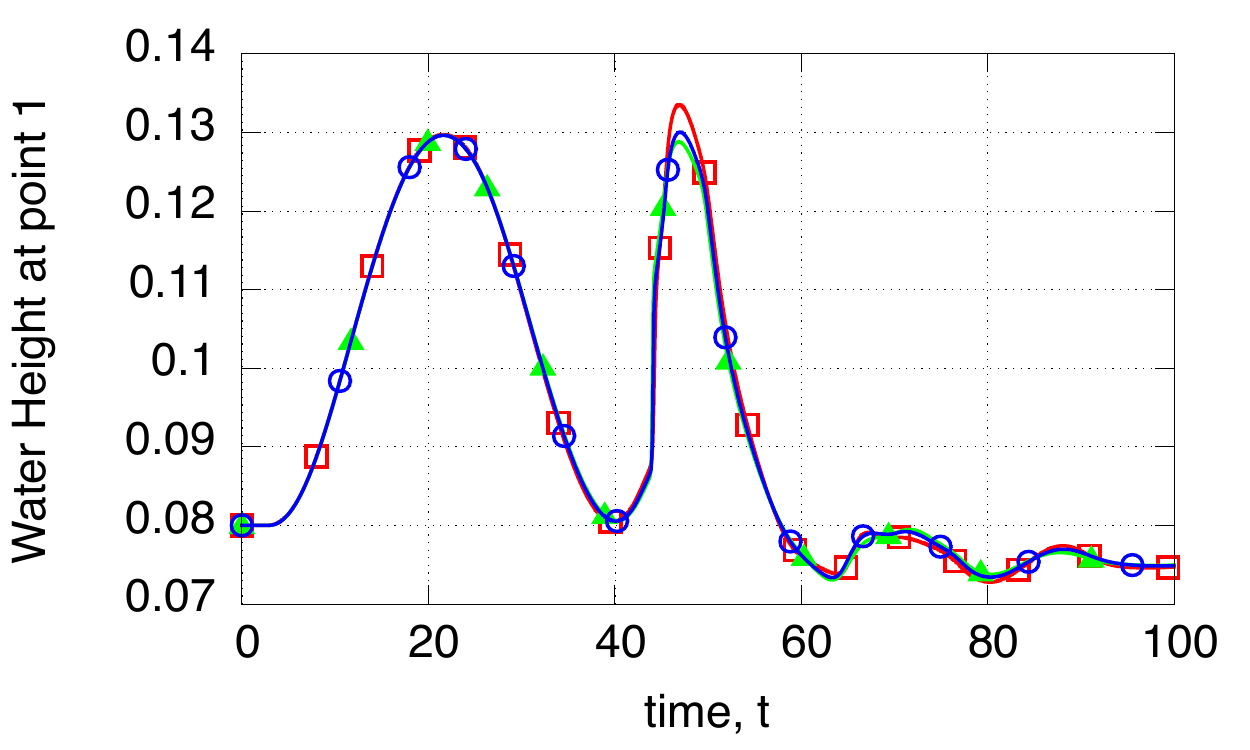}}
 \subfigure{\includegraphics[width=0.33\textwidth, height=0.30\textwidth]{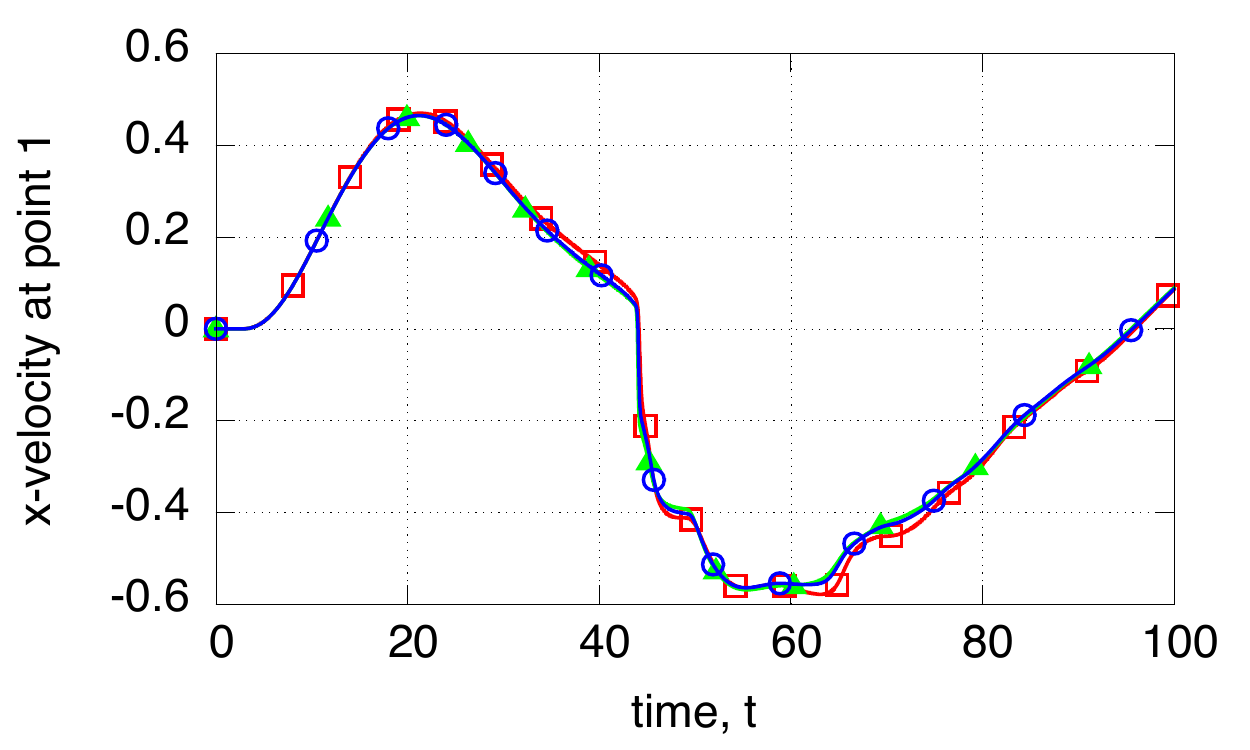}}
 \subfigure{\includegraphics[width=0.33\textwidth, height=0.30\textwidth]{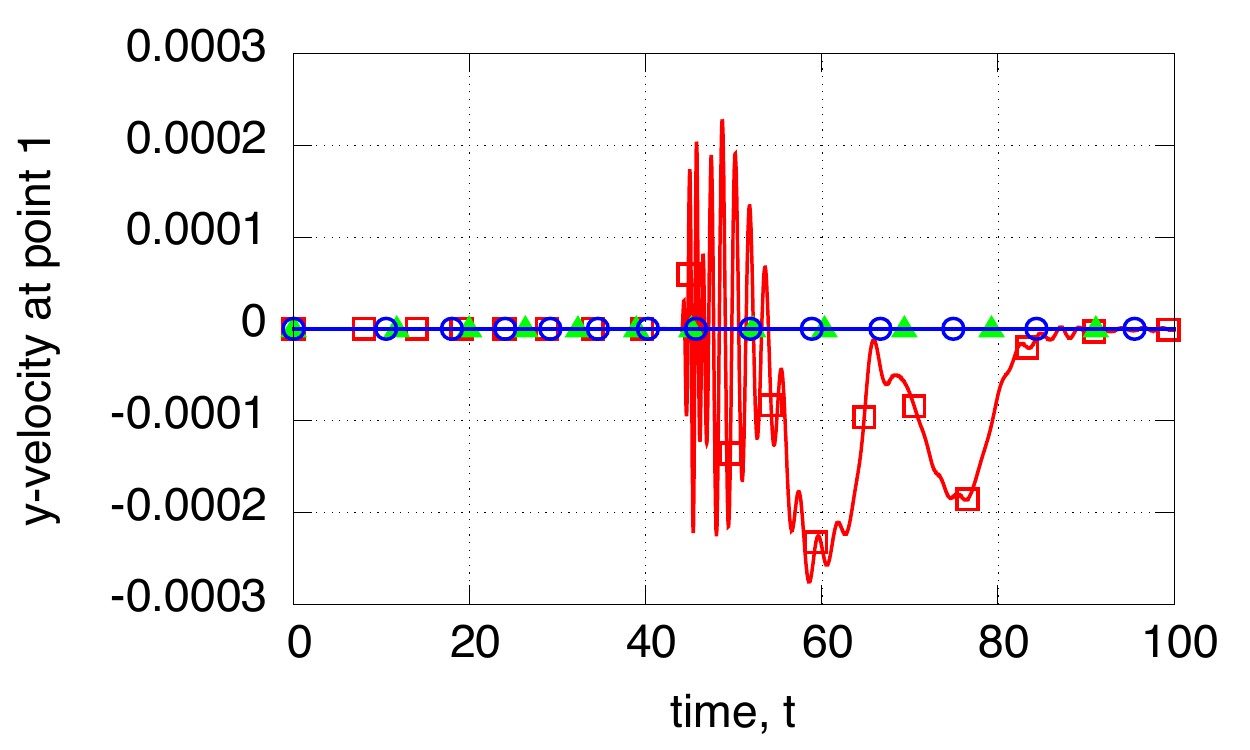}}\\

 \subfigure{\includegraphics[width=0.33\textwidth, height=0.30\textwidth]{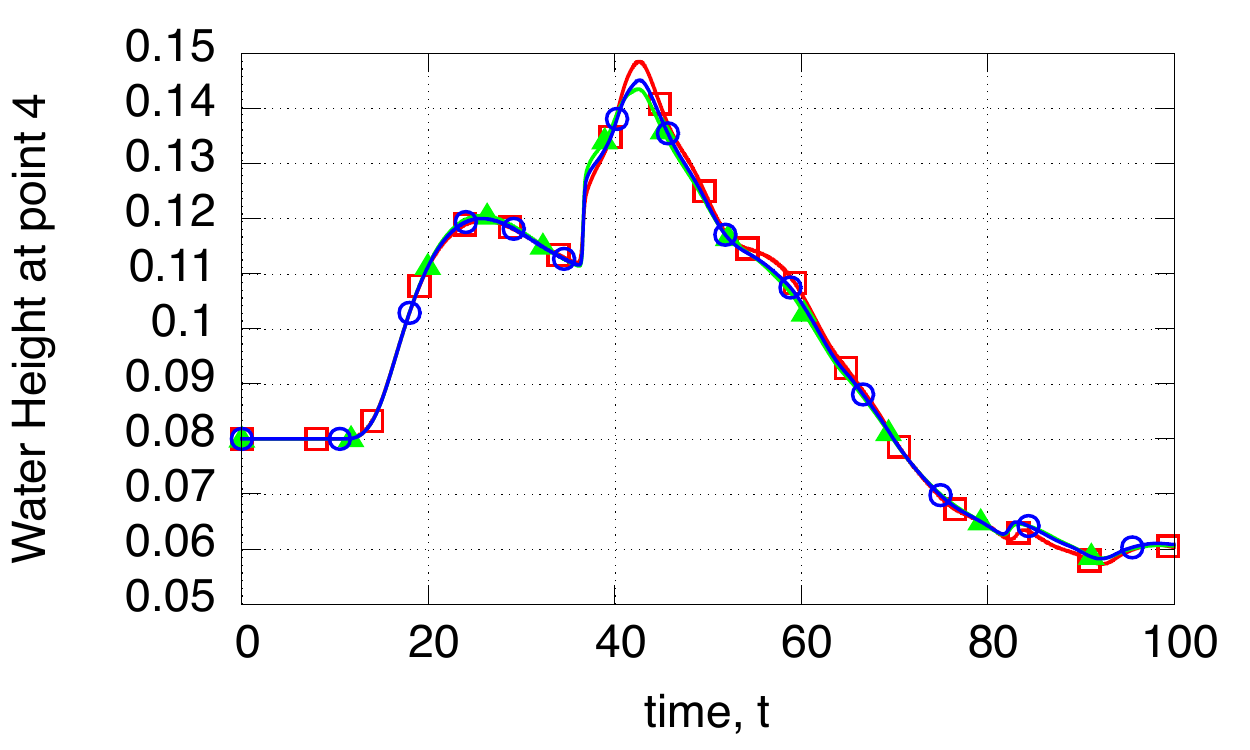}}
 \subfigure{\includegraphics[width=0.33\textwidth, height=0.30\textwidth]{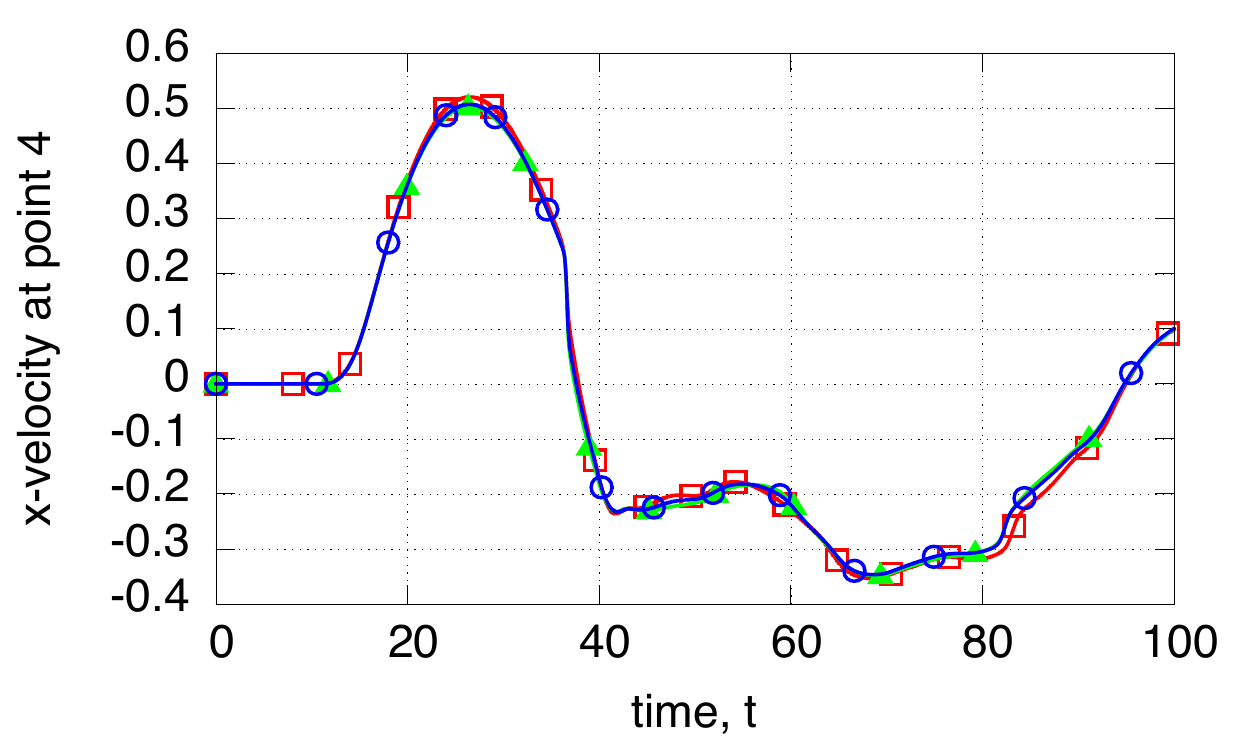}}
 \subfigure{\includegraphics[width=0.33\textwidth, height=0.30\textwidth]{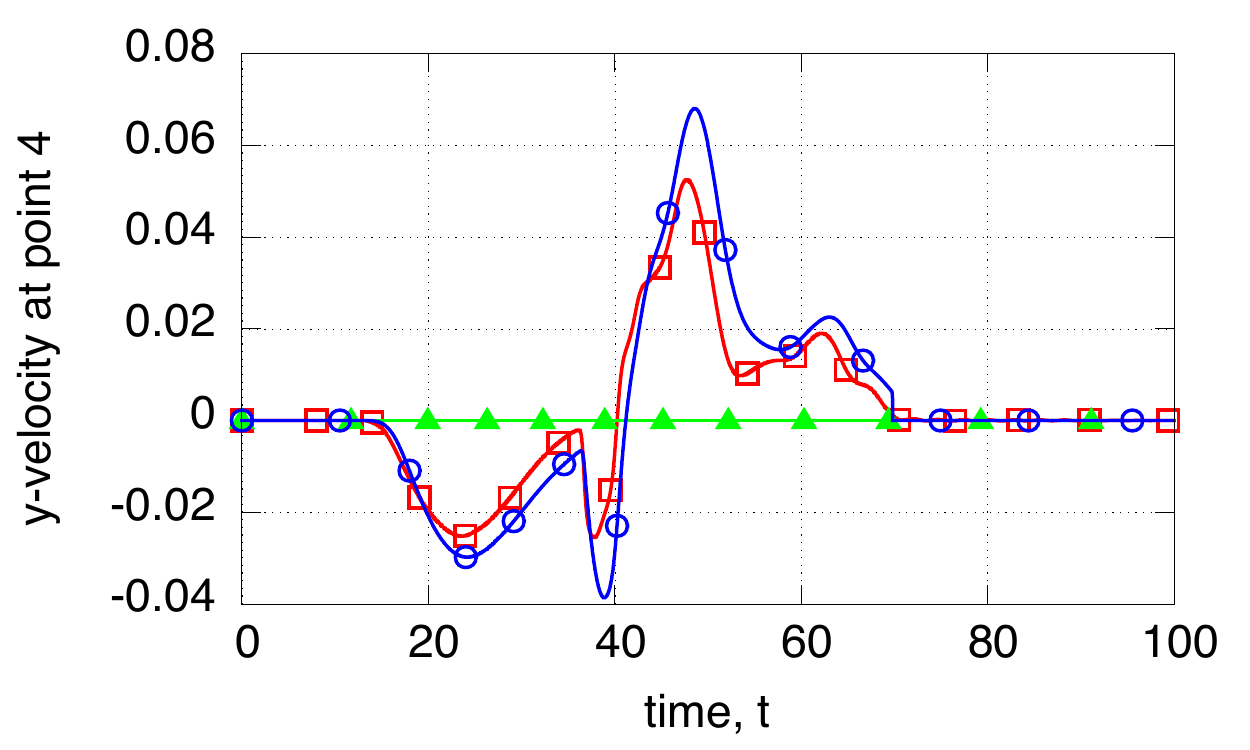}}\\

 \subfigure{\includegraphics[width=0.33\textwidth, height=0.30\textwidth]{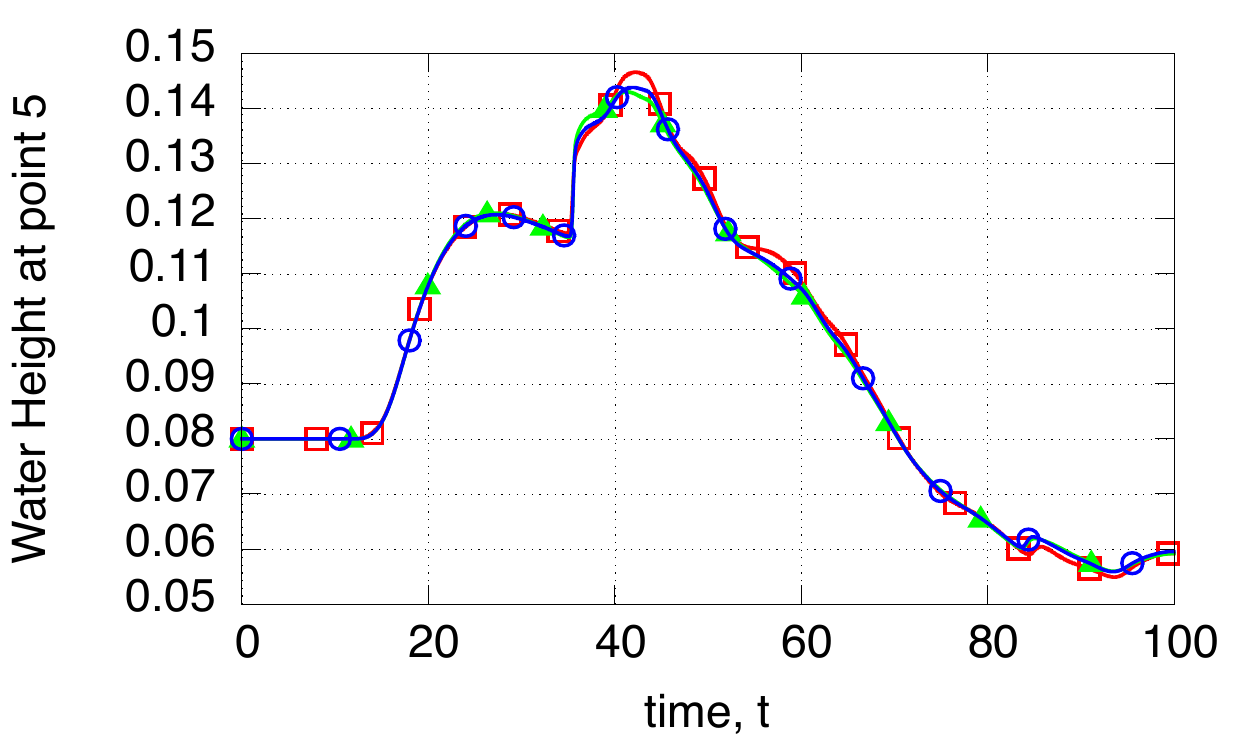}}
 \subfigure{\includegraphics[width=0.33\textwidth, height=0.30\textwidth]{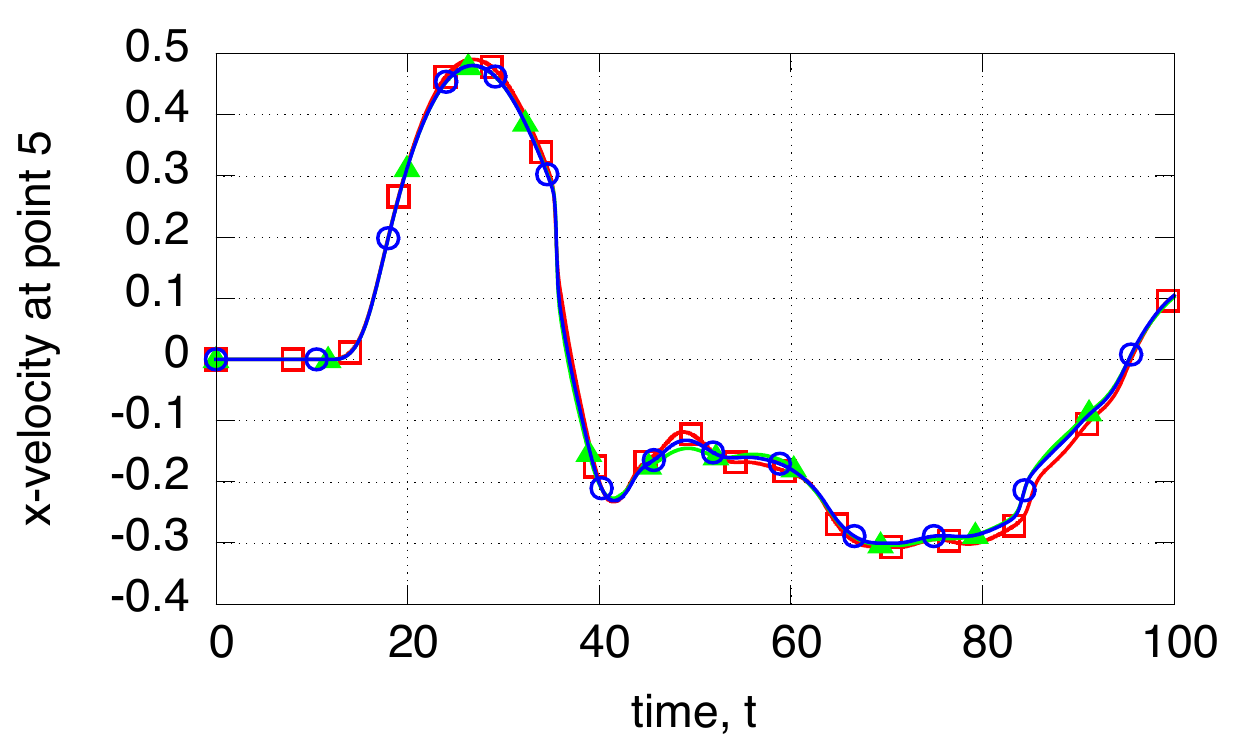}}
 \subfigure{\includegraphics[width=0.33\textwidth, height=0.30\textwidth]{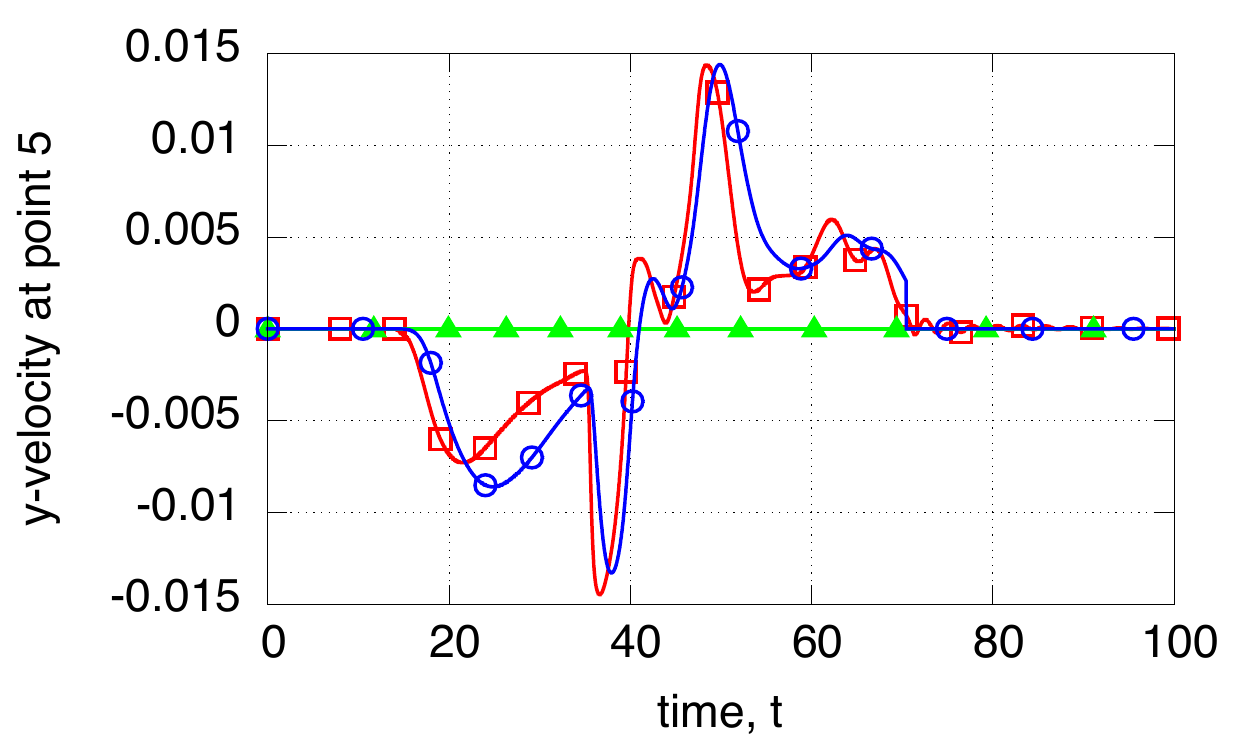}}\\

\subfigure{\includegraphics[width=0.33\textwidth,  height=0.30\textwidth]{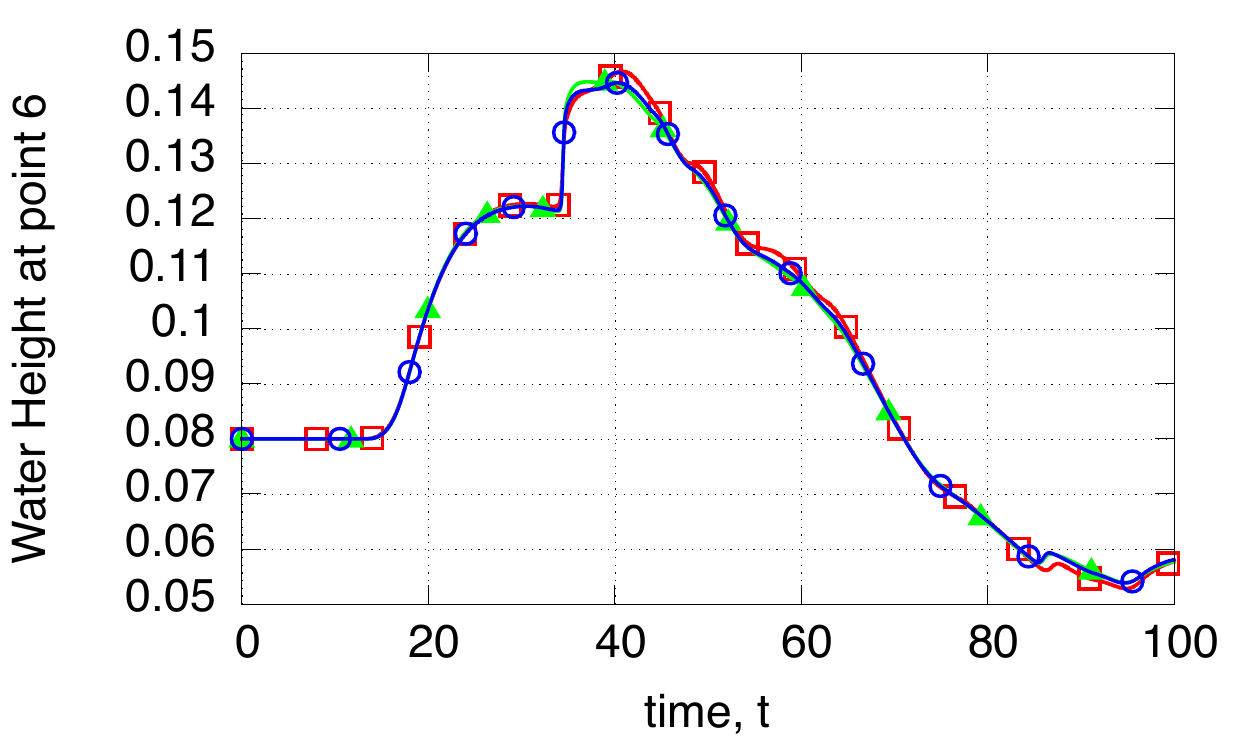}}
\subfigure{\includegraphics[width=0.33\textwidth,  height=0.30\textwidth]{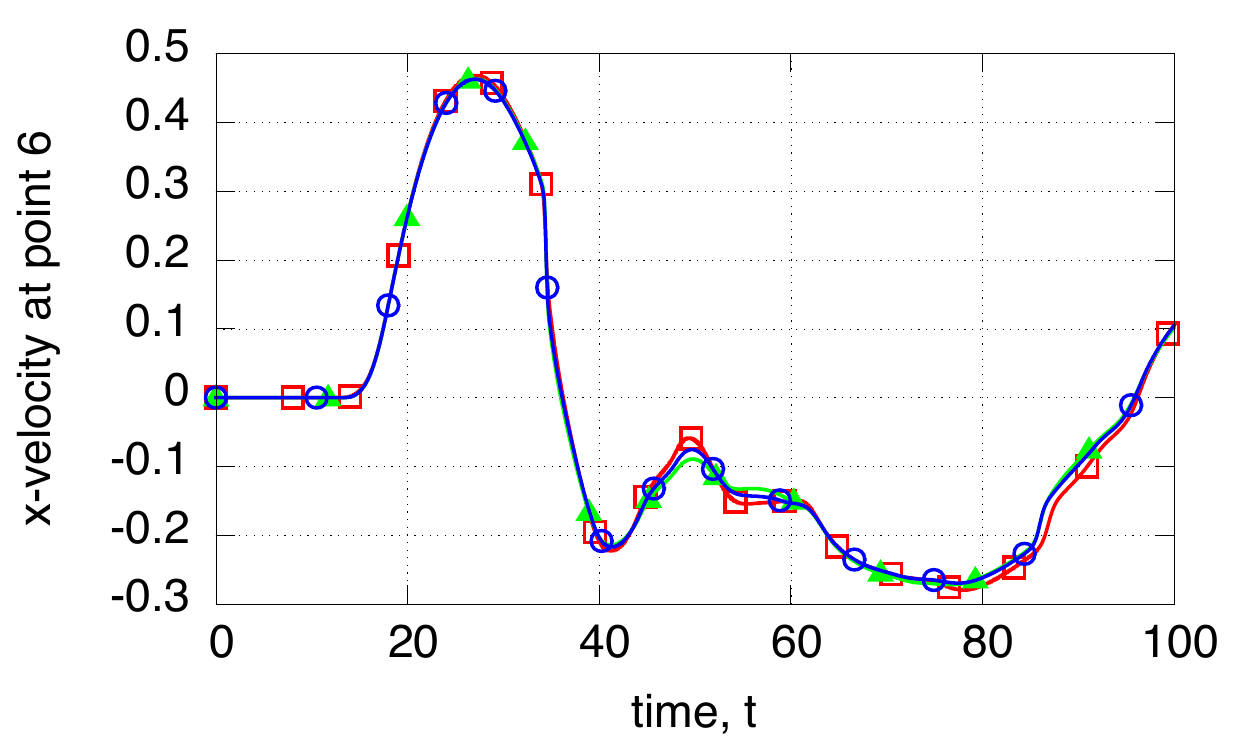}}
\subfigure{\includegraphics[width=0.33\textwidth,  height=0.30\textwidth]{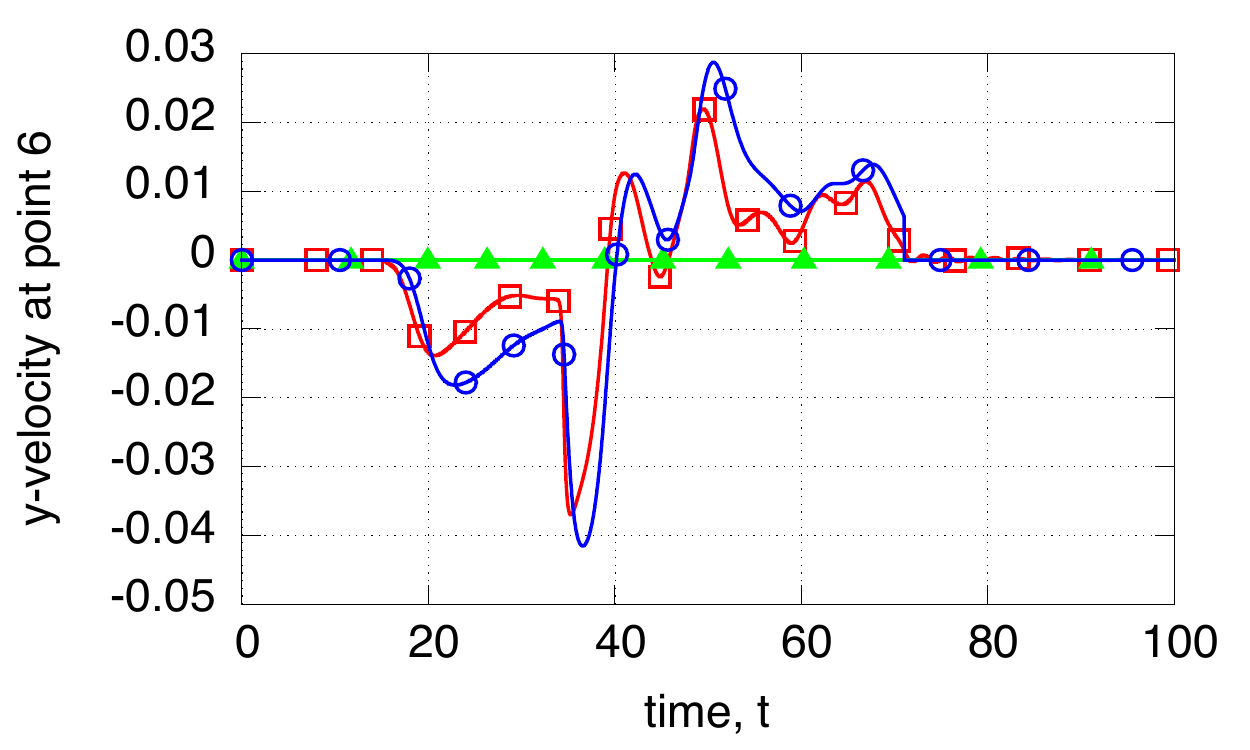}}\\
\subfigure{\includegraphics[width=0.5\textwidth]{legend-use}}
  \caption{Time variation of water depth $H$ (left column), $x$-velocity component (middle column) and
           $y$-velocity component (right column) at the indicated probe points within the channel for test case 3.
          Each row corresponds to one probe point.}

  \label{numfigtest3channelprobepoints}
\end{figure}

\begin{figure}[ht!] 
  \subfigure{\includegraphics[width=0.33\textwidth, height=0.30\textwidth]{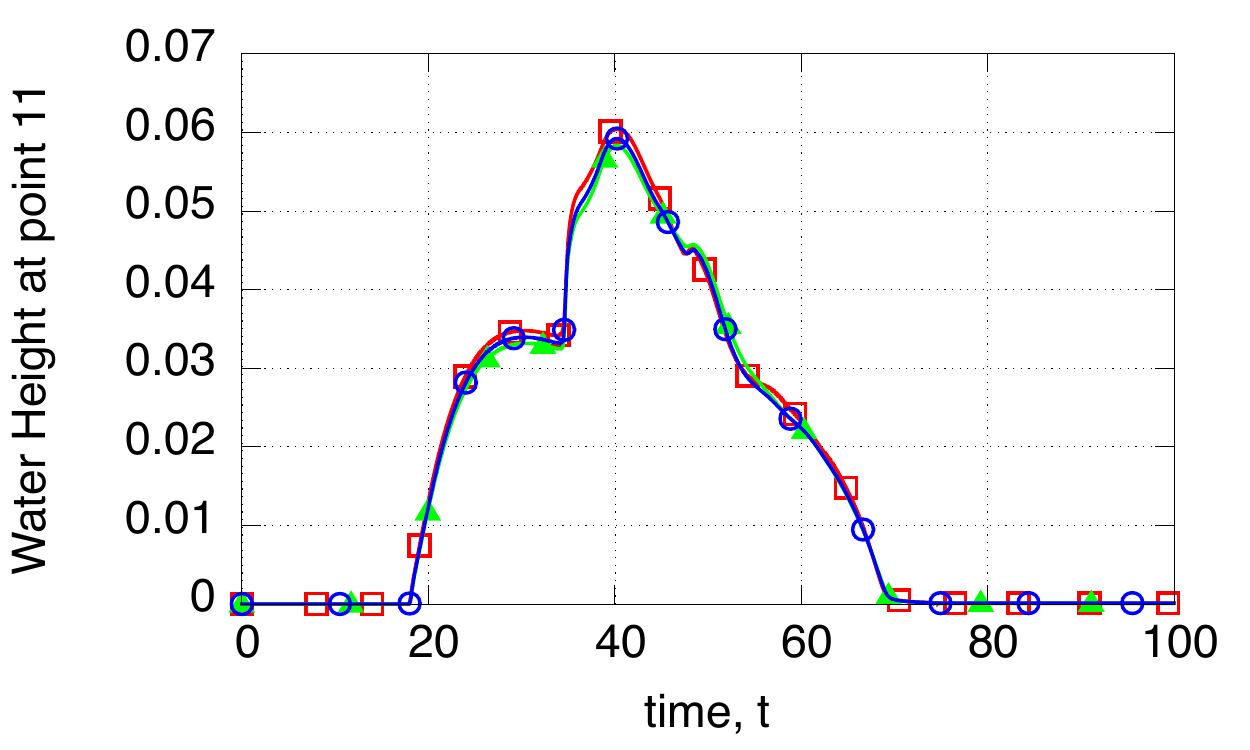}}
 \subfigure{\includegraphics[width=0.33\textwidth, height=0.30\textwidth ]{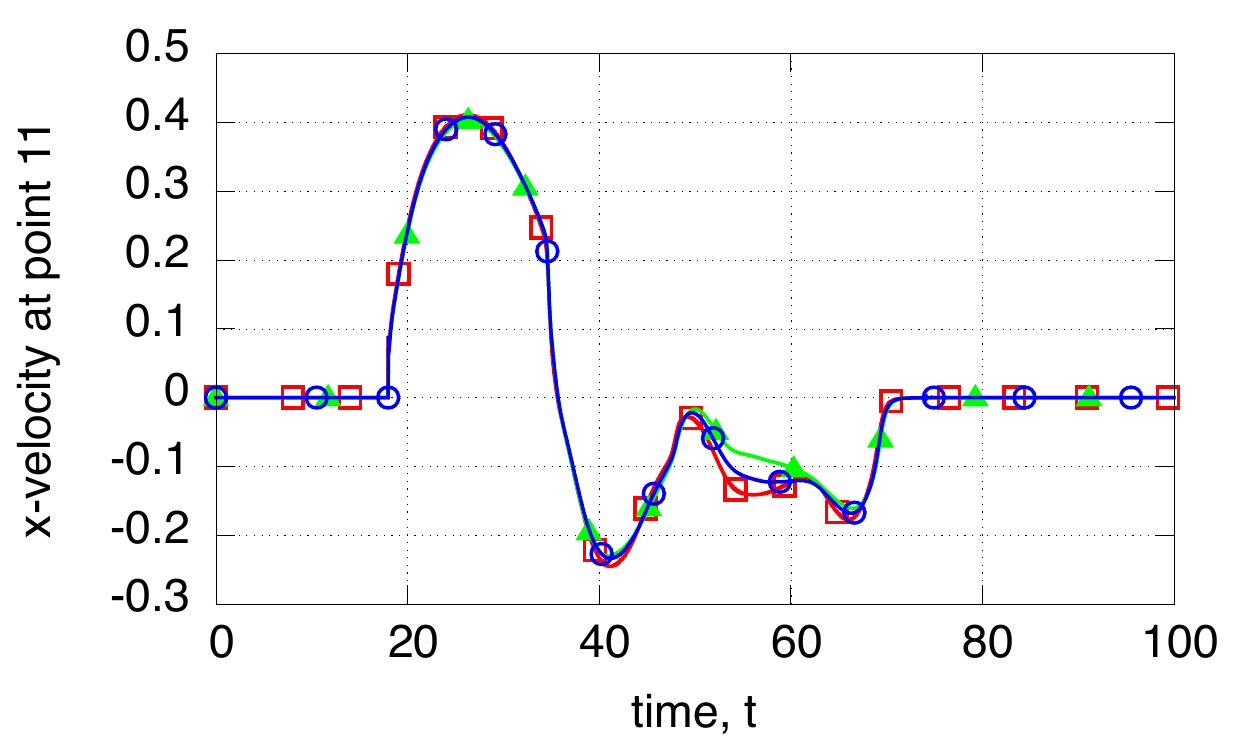}}
 \subfigure{\includegraphics[width=0.33\textwidth, height=0.30\textwidth]{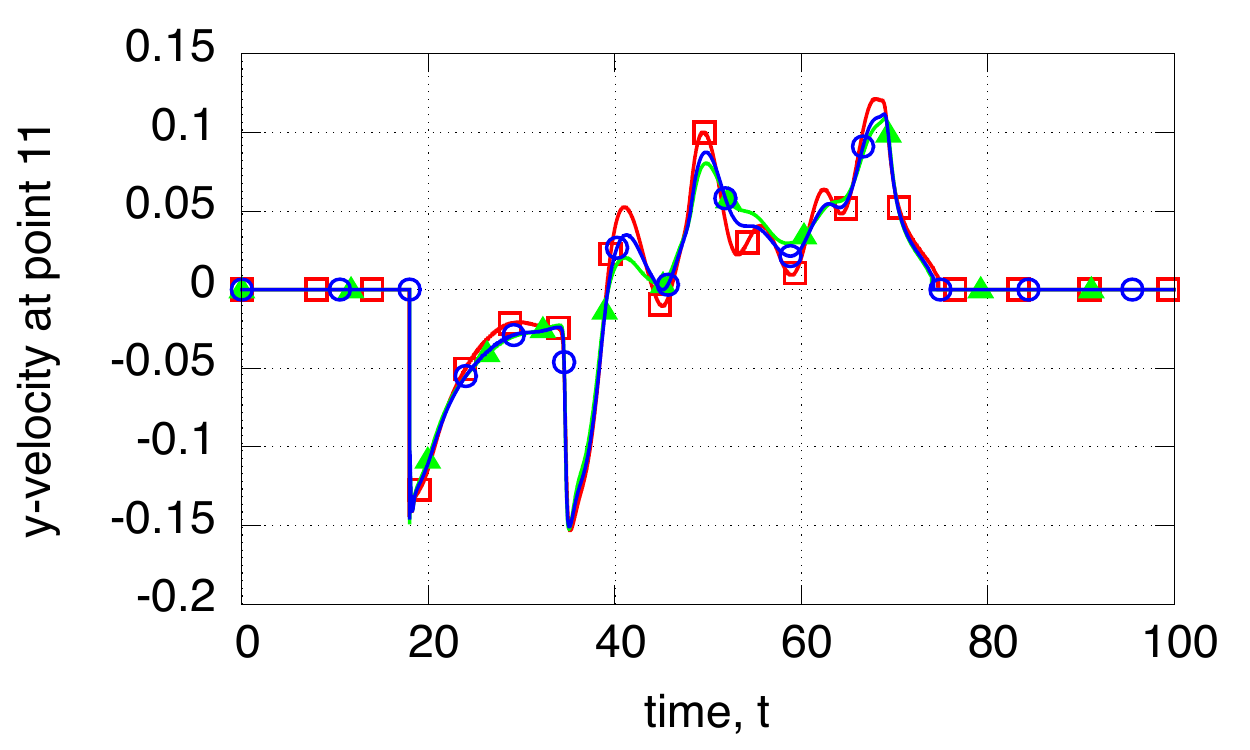}}\\

 \subfigure{\includegraphics[width=0.33\textwidth, height=0.30\textwidth]{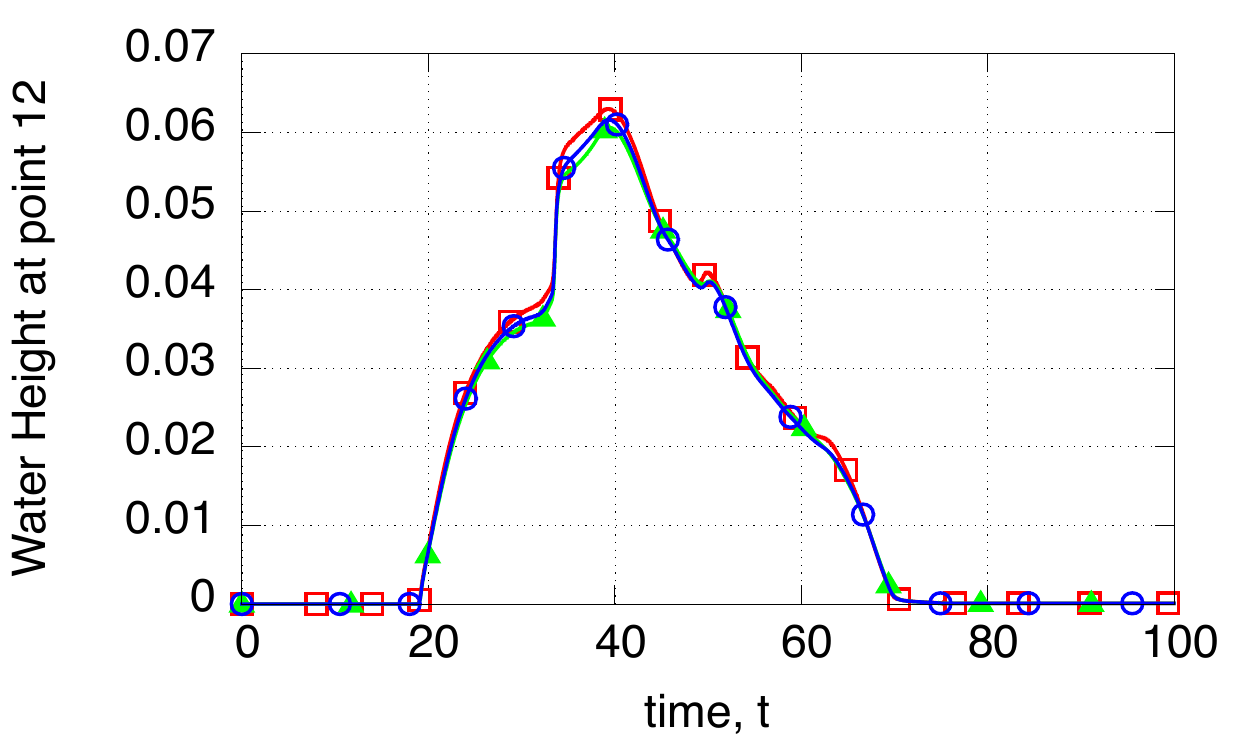}}
 \subfigure{\includegraphics[width=0.33\textwidth, height=0.30\textwidth ]{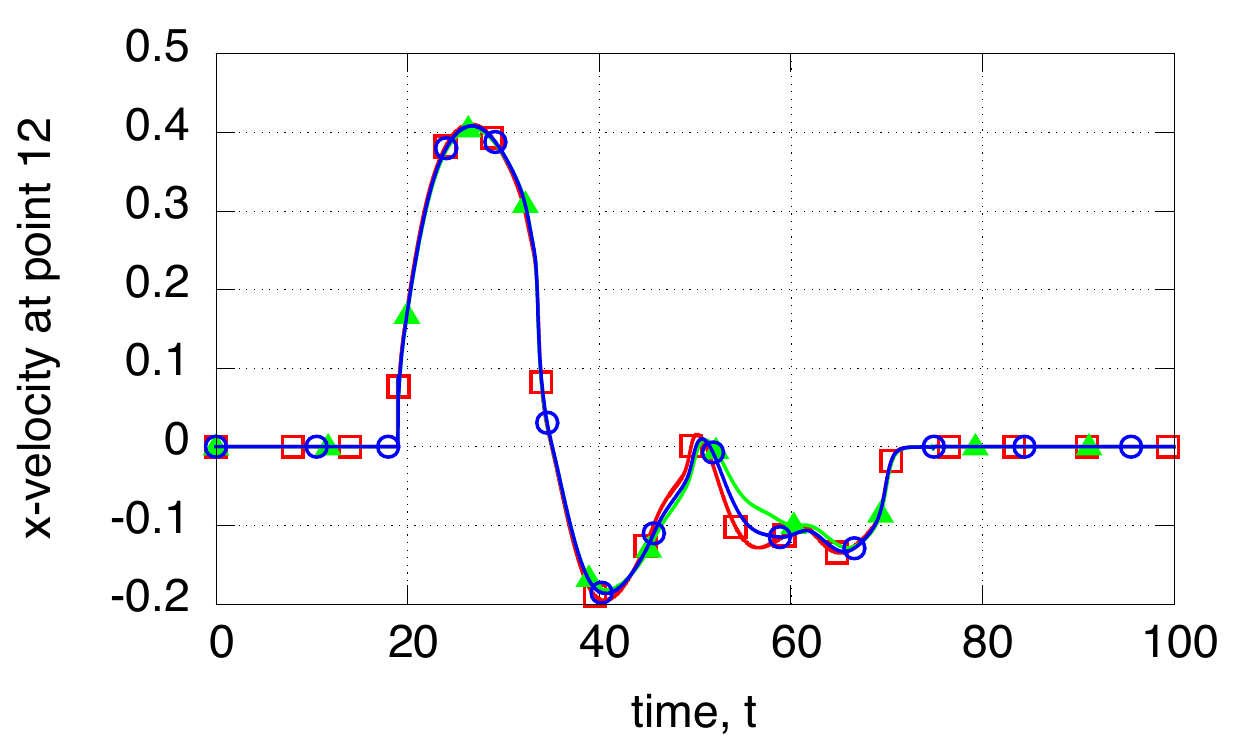}}
 \subfigure{\includegraphics[width=0.33\textwidth, height=0.30\textwidth]{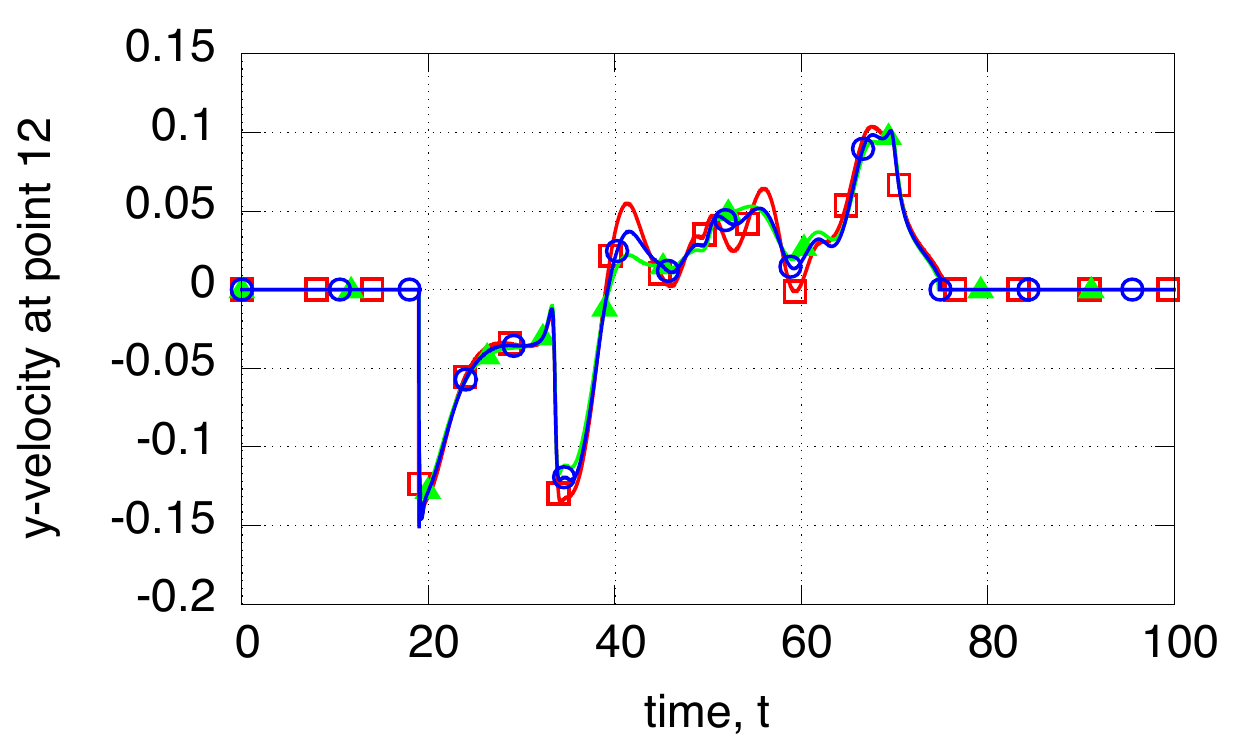}}\\

 \subfigure{\includegraphics[width=0.33\textwidth, height=0.30\textwidth]{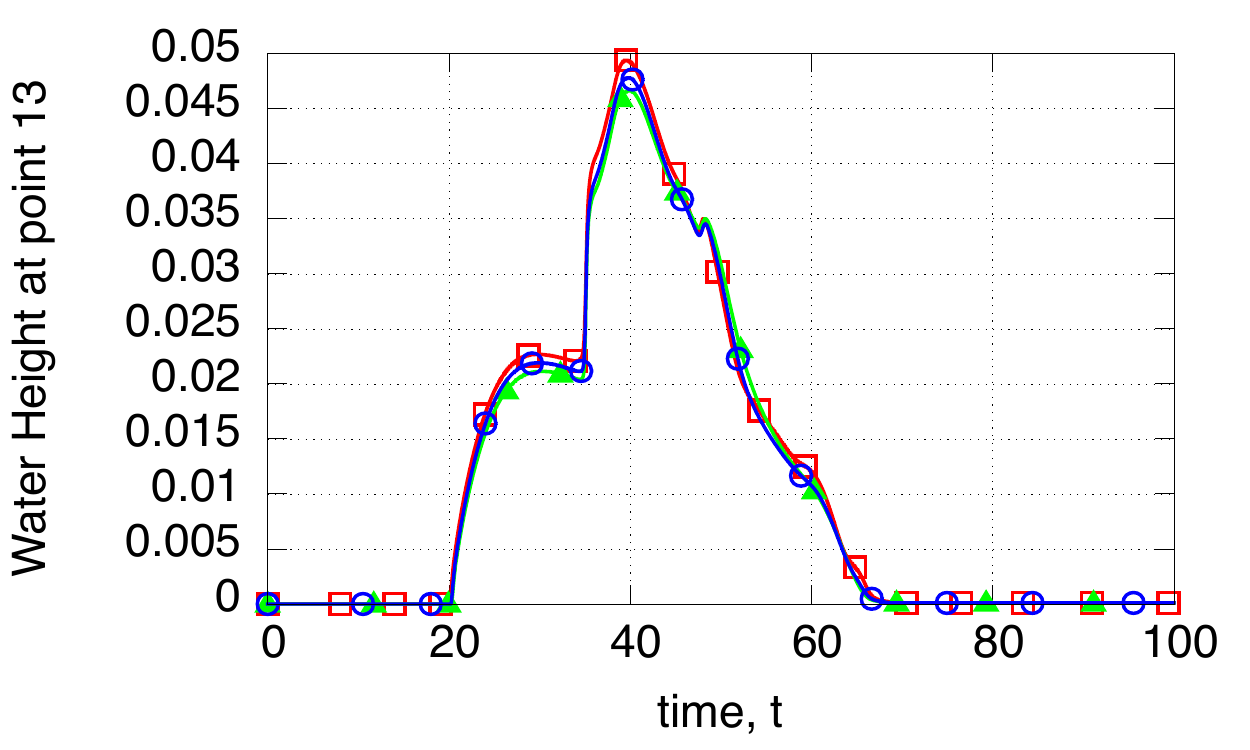}}
 \subfigure{\includegraphics[width=0.33\textwidth, height=0.30\textwidth]{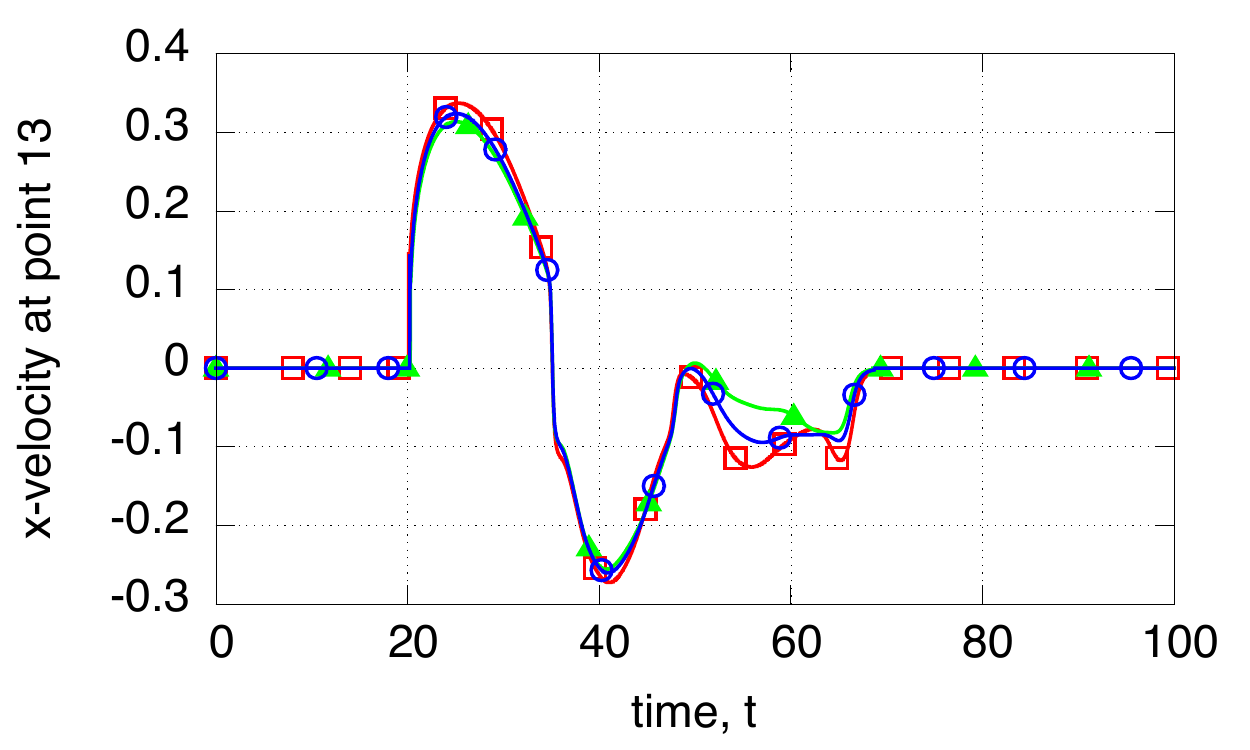}}
 \subfigure{\includegraphics[width=0.33\textwidth, height=0.30\textwidth]{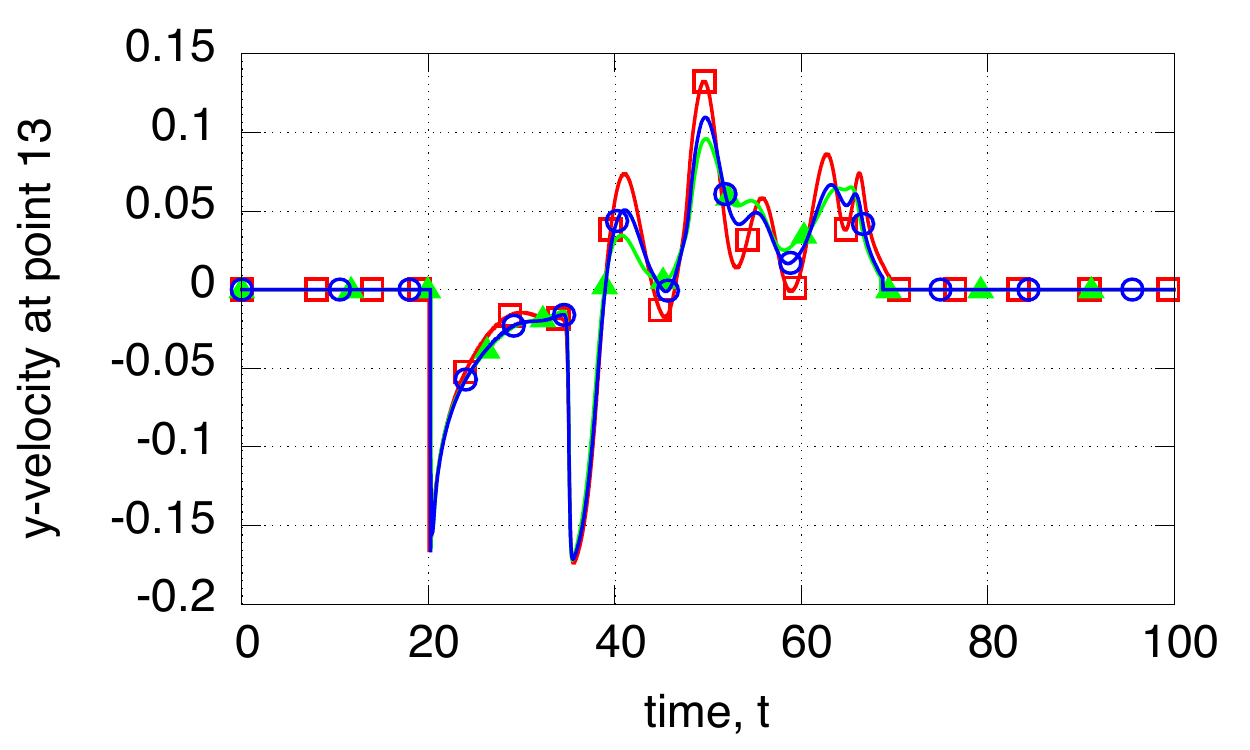}}\\

\subfigure{\includegraphics[width=0.33\textwidth,  height=0.30\textwidth]{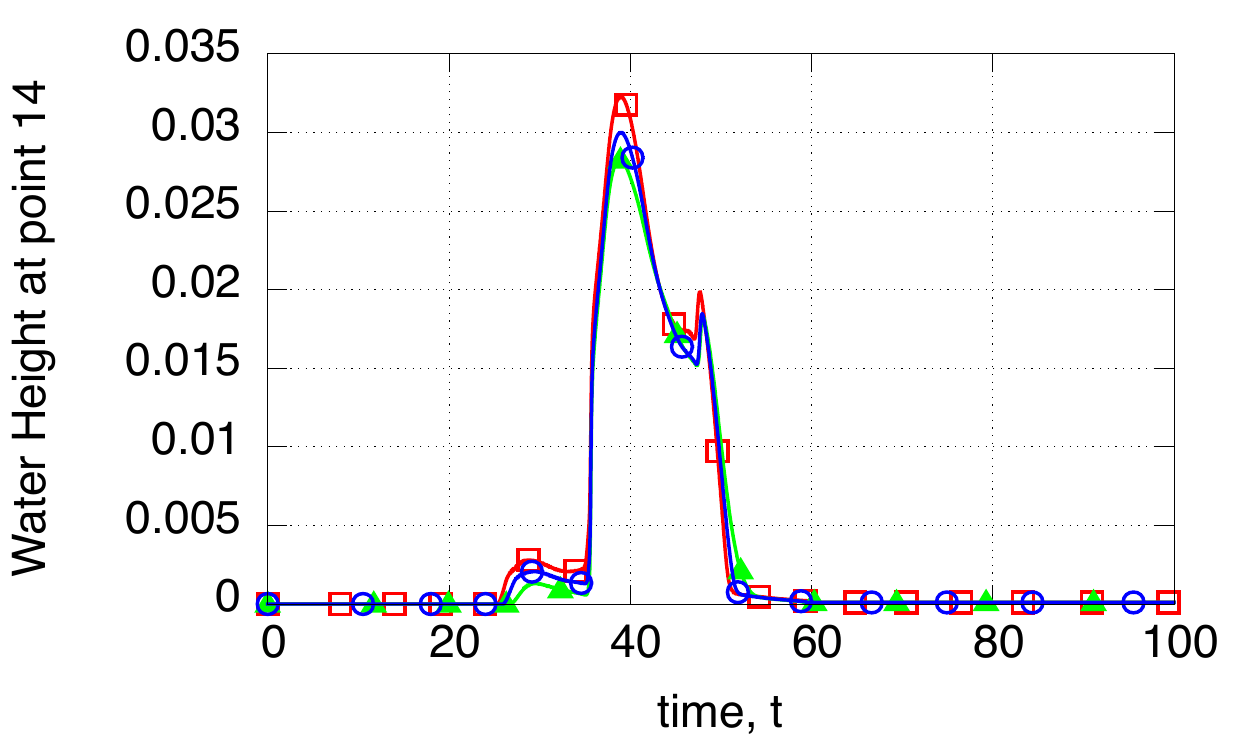}}
\subfigure{\includegraphics[width=0.33\textwidth,  height=0.30\textwidth]{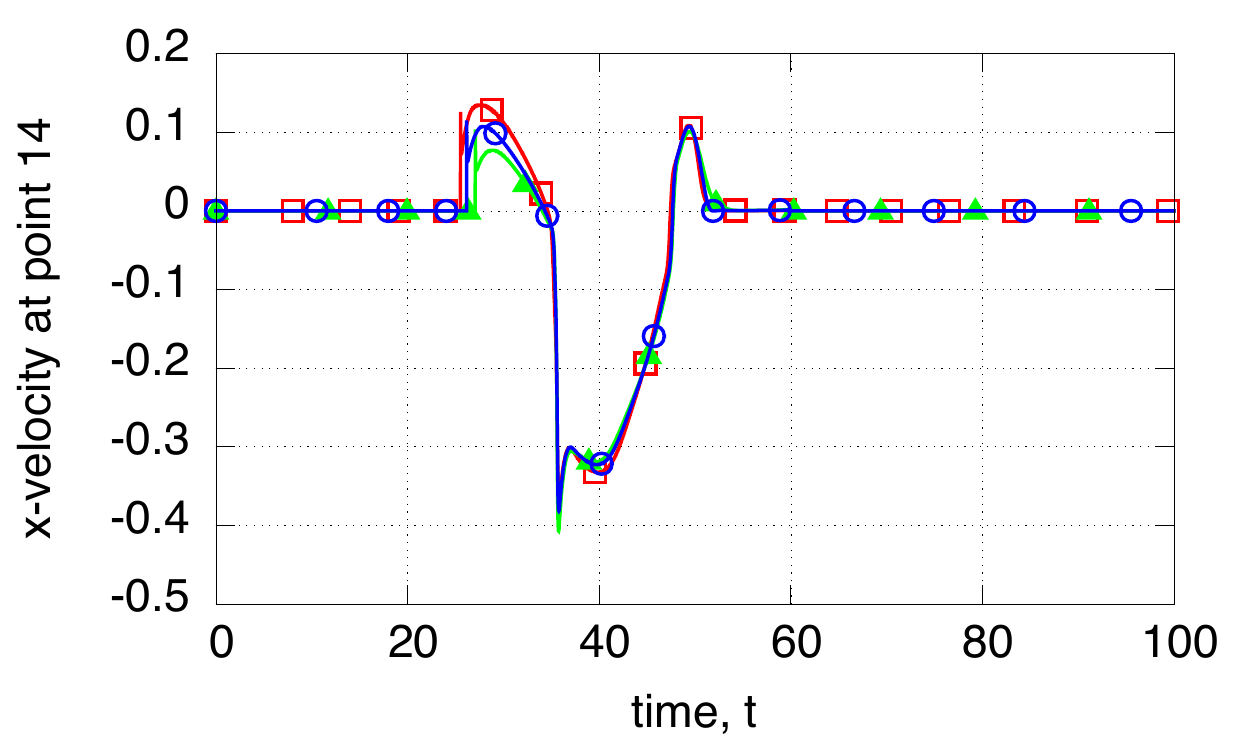}}
\subfigure{\includegraphics[width=0.33\textwidth,  height=0.30\textwidth]{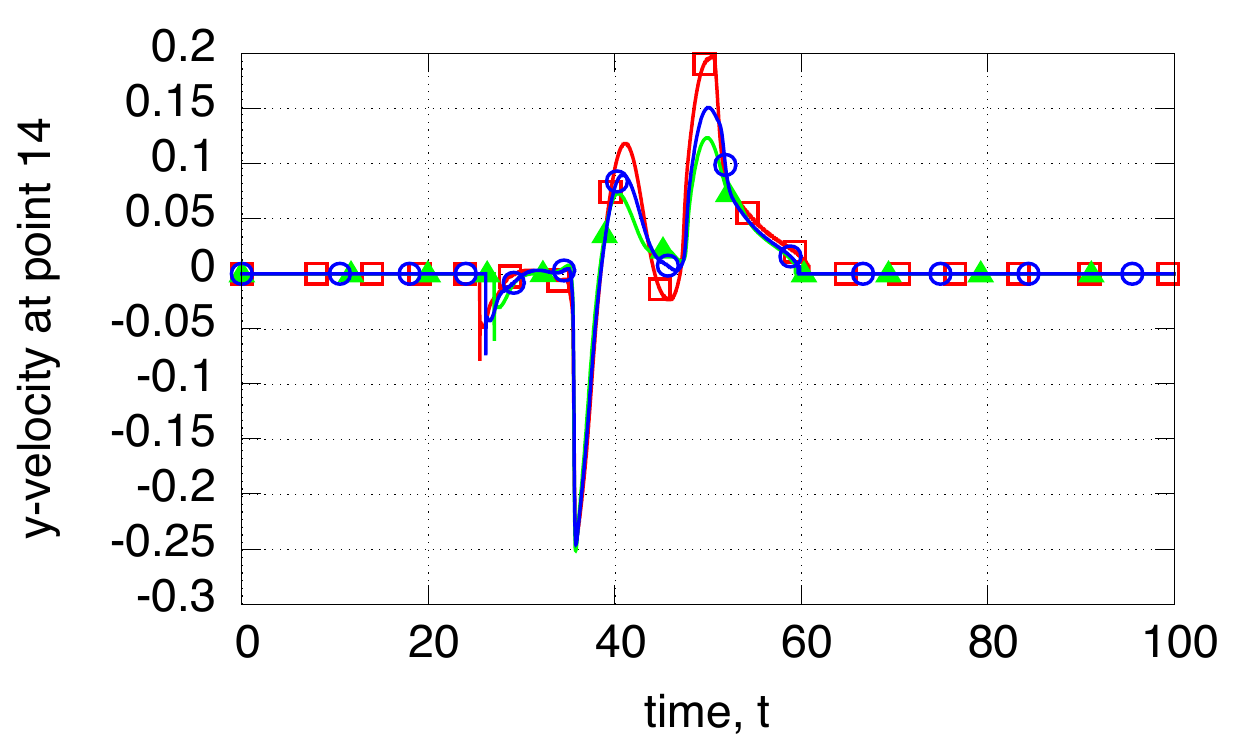}}\\

  \subfigure{\includegraphics[width=0.5\textwidth]{legend-use}}

  \caption{Time variation of water depth $H$ (left column), $x$-velocity component (middle column) and
           $y$-velocity component (right column) at the indicated probe points in the floodplain for test case 3.
          Each row corresponds to one probe point.}
  \label{numfigtest3floodplainprobepoints}
\end{figure}

Table \ref{numtabletest3} shows the domain discritization for both the channel and the floodplain for each method being discussed.
As before, all methods use the same grid for the floodplain but different grids for the channel. This problem was simulated
for $t=100$ seconds. We report, in figures \ref{numfigtest3etaview} - \ref{numfigtest3velomagview}, the results of the
simulation after 40 seconds and in figures \ref{numfigtest3channelprobepoints} and \ref{numfigtest3floodplainprobepoints}, we report
the results at selected probe points throughout the duration of the simulation.

As can been seen from the pictures, both coupling methods provide very good approximation of full 2D simulation results for both the 
free surface elevation (figure \ref{numfigtest3etaview}), the velocity components (figures \ref{numfigtest3xveloview} and \ref{numfigtest3yveloview})
and the velocity magnitude (figure \ref{numfigtest3velomagview}) for this test case.
And in terms of accuracy of y-velocity component,
the HCM provides better approximations as can be seen in figure \ref{numfigtest3yveloview}.

\begin{table}[h] 
\begin{center}
\begin{tabular}{ | l |  c | c | }
\hline
                 &  Channel Grid   &  Floodplain Grid    \\
\hline
Full 2D          &  $600\times30$  &  $600\times90$      \\ 
\hline
HCM              &  $600\times2$   &  $600\times90$      \\ 
\hline   
FBM              &  $600\times1$   &  $600\times90$     \\ 
\hline 
\end{tabular}
\end{center}
\caption{Grid cells, simulation times and number of time steps : Test 3}
\label{numtabletest3}
\end{table}

To further understand the results of the simulations, the time evolution of the flow  quantities at the
probe points, $P_1 - P_{15}$ have been examined. Here we report the results at the probe points
$P_1, P_4, P_5 $ and $P_6$ which are in the channel and the points, $P_{11}, P_{12}, P_{13}$
and $P_{14}$ in the floodplain. Figures \ref{numfigtest3channelprobepoints}  and \ref{numfigtest3floodplainprobepoints} 
show the results for the selected points in the channel and floodplain respectively.
In each figure, the left column displays the water depth, the second (middle) column shows the $x-$component 
of velocity, while the third(right) column shows
the $y-$component of velocity.

From figure \ref{numfigtest3channelprobepoints}, we can see that both coupling methods provide very
good approximation of the results of the full 2D simulations, especially for the water depth and
$x-$component of velocity. However, only the HCM is able to compute the variation in the
$y-$velocity and it does so with very good accuracy, see $P_4,P_5,P_6$ in figure \ref{numfigtest3channelprobepoints}.
This further verifies the ability of of the HCM to compute the lateral discharges within the channel.

From figure \ref{numfigtest3floodplainprobepoints}, we also see that for the points in the floodplain,
the coupling methods computed
very good approximations of results of the full 2D simulation with the HCM computing more accurate
results especially for the $y-$velocity. This figure also verify the no-numerical flooding property
of the methods. That is, the floodplain initially remained dry until the time when water height
rose above the channel banks. This is the reason why, for all points in the floodplain, the water depth
and velocity remained at zero for the first several seconds of the simulation. Another thing to note is 
that due to the time-dependent boundary condition for this problem,  water flowed onto the floodplain 
and after some time the water level in the channel decreased, hence the water in the floodplain drains back into the channel
leaving the floodplain dry again. The coupling methods truly capture this phenomenon as one can see in figure
\ref{numfigtest3floodplainprobepoints} where the water depth and velocity return to zero towards
the end of the simulation and remain at zero throughout the rest of the simulation.
This is true for all the points in the floodplain, even those not reported here.

\section{Conclusion}\label{hsecconc}
A horizontal coupling method has been proposed, implemented and tested in this paper.
It presents a strategy to overcome the difficulty in computing the channel lateral discharges,
circumvent the 1D assumption on the channel lateral discharge during flooding and propose a variant
of the coupling term of \cite{monnierMarin2009super} without the use of or imposing any 
restriction on the channel width variation. Numerical experiments show that the method computes
adequate results. 
Particularly, the channel lateral 
discharges are properly computed without adopting complicated/iterative procedures.
Finally, we note that  for all the numerical test cases considered in this paper,
the HCM would coincide with the FBM if the lateral discharges were not computed in the HCM. Therefore, the
improved solution observed in the HCM over the FBM, for these test cases, is a result of the lateral
discharges that are computed in the HCM. We, therefore,  conclude that properly computing and 
restoring the channel lateral discharge, improves the quality of the computed solution and
this can be done without introducing much computational overhead.

\section*{Acknowledgement}
We are grateful to the Petroleum Technology Development Fund (PTDF), Nigeria for funding
this study and to the Centre for Scientific Computing, University of Warwick for providing
the computing resources.

\bibliographystyle{abbrv}
\bibliography{mybib}            

\end{document}